# NEARLY UNBIASED VARIABLE SELECTION UNDER MINIMAX CONCAVE PENALTY

### By Cun-Hui Zhang[1]

*Rutgers University*

We propose MC+, a fast, continuous, nearly unbiased and accurate method of penalized variable selection in high-dimensional linear regression. The LASSO is fast and continuous, but biased. The bias of the LASSO may prevent consistent variable selection. Subset selection is unbiased but computationally costly. The MC+ has two elements: a minimax concave penalty (MCP) and a penalized linear unbiased selection (PLUS) algorithm. The MCP provides the convexity of the penalized loss in sparse regions to the greatest extent given certain thresholds for variable selection and unbiasedness. The PLUS computes multiple exact local minimizers of a possibly nonconvex penalized loss function in a certain main branch of the graph of critical points of the penalized loss. Its output is a continuous piecewise linear path encompassing from the origin for infinite penalty to a least squares solution for zero penalty. We prove that at a universal penalty level, the MC+ has high probability of matching the signs of the unknowns, and thus correct selection, without assuming the strong irrepresentable condition required by the LASSO. This selection consistency applies to the case of $p \gg n$, and is proved to hold for exactly the MC+ solution among possibly many local minimizers. We prove that the MC+ attains certain minimax convergence rates in probability for the estimation of regression coefficients in $\ell_r$ balls. We use the SURE method to derive degrees of freedom and $C_p$-type risk estimates for general penalized LSE, including the LASSO and MC+ estimators, and prove their unbiasedness. Based on the estimated degrees of freedom, we propose an estimator of the noise level for proper choice of the penalty level. For full rank designs and general sub-quadratic penalties, we provide necessary and sufficient conditions for the continuity of the penalized LSE. Simulation results

Received January 2009; revised June 2009.

[1]Supported in part by NSF Grants DMS-05-04387, DMS-06-04571, DMS-08-04626 and NSA Grant MDS-904-02-1-0063.

overwhelmingly support our claim of superior variable selection properties and demonstrate the computational efficiency of the proposed method.

**1. Introduction.** Variable selection is fundamental in statistical analysis of high-dimensional data. With a proper selection method and under suitable conditions, we are able to build consistent models which are easy to interpret, to avoid over fitting in prediction and estimation, and to identify relevant variables for applications or further study. Consider a linear model in which a response vector $\mathbf{y} \in \mathbb{R}^n$ depends on $p$ predictors $\mathbf{x}_j \in \mathbb{R}^n$, $j = 1, \ldots, p$, through a linear combination $\sum_{j=1}^p \beta_j \mathbf{x}_j$. For small $p$, subset selection methods can be used to find a good guess of the pattern $\{j : \beta_j \neq 0\}$. For example, one may impose a proper penalty on the number of selected variables based on the AIC [Akaike ([1973](#))], $C_p$ [Mallows ([1973](#))], BIC [Schwarz ([1978](#))], RIC [Foster and George ([1994](#))] or a data driven method. For large $p$, subset selection is not computationally feasible, so that continuous penalized or gradient threshold methods are typically used.

Let $\|\cdot\|$ be the Euclidean norm. Consider a penalized squared loss

$$(1.1) \qquad L(\mathbf{b}; \lambda) \equiv (2n)^{-1} \|\mathbf{y} - \mathbf{Xb}\|^2 + \sum_{j=1}^p \rho(|b_j|; \lambda)$$

with a penalty $\rho(t; \lambda)$ indexed by $\lambda \geq 0$, in the linear regression model

$$(1.2) \qquad \mathbf{y} = \sum_{j=1}^p \beta_j \mathbf{x}_j + \boldsymbol{\varepsilon},$$

where $\mathbf{X} \equiv (\mathbf{x}_1, \ldots, \mathbf{x}_p)$, $\boldsymbol{\beta} \equiv (\beta_1, \ldots, \beta_p)'$ and $\boldsymbol{\varepsilon} \sim N(0, \sigma^2 \mathbf{I}_n)$. Assume the penalty $\rho(t; \lambda)$ is nondecreasing in $t$ and has a continuous derivative $\dot{\rho}(t; \lambda) = (\partial/\partial t)\rho(t; \lambda)$ in $(0, \infty)$. Assume further $\dot{\rho}(0+; \lambda) > 0$, so that minimizers of (1.1) have variable selection features with zero components [Donoho et al. ([1992](#))]. Changing the index $\lambda$ if necessary, we assume $\dot{\rho}(0+; \lambda) = \lambda$ whenever $\dot{\rho}(0+; \lambda) < \infty$, so that $\lambda$ has the interpretation as the threshold level for the individual coefficients $\beta_j$ under the standardization $\|\mathbf{x}_j\|^2 = n$.

A widely used penalized lease squares estimator (LSE) is the LASSO [Tibshirani ([1996](#))] or equivalently Basis Pursuit [Chen and Donoho ([1994](#))], with $\rho(t; \lambda) = \lambda|t|$. The LASSO is easy to compute [Osborne, Presnell and Turlach ([2000a](#), [2000b](#)) and Efron et al. ([2004](#))] and has the interpretation as boosting [Schapire ([1990](#)), Freund and Schapire ([1996](#)) and Friedman, Hastie and Tibshirani ([2000](#))]. Throughout the paper, let

$$(1.3) \qquad A^o \equiv \{j : \beta_j \neq 0\} \quad \text{and} \quad d^o \equiv |A^o| = \#\{j : \beta_j \neq 0\}$$

unless otherwise stated. Under a strong irrepresentable condition on the normalized Gram matrix $\boldsymbol{\Sigma} \equiv \mathbf{X}'\mathbf{X}/n$, Meinshausen and Buhlmann ([2006](#)),

Tropp (2006), Zhao and Yu (2006) and Wainwright (2006) proved that the LASSO is variable selection consistent

$$(1.4) \qquad P\{\widehat{A} = A^o\} \to 1 \qquad \text{with } \widehat{A} \equiv \{j : \widehat{\beta}_j \neq 0\},$$

provided that $\min_{\beta_j \neq 0} |\beta_j|/\lambda$ is greater than the $\ell_\infty \to \ell_\infty$ norm of the inverse of a diagonal sub-matrix of $\mathbf{\Sigma}$ of rank $d^o$, among other regularity conditions on $\{\lambda, n, p, \varepsilon\}$. However, the strong irrepresentable condition, which essentially requires the $\ell_\infty \to \ell_\infty$ norm of a $(p - d^o) \times d^o$ matrix of $\mathbf{\Sigma}$ to be uniformly less than 1, is quite restrictive for moderately large $d^o$, and that due to the estimation bias, the condition is nearly necessary for the LASSO to be selection consistent. Here the bias of a penalized LSE is treated as its estimation error when $\varepsilon = \mathbf{0}$. Under a relatively mild sparse Riesz condition on the $\ell_2 \to \ell_2$ norm of sub-Gram matrices and their inverses up to a certain rank, Zhang and Huang (2008) proved that the dimension $|\widehat{A}|$ for the LASSO selection is of the same order as $d^o$ and that the LASSO selects all variables with $|\beta_j|$ above a certain quantity of the order $\sqrt{d^o}\lambda$. These results are still unsatisfactory in view of the possibility of incorrect selection and the extra factor $\sqrt{d^o}$ with the condition on the order of $|\beta_j|$ for correct selection, compared with the threshold level $\lambda$. Again, due to the estimation bias of the LASSO, the extra factor $\sqrt{d^o}$ cannot be completely removed under the sparse Riesz condition. From these points of view, the bias of the LASSO severely interferes with variable selection when $p$ and $d^o$ are both large.

Prior to the above mentioned studies about the interference of the bias of the LASSO with accurate variable selection or conditions to limit such interference, Fan and Li (2001) raised the concern of the effect of the bias of more general penalized estimators on estimation efficiency. They pointed out that the bias of penalized estimators can be removed almost completely by choosing a constant penalty beyond a second threshold level $\gamma\lambda$, and carefully developed the SCAD method [Fan (1997)] with the penalty $\lambda \int_0^t \min\{1, (\gamma - x/\lambda)_+/(\gamma - 1)\} \, dx$, $\gamma > 2$. Iterative algorithms were developed there and in Hunter and Li (2005) and Zou and Li (2008) to approximate a local minimizer of the SCAD penalized loss for fixed $(\lambda, \gamma)$. For penalized methods with unbiasedness and selection features, Fan and Peng (2004) proved the existence, variable selection consistency (1.4) and asymptotic estimation efficiency of some local minimizer of the penalized loss under the dimensionality constraint $p = o(n^r)$ with $r = 1/3, 1/4$ or $1/5$ depending on regularity conditions. Their results apply to general classes of loss and penalty functions but do not address the uniqueness of the solution or provide methodologies for finding or approximating the local minimizer with the stated properties, among potentially many local minimizers. A major cause of computational and analytical difficulties in these studies of nearly unbiased selection methods is the nonconvexity of the minimization problem.

A number of recent papers have considered LASSO-like or LASSO-based convex minimization procedures. Candés and Tao (2007) proposed a Dantzig selector and provided elegant probabilistic upper bounds for the $\ell_2$ loss for the estimation of $\boldsymbol{\beta}$. However, while the Dantzig selector and LASSO have been found to perform similarly in simulation studies [Efron, Hastie and Tibshirani (2007), Meinshausen, Rocha and Yu (2007) and Candés and Tao (2007), page 2401], little is known about the selection consistency of the Dantzig selector. Multiple-stage methods either share certain disadvantages of the LASSO for variable selection or require additional nontechnical side conditions, compared with our results. Current theory on such procedures has been focused on fixed $p$ or $d^o$, while the most interesting case is $p \gg n > d^o \to \infty$. Post-LASSO selection [Meinshausen (2007)] or bootstrapped LASSO [Bach (2008)] may not recover false nondiscovery of the LASSO (Section 6.5). Adaptive LASSO [Zou (2006), Huang, Ma and Zhang (2008) and Zou and Li (2008)] requires an initial estimator of $\boldsymbol{\beta}$ based on which small penalty levels could be assigned to most $\beta_j \neq 0$ and large penalty levels to most $\beta_j = 0$. The nonnegative garrotte estimator [Yuan and Lin (2007)] requires an initial estimator to be within $o(\lambda)$ from $\boldsymbol{\beta}$. For $p \gg n$, correlation screening [Fan and Lv (2008)] requires $A^o$ to be a subset of the indices of the $m$ largest values of $|\mathbf{x}_j'\mathbf{y}|/\|\mathbf{x}_j\|$ with a certain $m \leq n$.

The main purpose of this paper is to propose and study an MC+ methodology. The MC+ provides a fast algorithm for nearly unbiased concave penalized selection in the linear model (1.2). The selection consistency (1.4) holds for the computed MC+ solution at the *universal penalty level* $\lambda_{\mathrm{univ}} \equiv \sigma\sqrt{(2/n)\log p}$ [Donoho and Johnstone (1994b)], without assuming the strong irrepresentable condition or requiring $\min_{\beta_j \neq 0}|\beta_j|/\lambda_{\mathrm{univ}}$ to be greater than a quantity of the order $\sqrt{d^o}$ or the $\ell_\infty \to \ell_\infty$ norm of a matrix of rank $d^o$. This selection consistency holds up to dimension $d^o \leq d_*$, including the case of $p \gg n > d^o \to \infty$, and this upper bound $d_*$, determined by the sparse Riesz condition on $\mathbf{X}$, could be as large as $n/\log(p/n)$. We further prove that the $\ell_q$ loss of the MC+ attains minimax convergence rates in probability for the estimation of $\boldsymbol{\beta}$ in $\ell_r$ balls with $0 < r \leq 1 \wedge q \leq 2$. We also consider a general theory of penalized LSE, including the continuity of estimators, unbiased estimation of risk, and the estimation of noise level, in addition to variable selection and the estimation of $\boldsymbol{\beta}$. This paper is written based on Zhang (2007b), an April, 2007 Rutgers University Technical Report containing all the results in Sections 3, 4 and 5 with more extensive discussion of the PLUS algorithm and less explicit constants in the selection consistency theorems. A brief description of Zhang (2007b) can be found in Zhang (2008), which contains some additional simulation results.

**2. A sketch of main results.** We provide a brief description of the MC+ method and our main results, along with some crucial concepts, conditions and necessary notation.

2.1. *The MC+.* The MC+ has two components: a *minimax concave penalty* (MCP) and a *penalized linear unbiased selection* (PLUS) algorithm. The MCP is defined as

$$(2.1) \qquad \rho(t; \lambda) = \lambda \int_0^t (1 - x/(\gamma\lambda))_+ \, dx$$

with a regularization parameter $\gamma > 0$. It minimizes the maximum concavity

$$(2.2) \qquad \kappa(\rho) \equiv \kappa(\rho; \lambda) \equiv \sup_{0 < t_1 < t_2} \{\dot{\rho}(t_1; \lambda) - \dot{\rho}(t_2; \lambda)\}/(t_2 - t_1)$$

subject to the following unbiasedness and selection features:

$$(2.3) \qquad \dot{\rho}(t; \lambda) = 0 \qquad \forall t \geq \gamma\lambda, \qquad \dot{\rho}(0+; \lambda) = \lambda.$$

For $A \subseteq \{1, \dots, p\}$, define sub-design and sub-Gram matrices

$$(2.4) \qquad \mathbf{X}_A \equiv (\mathbf{x}_j, j \in A)_{n \times |A|}, \qquad \boldsymbol{\Sigma}_{A,B} \equiv \mathbf{X}'_A \mathbf{X}_B / n, \qquad \boldsymbol{\Sigma}_A \equiv \boldsymbol{\Sigma}_{A,A}.$$

Let $d^*$ be a positive integer. The penalized loss (1.1) is *sparse convex* with rank $d^*$ if it is convex in all models $\{\mathbf{b} : b_j = 0 \ \forall j \notin A\}$ with $|A| \leq d^*$. This sparse convexity condition holds if the convexity of the squared loss $\|\mathbf{y} - \mathbf{X}\mathbf{b}\|^2/(2n)$ overcomes the concavity of the penalty in all such sparse models with $|A| \leq d^*$, or equivalently

$$(2.5) \qquad \kappa(\rho; \lambda) < \min_{|A| \leq d^*} c_{\min}(\boldsymbol{\Sigma}_A) \qquad \text{where } c_{\min}(\boldsymbol{\Sigma}_A) \equiv \min_{\|\mathbf{u}\|=1} \|\boldsymbol{\Sigma}_A \mathbf{u}\|.$$

Although the unbiasedness and selection features (2.3) preclude convex penalties, the MCP provides the sparse convexity to the broadest extent by minimizing the maximum concavity (2.2). This is a natural motivation for the MCP. The MCP achieves $\kappa(\rho; \lambda) = 1/\gamma$. A larger value of its regularization parameter $\gamma$ affords less unbiasedness and more concavity. For each penalty level $\lambda$, the MCP provides a continuum of penalties with the $\ell_1$ penalty at $\gamma = \infty$ and the "$\ell_0$ penalty" as $\gamma \to 0+$.

Given a penalty $\rho(\cdot; \cdot)$, $\lambda \oplus \widehat{\boldsymbol{\beta}} \in \mathbb{R}^{1+p}$ is a critical point of the penalized loss in (1.1) if

$$(2.6) \qquad \begin{cases} \mathbf{x}'_j(\mathbf{y} - \mathbf{X}\widehat{\boldsymbol{\beta}})/n = \operatorname{sgn}(\widehat{\beta}_j)\dot{\rho}(|\widehat{\beta}_j|; \lambda), & \widehat{\beta}_j \neq 0, \\ |\mathbf{x}'_j(\mathbf{y} - \mathbf{X}\widehat{\boldsymbol{\beta}})/n| \leq \lambda, & \widehat{\beta}_j = 0, \end{cases}$$

where $\operatorname{sgn}(t) \equiv I\{t > 0\} - I\{t < 0\}$. For convex penalized loss, (2.6) is the Karush–Kuhn–Tucker (KKT) condition for the global minimization of (1.1). In general, solutions of (2.6) include all local minimizers of $L(\cdot; \lambda)$ for all $\lambda$. The graph of the solutions of (2.6) is studied in Section 3. Consider

$$(2.7) \qquad \lambda^{(x)} \oplus \widehat{\boldsymbol{\beta}}^{(x)} \equiv \begin{cases} \text{a continuous path of solutions of (2.6) in } \mathbb{R}^{1+p} \\ \text{with } \widehat{\boldsymbol{\beta}}^{(0)} = \mathbf{0} \text{ and } \lim_{x \to \infty} \lambda^{(x)} = 0. \end{cases}$$

For the MCP, we prove in Section 3.3 that almost everywhere in $(\mathbf{X}, \mathbf{y})$, a path (2.7) uniquely exists up to continuous transformations of $x$ from $[0, \infty)$ onto $[0, \infty)$ and that $\widehat{\boldsymbol{\beta}}^{(x)}$ ends at a point of global least squares fit as $x \to \infty$. Thus, in the graph of the solutions of (2.6), (2.7) provides a unique branch encompassing from the origin $\boldsymbol{\beta} = \mathbf{0}$ to an optimal fit. We call (2.7) the main branch of the solution graph. For concave penalties, solutions of (2.6) may form additional branches as loops not connected to the origin (Figure 3). In the PLUS algorithm, the integer part of $x$ in (2.7) represents the number of computational steps and the fraction part represents the linear interpolation between steps as in (3.8).

Given a penalty level $\lambda$, we propose as a variable selector and an estimator of $\boldsymbol{\beta}$

$$(2.8) \qquad \widehat{\boldsymbol{\beta}}(\lambda) \equiv \widehat{\boldsymbol{\beta}}^{(x_\lambda)} \qquad \text{in (2.7) with } x_\lambda = \inf\{x \geq 0 : \lambda^{(x)} \leq \lambda\},$$

or equivalently the solution when the penalty level $\lambda$ is first reached in the path. The estimator (2.8) and the global minimum of (1.1) may not be the same for nonconvex penalized loss. Still, the uniqueness of (2.7) implies that $\widehat{\boldsymbol{\beta}}(\lambda)$ is uniquely defined, including the case of $p > n$. We call (2.8) the MC+ if the MCP (2.1) is used in (2.6) and thus (2.7).

The PLUS algorithm computes the main branch (2.7) of the solution graph of (2.6) for quadratic spline penalty functions of the form $\rho(t; \lambda) = \lambda^2 \rho(t/\lambda)$. The PLUS is described in detail and studied in Section 3. For the quadratic spline penalties, the graph of solutions of (2.6) is piecewise linear and so is its main branch (2.7). The PLUS differs from existing nonconvex minimization algorithms in three important aspects: (i) it computes the exact value of local minimizers instead of iteratively approximating them; (ii) it computes a path of possibly multiple solutions for the entire range of the penalty level $\lambda \geq 0$ instead of a single solution for a fixed $\lambda$; (iii) it computes multiple local minimizers for individual $\lambda$ by tracking along its path of solutions for different values of $\lambda$ instead of trying to jump from the domain of attraction of one solution to another for a fixed $\lambda$. In each step, the PLUS computes one line segment in its path between two turning points, and its computational cost is the same as the LARS [Efron et al. (2004)] per step. The MC+ with larger regularization parameter $\gamma$ provides smoother estimators and computationally less complex path, but larger bias and less accurate variable selection. The MC+ path converges to the LASSO path as $\gamma \to \infty$.

2.2. *Some simulation results and heuristics.* The proposed MC+ provides fast, continuous, nearly unbiased and accurate variable selection in high-dimensional linear regression, as our theoretical and numerical results support.

Table 1

*Performance of LASSO, MC+ and SCAD+ in experiment 1. 100 replications, $n = 300$,*
*$p = 200$, $\beta_*/\sigma = 1/2$, $\gamma = 2/(1 - \max_{j \neq k} |\mathbf{x}_j' \mathbf{x}_k|/n)$, $\overline{\gamma} = 2.69$, CS $\equiv I\{\widehat{A} = A^o\}$,*
*$\mathrm{SE}_{\boldsymbol{\beta}} \equiv \|\widehat{\boldsymbol{\beta}} - \boldsymbol{\beta}\|^2$, $k \equiv \#(steps)$, $\log(\sigma\lambda/(\widehat{\sigma}\lambda_{\mathrm{univ}})) = integer/20$*

| | | $d^o = 10$ | | | $d^o = 20$ | | | $d^o = 40$ | | |
| $\lambda/\widehat{\sigma}$ | | **LASSO** | **MC+** | **SCAD+** | **LASSO** | **MC+** | **SCAD+** | **LASSO** | **MC+** | **SCAD+** |
|---|---|---|---|---|---|---|---|---|---|---|
| $\lambda/\widehat{\sigma}$ | $\overline{\mathrm{CS}}$ | 0.45 | **0.77** | 0.71 | 0.09 | **0.87** | 0.62 | 0.00 | **0.81** | 0.27 |
| $= \lambda_{\mathrm{univ}}/\sigma$ | $\overline{\mathrm{SE}_{\boldsymbol{\beta}}}$ | 0.340 | **0.063** | 0.131 | 0.831 | **0.160** | 0.480 | 2.097 | **0.452** | 1.842 |
| $= 0.188$ | $\overline{k}$ | 12 | 16 | 26 | 23 | 31 | 50 | 47 | 63 | 127 |
| Fixed $\lambda/\widehat{\sigma}$ | $\lambda/\widehat{\sigma}$ | 0.266 | 0.248 | 0.248 | 0.257 | 0.231 | 0.195 | 0.231 | 0.195 | 0.169 |
| for max $\overline{\mathrm{CS}}$ | $\overline{\mathrm{CS}}$ | 0.88 | **0.98** | 0.92 | 0.44 | **0.97** | 0.70 | 0.01 | **0.83** | 0.45 |
| | $\overline{k}$ | 11 | 11 | 17 | 21 | 23 | 47 | 44 | 60 | 149 |
| Fixed $\lambda/\widehat{\sigma}$ | $\lambda/\widehat{\sigma}$ | 0.076 | 0.153 | 0.138 | 0.060 | 0.138 | 0.124 | 0.042 | 0.138 | 0.120 |
| for min $\overline{\mathrm{SE}_{\boldsymbol{\beta}}}$ | $\overline{\mathrm{SE}_{\boldsymbol{\beta}}}$ | 0.154 | 0.043 | **0.041** | 0.287 | 0.082 | **0.080** | 0.502 | 0.167 | **0.161** |
| | $\overline{k}$ | 41 | 22 | 34 | 65 | 43 | 67 | 102 | 84 | 169 |

Table 1 presents the results of experiment 1 of our simulation study to demonstrate the superior selection accuracy and competitive computational complexity of the MC+, compared with the LASSO and SCAD. Since there are quite a few different ways of (approximately) computing possibly different SCAD local minimizers, we denote by SCAD+ the PLUS solution of the SCAD. We measure the selection accuracy by the proportion $\overline{\mathrm{CS}}$ of replications with the correct selection CS $\equiv I\{\widehat{A} = A^o\}$, and the computational complexity by the average $\overline{k}$ of the number of the PLUS steps. In this experiment, $(n, p) = (300, 200)$, $\mathbf{y}$ is generated with $\beta_j = \pm\beta_*$ for $j \in A^o$ and $\boldsymbol{\varepsilon} \sim N(0, \mathbf{I}_n)$ in (1.2), and $\mathbf{x}_j$ are generated by greedy sequential sampling (Section 6.1) of groups of 10 most correlated vectors from a pool of 600 vectors from the sphere $\{\mathbf{x} : \|\mathbf{x}\| = \sqrt{n}\}$. The design $\mathbf{X}$, $A^o$, the signs of $\boldsymbol{\beta}$ and $\boldsymbol{\varepsilon}$ are drawn independently for the 100 replications with $d^o \in \{10, 20, 40\}$. The $\widehat{\sigma}^2$ is the residual mean squares with 100 degrees of freedom in the full 200-dimensional model.

The dimensions $(n, p, d^o)$ in experiment 1 are moderate, but larger than those in some recent simulation studies of other nonconvex minimization algorithms. This modest setting allows us to demonstrate the significance of the impact of $d^o = \#\{j : \beta_j \neq 0\}$, and thus the estimation bias, on the selection consistency of the LASSO in the absence of difficulties involving ultrahigh dimensionality or the singularity with rank($\mathbf{X}$) $< p$. More simulation results are presented in Section 6 with $(n, p) = (300, 2000)$, $(600, 3000)$, $(100, 2000)$ and $(200, 10{,}000)$ to demonstrate the scalability of the PLUS algorithm, among other issues.

Why is the MC+ able to avoid both the interference of estimation bias with variable selection and the computational difficulties with nonconvex minimization? A short, heuristic explanation is that for standardized $\|\mathbf{x}_j\| = \sqrt{n}$ and a carefully chosen $\gamma > 1$, the condition

$$(2.9) \qquad \beta_* \equiv \min\{|\beta_j| : j \in A^o\} > \gamma\lambda \qquad \text{with } \lambda \geq \lambda_{\text{univ}} \equiv \sigma\sqrt{(2/n)\log p},$$

and the sparsity of $\boldsymbol{\beta}$ are allowed to match the extent of the unbiasedness and sparse convexity of the MC+. The lower bound for $\beta_*$ in (2.9) allows unbiased selection of all $j \in A^o$, while the lower bound for $\lambda$ prevents selection of variables outside $A^o$ given the selection of all variables in $A^o$. Thus, (2.9) guarantees with large probability that the LSE

$$(2.10) \qquad \widehat{\boldsymbol{\beta}}^o \equiv \arg\min_{\mathbf{b}}\{\|\mathbf{y} - \mathbf{X}\mathbf{b}\|^2 : b_j = 0 \ \forall j \notin B\}$$

with the oracular choice $B = A^o$, is one of the local minimizers of the penalized loss. Meanwhile, the sparse convexity (2.5) provides the uniqueness among sparse solutions of (2.6) and controls the computational complexity of the MC+.

This argument does not work with the LASSO due to the estimation bias. Let $\widetilde{\boldsymbol{\beta}}^o$ be the $\ell_1$ oracle with $\widetilde{\boldsymbol{\beta}}^o_{A^o}(\lambda) = \widehat{\boldsymbol{\beta}}^o_{A^o} - \lambda\boldsymbol{\Sigma}^{-1}_{A^o}\text{sgn}(\boldsymbol{\beta}_{A^o})$ and $\widetilde{\boldsymbol{\beta}}^o_j(\lambda) = 0$ for $j \notin A^o$. By the KKT condition, $\text{sgn}(\widehat{\boldsymbol{\beta}}(\lambda)) = \text{sgn}(\boldsymbol{\beta})$ for the LASSO if and only if (iff) $|\mathbf{x}'_j(\mathbf{y} - \mathbf{X}\widetilde{\boldsymbol{\beta}}^o(\lambda))|/n \leq \lambda$ and $\text{sgn}(\widetilde{\boldsymbol{\beta}}^o(\lambda)) = \text{sgn}(\boldsymbol{\beta})$. However, $\widetilde{\boldsymbol{\beta}}^o(\lambda)$ is biased with $E\widetilde{\boldsymbol{\beta}}^o_{A^o}(\lambda) - \boldsymbol{\beta}_{A^o} = -\lambda\boldsymbol{\Sigma}^{-1}_{A^o}\text{sgn}(\boldsymbol{\beta}_{A^o}) \neq 0$.

2.3. *Selection consistency.* We study the selection consistency of the penalized LSE under the sparse Riesz condition (SRC) on $\mathbf{X}$: for suitable $0 < c_* \leq c^* < \infty$ and rank $d^*$,

$$(2.11) \qquad c_* \leq \min_{|A| \leq d^*} c_{\min}(\boldsymbol{\Sigma}_A) \leq \max_{|A| \leq d^*} c_{\max}(\boldsymbol{\Sigma}_A) \leq c^*,$$

where $\boldsymbol{\Sigma}_A$ is as in (2.4) and $c_{\max}(\mathbf{M})$ is the largest eigenvalue of $\mathbf{M}$. Conditions on $\mathbf{X}$ and $\boldsymbol{\beta}$ must be configured to accommodate each other in our theorems. In this subsection, we study selection consistency for $d^o \leq d_* = d^*/(c^*/c_* + 1/2)$. In the next subsection, we study estimation by comparing $\widehat{\boldsymbol{\beta}}$ and the oracle estimator (2.10) with $|B| \leq d_*$. Section 4 covers more general configurations. Although $\{d^*, c_*, c^*\}$ are all allowed to depend on $n$, the SRC is easier to understand with fixed $\{c_*, c^*\}$ and large $d^* \equiv d_n^*$, which asserts the equivalence of the norms $\|\mathbf{X}\mathbf{b}\|/\sqrt{n}$ and $\|\mathbf{b}\|$ up to $\#\{j : b_j \neq 0\} = d^*$. Define $\widetilde{p}_\epsilon \equiv \widetilde{p}_{p,d^o,m,\epsilon}$ by

$$(2.12) \qquad 2\log\widetilde{p}_\epsilon - 1 - \log(2\log\widetilde{p}_\epsilon) = (2/m)\left\{\log\binom{p - d^o}{m} + \log(1/\epsilon)\right\}$$

for nonnegative integers $m \in [1, p - d^o]$ and reals $\epsilon \in (0, 1]$. Note that $2\log\widetilde{p}_\epsilon \geq 1$.

THEOREM 1. *Suppose (1.2) holds with $\|\mathbf{x}_j\|^2 = n$. Let $A^o$, $d^o$ and $\widehat{A}$ be as in (1.3) and (1.4) and $\widehat{\boldsymbol{\beta}}^o$ be as in (2.10) with $B = A^o$. Suppose (2.11) holds and $d^o \leq d_* = d^*/(c^*/c_* + 1/2)$. Let $\lambda_{1,\epsilon} = \sigma \sqrt{(2/n)\log((p - d^o)/\epsilon)}$ and $\lambda_{2,\epsilon} \geq \max\{2\sqrt{c^*}\sigma\sqrt{(2/n)\log\widetilde{p}_\epsilon}, \lambda_{1,\epsilon}\}$, where $\epsilon \in (0, 1]$ is fixed and $\widetilde{p}_\epsilon$ is defined with $m = d^* - d^o$. Let $w^o$ be the largest diagonal element of $\boldsymbol{\Sigma}_{A^o}^{-1}$. Let $\widehat{\boldsymbol{\beta}} = \widehat{\boldsymbol{\beta}}(\widehat{\lambda})$ with a deterministic or random $\widehat{\lambda}$, where $\widehat{\boldsymbol{\beta}}(\lambda)$ is the MC+ selector (2.8) with $\gamma \geq c_*^{-1}\sqrt{4 + c_*/c^*}$. Then*

$$(2.13) \quad P\{\widehat{\boldsymbol{\beta}} \neq \widehat{\boldsymbol{\beta}}^o \text{ or } \text{sgn}(\widehat{\boldsymbol{\beta}}) \neq \text{sgn}(\boldsymbol{\beta})\} \leq P\{\widehat{\lambda} \notin [\lambda_{1,\epsilon}, \lambda_{2,\epsilon}]\} + (3/2 + 1/\sqrt{2})\epsilon,$$

*provided that $\beta_* \equiv \min_{j \in A^o}|\beta_j| \geq \sigma\sqrt{w^o(2/n)\log(d^o/\epsilon)} + \gamma\lambda_{2,\epsilon}$. Moreover, (1.4) holds and the MC+ estimator $\widehat{\boldsymbol{\beta}}$ achieves the estimation efficiency of the oracle LSE $\widehat{\boldsymbol{\beta}}^o$, provided that $P\{\lambda_{1,\epsilon} \leq \widehat{\lambda} \leq \lambda_{2,\epsilon}\} \to 1$ and $\epsilon^{-1} \vee \min\{p - d^o, \widetilde{p}_1, \sqrt{n/w^o}(\beta_* - \gamma\lambda_{2,\epsilon})/\sigma\} \to \infty$.*

COROLLARY 1. *Let $\lambda_{\text{univ}} \equiv \sigma\sqrt{(2/n)\log p}$. Suppose $\|\mathbf{x}_j\|^2 = n$, $d^*/(c^*/c_* + 1/2) \geq d^o \to \infty$, $\gamma \geq c_*^{-1}\sqrt{4 + c_*/c^*}$ and $\beta_* \geq \sigma\sqrt{w^o(2/n)\log d^o} + \gamma\max\{2 \times \sqrt{c^*}\sigma\sqrt{(2/n)\log\widetilde{p}_1}, \lambda_{\text{univ}}\}$ in (1.2), (2.11) and (2.1). Then $P\{\widehat{\boldsymbol{\beta}}(\lambda_{\text{univ}}) \neq \widehat{\boldsymbol{\beta}}^o$ or $\text{sgn}(\widehat{\boldsymbol{\beta}}(\lambda_{\text{univ}})) \neq \text{sgn}(\boldsymbol{\beta})\} \to 0$.*

A lower bound condition on $\beta_*$ can be viewed as an information requirement for selection consistency. A variation below of Proposition 1 in Zhang (2007c) asserts that the condition on $\beta_*$ in Theorem 1 is optimal up to a factor of $4\gamma\sqrt{c^*}(1 + o(1))$ when $\log d^o = o(\log p)$.

PROPOSITION 1. *For $\boldsymbol{\beta} \in \mathbb{R}^p$ let $A^o$ and $d^o$ be as in (1.3), $\beta_* \equiv \min_{j \in A^o}|\beta_j|$, and $\mathbf{y}$ be as in (1.2) with $\|\mathbf{x}_j\|^2 = n$. Let $p$, $d^o$ and $\sigma > 0$ be dependent on $n$ with $p - d^o \to \infty$. Then*

$$\liminf_{n \to \infty} \inf_{(\mathbf{X}, \mathbf{y}) \to \widehat{A}} \sup_{|A^o| = d^o} \sup_{\beta_* = c\lambda_{1,1}} P\{\widehat{A} \neq A^o\} \geq 1 - 4c^2 \qquad \forall c > 0,$$

*where $\lambda_{1,1} = \sigma\sqrt{(2/n)\log(p - d^o)}$ and the infimum is taken over all Borel mappings.*

REMARK 1. Since (2.13) is nonasymptotic, $\{p, d^*, c_*, c^*, d^o, \boldsymbol{\beta}, \sigma, \epsilon\}$ are all allowed to depend on $n$. The requirement $d^o \leq d_* = d^*/(c^*/c_* + 1/2)$ could be viewed as a condition on the sparsity of $\boldsymbol{\beta}$ given $\{d^*, c_*, c^*\}$. On the other hand, for given $d^o \equiv d^o_n$ it is closely related to the restricted isometry constant $\delta_d \equiv \max\{|\|\boldsymbol{\Sigma}_A\mathbf{u}\| - 1| : |A| = d, \|\mathbf{u}\| = 1\}$ [Candés and Tao (2005)], although $c^* > 2$ is allowed in (2.11). For example, $d^o \leq d^*/(c^*/c_* + 1/2)$ is a consequence of $\delta_{3d^o} \leq 3/7$ with explicit $d^* = 3d^o$, $c_* = 4/7$ and $c^* = 10/7$. With larger $\lambda_{2,\epsilon}/\sqrt{\sigma^2\log\widetilde{p}_\epsilon}$ and $\gamma$, Theorem 5 allows fixed $d^*/d_* > (c^*/c_* +$ $1)/2$, which is a consequence of $\delta_{2d^o} < 1/2$ or $\delta_{3d^o} < 2/3$. See Remark 5 in Section 4.

REMARK 2. For $p \gg n$, random matrix theory provides the possibility of $d^o \asymp n/\log(p/n)$. For example, if the rows of $\mathbf{X}$ are i.i.d. Gaussian vectors with $E\mathbf{X} = 0$ and $c_1 \leq E\|\mathbf{Xb}\|^2/n \leq c_2$ for all $\|\mathbf{b}\| = 1$, then $P\{(2.11) \text{ holds}\} \to 1$ with fixed $c_* = (1-\delta)^2 c_1$, $c^* = (1+\delta)^2 c_2$ and $d^* = \max\{d : \sqrt{d/n}(1 + \sqrt{2 + 2\log(p/d)}) \leq \delta\}$, where $0 < \delta < 1$ is fixed [Davidson and Szarek (2001), Candes and Tao (2007), Wainwright (2006) and Zhang and Huang (2008)].

REMARK 3. The condition $\boldsymbol{\varepsilon} \sim N(0, \sigma^2 \mathbf{I}_n)$ is not essential. In particular, Corollary 1 holds if the normality assumption is replaced by the sub-Gaussian condition $Ee^{\mathbf{x}'\boldsymbol{\varepsilon}} \leq e^{\sigma_1^2 \|\mathbf{x}\|^2/2} \ \forall \mathbf{x}$, provided $\sigma_1^2 < \sigma^2$. See Section 7.3 and Lemma 2.

Theorem 1 compares favorably with existing results in the required regularity of $\mathbf{X}$ and the information content in the data as measured in $\beta_* \equiv \min_{\beta_j \neq 0} |\beta_j|$. For the LASSO, a bound similar to (2.13) on selection consistency essentially requires

$$(2.14) \qquad \beta_* \geq \sigma\sqrt{w^o(2/n)\log(d^o/\epsilon)} + \theta_1^*\lambda \quad \text{and} \quad \lambda \geq \lambda_{1,\epsilon}/(1-\theta_2^*)_+,$$

where $\theta_1^* \equiv \|\boldsymbol{\Sigma}_{A^o}^{-1}\operatorname{sgn}(\boldsymbol{\beta}_{A^o})\|_\infty$ and $\theta_2^* \equiv \|\boldsymbol{\Sigma}_{(A^o)^c,A^o}\boldsymbol{\Sigma}_{A^o}^{-1}\operatorname{sgn}(\boldsymbol{\beta}_{A^o})\|_\infty$ [Meinshausen and Buhlmann (2006), Tropp (2006), Zhao and Yu (2006) and Wainwright (2006)]. The maxima of $\theta_1^*$ and $\theta_2^*$ over the unknown $\operatorname{sgn}(\boldsymbol{\beta}_{A^o})$ are, respectively, the norms $\|\boldsymbol{\Sigma}_{A^o}^{-1}\|_\infty$ and $\|\boldsymbol{\Sigma}_{(A^o)^c,A^o}\boldsymbol{\Sigma}_{A^o}^{-1}\|_\infty$ for linear mappings between $\ell_\infty$ spaces. Consider the case of $d^o \to \infty$. The strong irrepresentable condition, which requires $\theta_2^* < 1$ uniformly strictly, is restrictive since $\|\boldsymbol{\Sigma}_{(A^o)^c,A^o}\boldsymbol{\Sigma}_{A^o}^{-1}\|_\infty$ is not length normalized. For $\log(p/d^*) \asymp \log p \to \infty$, $\sigma\sqrt{(2/n)\log\widetilde{p}_\epsilon} = (1 + o(1))\lambda_{1,\epsilon}$ by (2.12), so that Theorem 1 replaces $\theta_1^*/(1-\theta_2^*)_+$ in (2.14) by $2\gamma\sqrt{c^*}$ as a required lower bound for the signal-to-noise ratio (SNR) $\beta_*/\lambda_{1,\epsilon}$. For $\log p = (1+o(1))\log d^*$, for example, $p \asymp n\log n$ and $d^* \asymp n/\log\log n$, $\log\widetilde{p}_1 = o(1)\log p$, so that Corollary 1 simply requires $\beta_* \geq \sigma\sqrt{w^o(2/n)\log d^o} + \gamma\lambda_{\text{univ}}$ for (1.4) when $c^* = O(1)$. A commonly used bound is $\theta_1^* \leq \sqrt{d^o}/c_{\min}(\boldsymbol{\Sigma}_{A^o})$. Wainwright (2006) proved $\|\boldsymbol{\Sigma}_{A^o}^{-1}\|_\infty = O_P(1)$ when the rows of $\mathbf{X}_{A^o}$ are i.i.d. Gaussian vectors with $\|(E\boldsymbol{\Sigma}_{A^o})^{-1}\|_\infty = O(1)$. The adverse effects of large $d^o$ on the LASSO selection are evident in our simulation experiments.

In addition to conditions on $\mathbf{X}$ and $\boldsymbol{\beta}$, Theorem 1 makes significant advances by allowing the exact universal penalty level $\lambda_{\text{univ}}$ for selection consistency (Corollary 1) in the case of a known $\sigma^2$ or $\widetilde{\lambda} = \widehat{\sigma}\sqrt{(2/n)\log p}$ with any consistent upper confidence bound $\widehat{\sigma}$ in the case of unknown $\sigma$, while the penalty level $\lambda$ in (2.14) depends on $A^o$ via the $\ell_\infty$ norm $\theta_2^*$.

From these points of view, the thrust of Theorem 1 is to replace the strong irrepresentable condition by the SRC with $d^o \leq d^*/(c^*/c_* + 1/2)$, to replace the $\ell_\infty \to \ell_\infty$ norm of matrices of rank $d^o$ by the $\ell_2 \to \ell_2$ norm of matrices of rank no greater than $d^o(c^*/c_* + 1/2)$ in the requirement on $\beta_*$, and to completely remove the factor $1/(1 - \theta_2^*)_+$ on $\lambda$, compared with (2.14).

2.4. *Estimation of regression coefficients.* We have shown the selection consistency of the MC+ up to $|A^o| \leq d_* = d^*/(c^*/c_* + 1/2)$ under (2.11). This selection consistency is proved via an upper bound on the false positive in Theorem 6 which naturally leads to performance bounds for the estimation of $\boldsymbol{\beta}$. Although we do not fully address the topic here, we present a theorem to highlight the consequences of our oracle inequalities.

Let $\|\mathbf{b}\|_q = (\sum_{j=1}^p |b_j|^q)^{1/q}$ be the $\ell_q$ norm with the usual extension to $q = \infty$ and $\Theta_{r,R} \equiv \{\mathbf{b} : \|\mathbf{b}\|_r \leq R\}$ be the $\ell_r$ ball. It was proved recently in Ye and Zhang (2009) that for all $1 < r \vee 1 \leq q$ and $0 < \varepsilon < 1$

$$(2.15) \quad \liminf_{p \to \infty} \inf_{\mathbf{X}} \inf_{(\mathbf{X}, \mathbf{y}) \to \widehat{\boldsymbol{\beta}}} \sup_{\boldsymbol{\beta} \in \Theta_{r,R}} P\{\|\widehat{\boldsymbol{\beta}} - \boldsymbol{\beta}\|_q^q \geq (1 - \epsilon) R^r \lambda_{\mathrm{mm}}^{q-r}\} \geq \frac{\epsilon}{3^q}$$

subject to $\|\mathbf{x}_j\|^2 = n$ in (1.2), where the second infimum is taken over all Borel mappings of proper dimension and

$$\lambda_{\mathrm{mm}} \equiv \sigma\{(2/n)\log(\sigma^r p/(n^{r/2} R^r))\}^{1/2},$$

provided that $R^r/\lambda_{\mathrm{mm}}^r \to \infty$ and $n\lambda_{\mathrm{mm}}^2/\sigma^2 \to \infty$. This minimax lower bound for the $\ell^q$ loss is an extension of the lower bound for the minimax $\ell^q$ risk in Donoho and Johnstone (1994a). The following theorem provides sufficient conditions for the PLUS estimator (2.8) to attain this minimax rate.

THEOREM 2. *Let $\kappa \geq 0$ and $\rho(t; \lambda)$ be a penalty satisfying $\lambda(1 - \kappa|t|/\lambda)_+ \leq \dot{\rho}(|t|; \lambda) \leq \lambda$. Suppose (2.11) holds with certain $d^*$ and $c^* \geq c_* \geq \kappa\sqrt{4 + c_*/c^*}$. Let $B$ be a deterministic subset of $\{1, \ldots, p\}$ with $|B| \leq d_* = d^*/(c^*/c_* + 1/2)$. Let $\widehat{\boldsymbol{\beta}}(\lambda)$ be as in (2.8) and $\widehat{\boldsymbol{\beta}}^o$ as in (2.10). Let $\theta_B \equiv \|\mathbf{X}(\boldsymbol{\beta} - E\widehat{\boldsymbol{\beta}}^o)\|/\sqrt{n}$ and $\widetilde{p}_\epsilon$ be as in (2.12) with $m = d^* - |B|$ and $d^o = |B|$.*

(i) *Let $\lambda \geq 2\sqrt{c^*}(\sigma\sqrt{(2/n)\log\widetilde{p}_\epsilon} + \theta_B/\sqrt{m})$. Then, with at least probability $1 - \epsilon/\sqrt{4\log\widetilde{p}_\epsilon}$,*

$$(2.16) \quad c_*\|\widehat{\boldsymbol{\beta}}(\lambda) - \widehat{\boldsymbol{\beta}}^o\| \leq \left\{\sum_{j \in B} \dot{\rho}^2(|\widehat{\beta}_j|; \lambda)\right\}^{1/2} + (\lambda/2)\sqrt{|B|} \leq (3/2)\lambda\sqrt{|B|}.$$

(ii) *Suppose $R^r/\lambda_{\mathrm{mm}}^r = |B|$. Let $\lambda = 2\sqrt{c^*}\{\lambda_{\mathrm{mm}}(1 + \sqrt{2c_*}) + \epsilon_1\sigma/\sqrt{n}\}$ with the $\lambda_{\mathrm{mm}}$ in (2.15) and a fixed $\epsilon_1 > 0$. Let $0 < r \leq 1 \wedge q \leq 2$ and $M_q = (M_{1,q}^{q \wedge 1} + M_{2,q}^{q \wedge 1})^{(1/q) \vee 1}$, where $M_{1,q} = (c^*/c_* + 1/2)^{1/q - 1/2} 3(\sqrt{c^*}/c_*)(1 + \sqrt{2c_*} + $* $\epsilon_2$) and $M_{2,q} = \{(c^*/c_* + \epsilon_2/c_*)^{q/2} + 1\}^{1/q}$ with a fixed $\epsilon_2 > 0$. Then, with $\{p, R, \sigma, d^*, c_*, c^*, M\}$ all allowed to depend on $n$,

$$(2.17) \qquad \sup_{\beta \in \widetilde{\Theta}_{r,R}} P\{\|\widehat{\boldsymbol{\beta}}(\lambda) - \boldsymbol{\beta}\|_q^q \geq M_q^q R^r \lambda_{\mathrm{mm}}^{q-r}\} \to 0$$

as $n\lambda_{\mathrm{mm}}^2/\sigma^2 \to \infty$, where $\widetilde{\Theta}_{r,R} \equiv \{\boldsymbol{\beta} : \sum_{j=1}^p |\beta_j|^r \wedge \lambda_{\mathrm{mm}}^r \leq R^r\} \supseteq \Theta_{r,R}$.

REMARK 4. We may choose the set $B$ in Theorem 2(i) to minimize $\theta_B$ or $\sum_{j \notin B} |\beta_j|$ given a size $|B|$, but we are not confined to these examples. The condition $R^r/\lambda_{\mathrm{mm}}^r = |B|$ in Theorem 2(ii) is not restrictive, since $\mathbf{X}$ and $\sigma$ could be scaled by a bounded factor to meet it. Remarks 1, 2 and 3 are applicable to Theorem 2 with $\gamma = 1/\kappa$.

The oracle inequality (2.16) exhibits the advantage of the MC+ when a fraction of $|\beta_j|$ are of the order $\lambda$, since the MCP with $\gamma = 1/\kappa$ has the smallest possible $\dot{\rho}(t; \lambda) = (1 - \kappa t/\lambda)_+$ under the assumption on the penalty.

For fixed $0 < c_* \leq c^* < \infty$, (2.17) provides the convergence of $\widehat{\boldsymbol{\beta}}(\lambda)$ based on $(\mathbf{X}, \mathbf{y})$ at the minimax rate (2.15), up to $\#\{\text{significant } \beta_j\} \asymp R^r/\lambda_{\mathrm{mm}}^r \leq d_* = d^*/(c^*/c_* + 1/2)$, including the case of $p \gg n \geq d^o \to \infty$. Such uniform convergence rates in $\ell_r$ balls cannot be obtained from existing results requiring penalty levels $\lambda \geq \lambda_{\mathrm{univ}} \equiv \sigma\sqrt{(2/n)\log p}$ in the case of $\lambda_{\mathrm{mm}}/\lambda_{\mathrm{univ}} \to 0$. Theorem 2 closes this gap by allowing $\lambda \asymp \lambda_{\mathrm{mm}}$. We observe that $\lambda_{\mathrm{mm}} < \lambda_{\mathrm{univ}}$ in (2.15) whenever $R > \sigma/\sqrt{n}$. The relevance of smaller $\lambda_{\mathrm{mm}}$ is evident in our simulation experiments where the best penalty levels for estimation are all less than or equal to $\lambda_{\mathrm{univ}}$. See Section 6.1 in addition to Table 1. For recent advances in the LASSO or LASSO-like estimations of $\mathbf{X}\boldsymbol{\beta}$ and $\boldsymbol{\beta}$, we refer to Greenshtein and Ritov (2004), Candés and Tao (2007), Bunea, Tsybakov and Wegkamp (2007), van de Geer (2008), Zhang and Huang (2008) and Meinshausen and Yu (2009).

2.5. *Organization of the rest of the paper.* Section 3 provides an explicit description of the PLUS algorithm and studies the geometry of the solutions of the estimating equation (2.6). Section 4 studies the selection consistency of both the global minimizer of (1.1) and the local solution (2.8) for general penalties. Section 5 develops methods for the estimations of the mean squared error (MSE) of the penalized LSE and the noise level in the linear model (1.2). Section 6 reports simulation results. Section 7 contains some discussion.

**3. The PLUS algorithm and quadratic spline penalties.** We divide this section into three subsections to cover quadratic spline penalties, the PLUS algorithm and the existence and uniqueness of the MC+ path. An R package "plus" has been released.

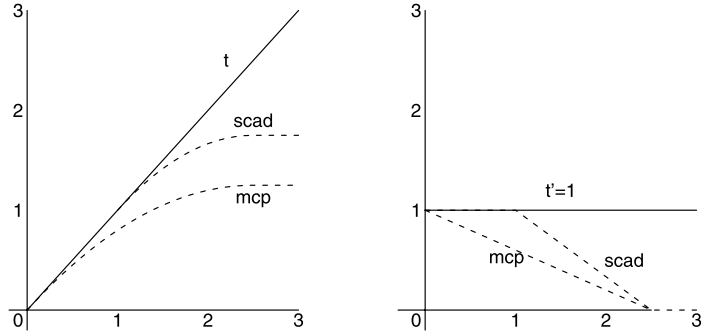

Fig. 1. *The $\ell_1$ penalty $\rho_1(t) = t$ for the LASSO along with the MCP $\rho_2(t)$ and the SCAD penalty $\rho_3(t)$, $t > 0$, $\gamma = 5/2$. Left: penalties $\rho_m(t)$. Right: their derivatives $\dot{\rho}_m(t)$.*

3.1. *Quadratic spline penalties and the MCP.* The PLUS algorithm assumes that the penalty function is of the form $\rho(t; \lambda) = \lambda^2 \rho(t/\lambda)$, where $\rho(t)$ is a nondecreasing quadratic spline in $[0, \infty)$. Such $\rho(t)$ must have a piecewise linear nonnegative continuous derivative $\dot{\rho}(t)$ for $t \geq 0$, so that the solution graph of (2.6) is piecewise linear. The maximum concavity $\kappa(\rho) \equiv \kappa(\rho; \lambda)$ does not depend on $\lambda$. We index $\rho(t)$ by the number of threshold levels $m$, or equivalently the number of knots in $[0, \infty)$, including zero as a knot. Thus,

$$\rho(t; \lambda) = \lambda^2 \rho_m(t/\lambda), \dot{\rho}_m(t) \equiv (d\rho_m/dt)(t)$$

$$(3.1) \qquad = \sum_{i=1}^{m} (u_i - v_i t) I\{t_i \leq t < t_{i+1}\}$$

with $u_1 = 1$, $v_m = 0$, $t_{m+1} = \infty$ and knots $t_1 = 0 < t_2 < \cdots < t_m = \gamma$, satisfying $u_i - v_i t_{i+1} = u_{i+1} - v_{i+1} t_{i+1} \geq 0$, $1 \leq i < m$.

We set $\dot{\rho}_m(0+) = u_1 = 1$ to match the standardization $\dot{\rho}(0+; \lambda) = \lambda$ in (2.3), and $v_m = 0$ for the uniform boundedness of $\dot{\rho}(t; \lambda)$. The unbiasedness feature $\lim_{t \to \infty} \dot{\rho}(t; \lambda) = 0$ demands $t_m = \gamma > 0 = u_m = v_m$ and thus $m > 1$, but the PLUS includes the LASSO with $m = 1$. For $\|\mathbf{x}_j\|^2 = n$, $c_{\min}(\mathbf{\Sigma}_A) \leq 1$, so that (2.5) becomes $\kappa(\rho_m) = \max_{i \leq m} v_i < c_* \leq 1$ under (2.11).

The penalty class (3.1) includes the $\ell_1$ penalty with $m = 1$ and $\kappa(\rho_1) = 0$, the MCP with $m = 2$ and $\kappa(\rho_2) = v_1 = 1/\gamma$, and the SCAD penalty with $m = 3$, $v_1 = 0$, $t_2 = 1$ and $\kappa(\rho_3) = v_2 = 1/(\gamma - 1)$. We plot these three penalty functions $\rho_m$, $m = 1, 2, 3$ and their derivatives in Figure 1, with $\gamma = 5/2$ for the MCP and SCAD penalty.

As mentioned in the Introduction, we propose the MCP (2.1) as the default penalty for the PLUS, and thus the acronym MC+. The MCP corresponds to (3.1) with

$$(3.2) \qquad \rho_2(t) = \min\{t - t^2/(2\gamma), \gamma/2\}, \qquad \dot{\rho}_2(t) = (1 - t/\gamma)_+, \qquad t \geq 0.$$

Among spline penalties satisfying (2.3), the MCP has the smallest number of threshold levels $m = 2$. It follows from (2.6) and (3.1) that the piecewise-linear PLUS path makes a turn whenever $|\widehat{\beta}_j(\lambda)/\lambda|$ hits one of the $m$ thresholds for any $j \le p$. From this point of view, MC+ is the simplest for the PLUS to compute except for the LASSO with $m = 1$.

3.2. *Explicit description of the PLUS algorithm.* Let $\widetilde{\mathbf{z}} \equiv \mathbf{X}'\mathbf{y}/n$. For penalty functions of the form $\rho(t; \lambda) = \lambda^2 \rho_m(t/\lambda)$ with the $\rho_m$ in (3.1), the estimating equation (2.6) is equivalent to the following rescaled version:

$$(3.3) \qquad \begin{cases} z_j - \boldsymbol{\chi}_j'\mathbf{b} = \mathrm{sgn}(b_j)\dot{\rho}_m(|b_j|), & b_j \ne 0, \\ |z_j' - \boldsymbol{\chi}_j'\mathbf{b}| \le 1 = \dot{\rho}_m(0+), & b_j = 0, \end{cases}$$

through the scale change $\widetilde{\mathbf{z}}/\lambda \to \mathbf{z}$ and $\boldsymbol{\beta}/\lambda \to \mathbf{b}$, where $\boldsymbol{\chi}_j \equiv \mathbf{X}'\mathbf{x}_j/n$ are the columns of $\boldsymbol{\Sigma} \equiv \mathbf{X}'\mathbf{X}/n$. The solution $\mathbf{b}(\mathbf{z})$ of (3.3) along the ray $\{\widetilde{\mathbf{z}}/\lambda, \lambda > 0\}$ provides the solution of (2.6) with the inverse transformation $\widehat{\boldsymbol{\beta}}(\lambda) = \lambda \mathbf{b}(\widetilde{\mathbf{z}}/\lambda)$.

We shall "plot" the solution $\mathbf{b}(\mathbf{z})$ of (3.3) against $\mathbf{z}$ to allow multiple solutions, instead of directly solving it for a given $\mathbf{z} = \widetilde{\mathbf{z}}/\lambda = \mathbf{X}'\mathbf{y}/(n\lambda)$. In the univariate case $p = 1$, we plot functions in $\mathbb{R}^2$. For $p > 1$, we need to consider $\mathbf{b}$ versus $\mathbf{z}$ in $\mathbb{R}^{2p}$. Let $H = \mathbb{R}^p$, $H^*$ be its dual, and $\mathbf{z} \oplus \mathbf{b}$ be members of $H \oplus H^* = \mathbb{R}^{2p}$. Define

$$(3.4) \quad u(i) \equiv u_{|i|}, \qquad v(i) \equiv v_{|i|}, \qquad t(i) \equiv \begin{cases} t_i, & 0 < i \le m + 1, \\ -t_{|i|+1}, & -m \le i \le 0, \end{cases}$$

where $u_i, v_i$ and $t_i$ specify $\rho_m$ as in (3.1). For indicators $\boldsymbol{\eta} \in \{-m, \ldots, m\}^p$, let

$$S(\boldsymbol{\eta}) \equiv \text{ the set of all } \mathbf{z} \oplus \mathbf{b}$$

$$(3.5)$$

$$\text{satisfying } \begin{cases} z_j - \boldsymbol{\chi}_j'\mathbf{b} = \mathrm{sgn}(\eta_j)u(\eta_j) - b_j v(\eta_j), & \eta_j \ne 0, \\ -1 \le z_j - \boldsymbol{\chi}_j'\mathbf{b} \le 1, & \eta_j = 0, \\ t(\eta_j) \le b_j \le t(\eta_j + 1), & \eta_j \ne 0, \\ b_j = 0, & \eta_j = 0. \end{cases}$$

Since $\mathrm{sgn}(b_j)\dot{\rho}_m(|b_j|) = \mathrm{sgn}(\eta_j)u(\eta_j) - b_j v(\eta_j)$ for $t(\eta_j) \le b_j \le t(\eta_j + 1)$, (3.3) holds iff (3.5) holds for a certain $\boldsymbol{\eta}$. For each $\boldsymbol{\eta}$, the linear system in (3.5) is of rank $2p$, since one can always uniquely solve for $\mathbf{b}$ and then $\mathbf{z}$ if the inequalities are replaced by equations. Thus, since (3.5) has $p$ equations and $p$ pairs of parallel inequalities, $S(\boldsymbol{\eta})$ are $p$-dimensional parallelepipeds living in $H \oplus H^* = \mathbb{R}^{2p}$. Due to the continuity of $\dot{\rho}_m(t) = (d/dt)\rho_m(t)$ in $t$ by (3.1) and that of $z_j - \boldsymbol{\chi}_j'\mathbf{b}$ in both $\mathbf{z}$ and $\mathbf{b}$, the solutions of (3.5) are identical in the intersection of any given pair of $S(\boldsymbol{\eta})$ with adjacent $\boldsymbol{\eta}$. Furthermore, the $p$-dimensional interiors of different $S(\boldsymbol{\eta})$ are disjoint in view of the constraints of (3.5) on $\mathbf{b}$. Thus, the union of all the $p$-parallelepipeds $S(\boldsymbol{\eta})$ forms a continuous $p$-dimensional surface $S \equiv \cup\{S(\boldsymbol{\eta}) : \boldsymbol{\eta} \in \{-m, \ldots, m\}^p\}$ in $H \oplus H^* = \mathbb{R}^{2p}$. This continuous surface $S$ is the solution set (or the "plot") of all $\mathbf{z} \oplus \mathbf{b} \in H \oplus H^*$ satisfying the rescaled estimating equation (3.3).

Given $\widetilde{\mathbf{z}} = \mathbf{X}'\mathbf{y}/n$, the solution set of (3.3) for all $\mathbf{z} = \tau\widetilde{\mathbf{z}}$ and $\tau > 0$, or equivalently that of (2.6) for all $\lambda$, is identical to the intersection of the surface $S$ and the $(p+1)$-dimensional open half subspace $\{(\tau\widetilde{\mathbf{z}}) \oplus \mathbf{b} : \tau > 0, \mathbf{b} \in H^*\}$ in $\mathbb{R}^{2p}$. Figure 2 depicts the MC+ and LASSO solution sets and the projections of $S(\boldsymbol{\eta})$ to $H$ in the nonoverlapping scenario [under the convexity condition (2.5) with full rank $d^* = p = 2$]. Figure 3 depicts an overlapping scenario in which the complete solution set of (2.6) contains the main branch covered by the MC+ path and a loop not covered.

The rescaled PLUS path in $H \oplus H^*$ is a union of connected line segments

$$(3.6) \qquad \bigcup_{k=0}^{k^*} \ell(\boldsymbol{\eta}^{(k)} | \widetilde{\mathbf{z}}), \qquad \ell(\boldsymbol{\eta} | \mathbf{z}) \equiv S(\boldsymbol{\eta}) \cap \{(\tau\mathbf{z}) \oplus \mathbf{b} : \tau > 0, \mathbf{b} \in H^*\},$$

beginning with $\ell(\boldsymbol{\eta}^{(0)} | \widetilde{\mathbf{z}}) = \{(\tau\widetilde{\mathbf{z}}) \oplus \mathbf{0} : 0 < \tau \le \tau^{(0)}\}$, $\boldsymbol{\eta}^{(0)} = \mathbf{0}$ and connected at

$$(3.7) \quad \{(\tau^{(k-1)}\widetilde{\mathbf{z}}) \oplus \mathbf{b}^{(k-1)}\} = \ell(\boldsymbol{\eta}^{(k-1)} | \widetilde{\mathbf{z}}) \cap \ell(\boldsymbol{\eta}^{(k)} | \widetilde{\mathbf{z}}), \qquad \widetilde{\mathbf{z}} \equiv \mathbf{X}'\mathbf{y}/n.$$

Given $(\tau^{(k-1)}\widetilde{\mathbf{z}}) \oplus \mathbf{b}^{(k-1)}$, we find a *new line segment* $\ell(\boldsymbol{\eta}^{(k)} | \widetilde{\mathbf{z}})$ and compute the other end of it as $(\tau^{(k)}\widetilde{\mathbf{z}}) \oplus \mathbf{b}^{(k)}$, $k \ge 1$. Given $\widetilde{\mathbf{z}}$, we write the turning points in the simpler form $\tau^{(k)} \oplus \mathbf{b}^{(k)} \in \mathbb{R}^{1+p}$. The PLUS path (2.7) is defined through the linear interpolation of $\tau^{(k)} \oplus \mathbf{b}^{(k)}$ and reverse scale change from

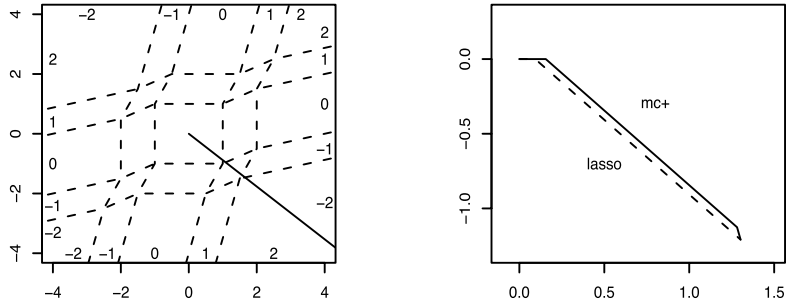

Fig. 2. *Left: the solid ray as $\tau\widetilde{\mathbf{z}}$ and the projections of the $5^2 = 25$ parallelograms $S(\boldsymbol{\eta})$ for the MCP to the $\mathbf{z}$-space $H$ with dashed-edges, labeled by $\eta_1$ and $\eta_2$ along the margins inside the box. Right: the MC+ path (solid) as the entire solution set of (2.6) in the $\boldsymbol{\beta}$-space, along with the LASSO path (dashed). Data: $\|\mathbf{x}_j\|^2/2 = 1$, $\mathbf{x}_1'\mathbf{x}_2/2 = 1/4$, $(\widetilde{z}_1, \widetilde{z}_2) = (1, -0.883)$ and $p = \gamma = 2$.*

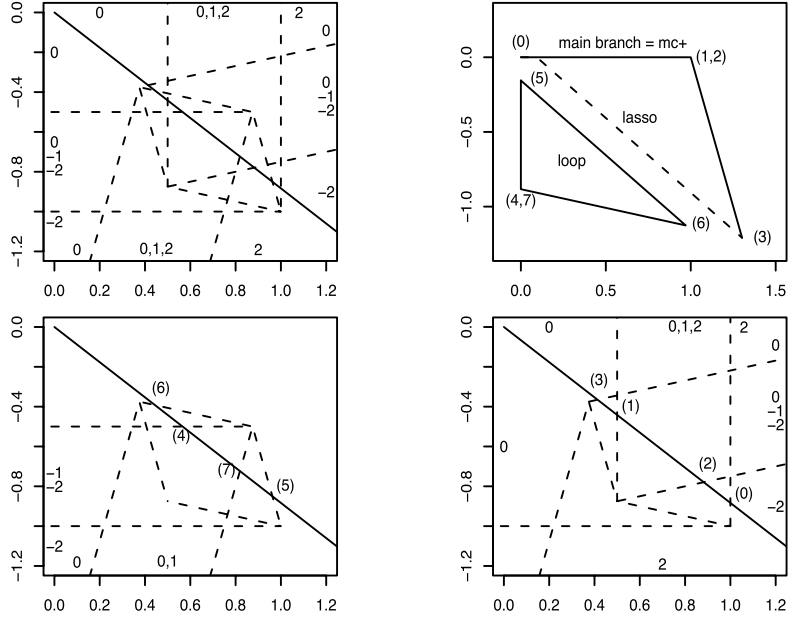

Fig. 3. *Plots for the same data as in Figure 2 with $\gamma = 1/2$ for the MCP. Clockwise from the top left: the $\mathbf{z}$-space plot with overlapping areas marked by multiple values of $\eta_j$; the main branch and one loop as the entire MCP solution set of (2.6) in the $\boldsymbol{\beta}$-space, along with the LASSO; the segments of the main branch with $\tau^{(k)}\widetilde{\mathbf{z}}$, $k = 0, 1, 2, 3$, representing transitions $\boldsymbol{\eta} = \binom{0}{0} \to \binom{1}{0} \to \binom{2}{0} \to \binom{2}{-2} \to \binom{2}{-2}$; the loop with $\tau^{(k)}\widetilde{\mathbf{z}}$, $k = 4, 5, 6, 7$, representing transitions $\boldsymbol{\eta} = \binom{0}{-2} \to \binom{0}{-1} \to \binom{1}{-1} \to \binom{1}{-2} \to \binom{0}{-2}$. For $\boldsymbol{\eta} \in \{-2, 0, 2\}^p$, $\mathbf{z}$-segments turn into $\boldsymbol{\beta}$-points in the MC+ path. A topologically equivalent way of creating the main branch and loop is to fold a piece of paper twice parallel to the horizontal axis and then twice parallel to the vertical axis, cut through the fold and then unfold.*

$\tau \oplus \mathbf{b}$ to $\lambda \oplus \boldsymbol{\beta}$:

$$(3.8) \quad \begin{cases} \tau^{(x)} \oplus \mathbf{b}^{(x)} \equiv (k - x)(\tau^{(k-1)} \oplus \mathbf{b}^{(k-1)}) + (x - k + 1)(\tau^{(k)} \oplus \mathbf{b}^{(k)}), \\ \qquad\qquad\qquad\qquad\qquad\qquad\qquad\qquad\qquad\quad k - 1 < x \leq k, \\ \lambda^{(x)} \oplus \widehat{\boldsymbol{\beta}}^{(x)} \equiv (1 \oplus \mathbf{b}^{(x)})/\tau^{(x)}, \qquad 0 \leq x \leq k^*, \end{cases}$$

with the initialization $\boldsymbol{\eta}^{(0)} = \mathbf{b}^{(0)} = \mathbf{0}$ and $\tau^{(0)} = 1/\max_{j \leq p} |\widetilde{z}_j|$. The PLUS path ends at step $k^*$ if $\widehat{\boldsymbol{\beta}}^{(k^*)}$ provides a global least squares fit with $\mathbf{X}'(\mathbf{y} - \mathbf{X}\widehat{\boldsymbol{\beta}}^{(k^*)}) = \mathbf{0}$. We define $\widehat{\boldsymbol{\beta}}^{(x)} \equiv \widehat{\boldsymbol{\beta}}^{(k^*)}$ and $\lambda^{(x)} = (k^*/x)\lambda^{(k^*)}$ for $x > k^*$. Clearly, $x$ interpolates the number of steps in $[0, k^*]$.

We compute the turning points $\tau^{(k)} \oplus \mathbf{b}^{(k)}$ in (3.8) by finding the "state" $\boldsymbol{\eta}^{(k)}$, the slope $\mathbf{s}^{(k)} \equiv (d\mathbf{b}^{(x)}/d\tau^{(x)})$, $k - 1 < x < k$, the sign $\xi^{(k)} \equiv \text{sgn}(\tau^{(k)} - \tau^{(k-1)})$ and the "length" $\Delta^{(k)} \equiv |\tau^{(k)} - \tau^{(k-1)}|$ for the new segment. We now provide algebraic formulas for the computation of these quantities in a certain "one-at-a-time" scenario. We prove that the PLUS path is one-at-a-time almost everywhere in $(\mathbf{X}, \mathbf{y})$ in the next subsection.

At $\tau^{(k-1)} \oplus \mathbf{b}^{(k-1)}$, (3.8) must hit one of the inequalities in (3.5) for a certain index

$$
\begin{aligned}
(3.9) \qquad & j^{(k-1)} \in \{j : |b_j^{(k-1)}| \in \{t_1, \ldots, t_m\} \\
& \qquad \text{with } \eta_j^{(k-1)} \neq 0, \text{ or } |\tau^{(k-1)} \widetilde{z}_j - \boldsymbol{\chi}_j' \mathbf{b}^{(k-1)}| = 1\},
\end{aligned}
$$

where $t_1, \ldots, t_m$ are the knots of (3.1). If $j^{(k-1)}$ is unique, $\eta_j^{(k)} = \eta_j^{(k-1)}$ for $j \neq j^{(k-1)}$ and

$$
(3.10) \qquad \eta_j^{(k)} = \begin{cases} \operatorname{sgn}(\tau^{(k-1)} \widetilde{z}_j - \boldsymbol{\chi}_j' \mathbf{b}^{(k-1)}), & \eta_j^{(k-1)} = 0, \\ \eta_j^{(k-1)} + \operatorname{sgn}(b_j^{(k-1)} - b_j^{(k-2)}), & \eta_j^{(k-1)} \neq 0, \end{cases}
$$

for $j = j^{(k-1)}$. Let $\boldsymbol{\Sigma}_A$ be as in (2.4) and $A(\boldsymbol{\eta}) \equiv \{j : \eta_j \neq 0\}$. Define

$$
\begin{aligned}
(3.11) \qquad & \boldsymbol{\Sigma}(\boldsymbol{\eta}) \equiv \boldsymbol{\Sigma}_{A(\boldsymbol{\eta})}, \qquad \mathbf{Q}(\boldsymbol{\eta}) \equiv \boldsymbol{\Sigma}(\boldsymbol{\eta}) - \operatorname{diag}(v(\eta_j), \eta_j \neq 0), \\
& d(\boldsymbol{\eta}) \equiv |A(\boldsymbol{\eta})|.
\end{aligned}
$$

Since the $\boldsymbol{\chi}_j$ in (3.3) are the columns of $\boldsymbol{\Sigma}$, for $k - 1 < x < k$ the first equation of (3.5) can be written as $\mathbf{Q}(\boldsymbol{\eta}^{(k)}) \mathbf{P}(\boldsymbol{\eta}^{(k)}) \mathbf{b}^{(x)} = \mathbf{P}(\boldsymbol{\eta}^{(k)})(\tau^{(x)} \widetilde{\mathbf{z}} - \operatorname{sgn}(\boldsymbol{\eta}^{(k)}) u(\boldsymbol{\eta}^{(k)}))$, where $\mathbf{P}(\boldsymbol{\eta})$ is the projection $\mathbf{b} \to (b_j, \eta_j \neq 0)'$ and $u(\cdot)$ is as in (3.4). Differentiating this identity, we find

$$
(3.12) \quad \mathbf{Q}(\boldsymbol{\eta}^{(k)}) \mathbf{P}(\boldsymbol{\eta}^{(k)}) \mathbf{s}^{(k)} = \mathbf{P}(\boldsymbol{\eta}^{(k)}) \widetilde{\mathbf{z}}, \qquad \eta_j^{(k)} = 0 \quad \Rightarrow \quad s_j^{(k)} = 0,
$$

so that $\mathbf{s}^{(k)}$ is solved by inverting $\mathbf{Q}(\boldsymbol{\eta}^{(k)})$. If the segment $\ell(\boldsymbol{\eta}^{(k)} | \widetilde{\mathbf{z}})$ does not live in the boundary of $S(\boldsymbol{\eta}^{(k)})$, the path has to move into its interior from side $j^{(k-1)}$, so that

$$
(3.13) \qquad \xi^{(k)} = \begin{cases} (\eta_j^{(k)} - \eta_j^{(k-1)}) \operatorname{sgn}(s_j^{(k)}), & \eta_j^{(k)} \neq 0, j = j^{(k-1)}, \\ \eta_j^{(k-1)} \operatorname{sgn}(\boldsymbol{\chi}_j' \mathbf{s}^{(k)} - \widetilde{z}_j), & \eta_j^{(k)} = 0, j = j^{(k-1)}. \end{cases}
$$

Given the slope $\mathbf{s}^{(k)}$ and the sign $\xi^{(k)}$ of $d\tau$ for the segment, there are at most $p$ possible ways for $(\tau \widetilde{\mathbf{z}}) \oplus \mathbf{b}(\tau \widetilde{\mathbf{z}})$ to hit a new side of the boundary of the $p$-parallelepiped $S(\boldsymbol{\eta}^{(k)})$ in (3.5). If it first hits the boundary indexed by

$\eta_j^{(k)}$, by (3.5) and (3.8) $\Delta^{(k)}$ would be

$$(3.14) \quad \Delta_j^{(k)} = \begin{cases} \xi_j^{(k)}\{t(\eta_j^{(k)}+1)-b_j^{(k-1)}\}/s_j^{(k)}, \\ \qquad \xi_j^{(k)}s_j^{(k)} > 0 \neq \eta_j^{(k)}, \\ \xi_j^{(k)}\{t(\eta_j^{(k)})-b_j^{(k-1)}\}/s_j^{(k)}, \\ \qquad \xi_j^{(k)}s_j^{(k)} < 0 \neq \eta_j^{(k)}, \\ \xi_j^{(k)}\{1-g_j^{(k-1)}\}/\{\widetilde{z}_j - \boldsymbol{\chi}_j'\mathbf{s}^{(k)}\}, \\ \qquad \xi_j^{(k)}(\widetilde{z}_j - \boldsymbol{\chi}_j'\mathbf{s}^{(k)}) > 0 = \eta_j^{(k)}, \\ \xi_j^{(k)}\{-1-g_j^{(k-1)}\}/\{\widetilde{z}_j - \boldsymbol{\chi}_j'\mathbf{s}^{(k)}\}, \\ \qquad \xi_j^{(k)}(\widetilde{z}_j - \boldsymbol{\chi}_j'\mathbf{s}^{(k)}) < 0 = \eta_j^{(k)}, \end{cases}$$

where $t(\cdot)$ is as in (3.4) and $g_j^{(k-1)} \equiv \tau^{(k-1)}\widetilde{z}_j - \boldsymbol{\chi}_j'\mathbf{b}^{(k-1)}$. It follows that

$$(3.15) \qquad \tau^{(k)} = \tau^{(k-1)} + \xi^{(k)}\Delta^{(k)}, \qquad \Delta^{(k)} = \min_{1 \leq j \leq p} \Delta_j^{(k)},$$

with the minimum attained at $j = j^{(k)}$ as in (3.9). We formally write the PLUS as follows.

THE PLUS ALGORITHM.
*Initialization*: $\boldsymbol{\eta}^{(0)} \leftarrow \mathbf{0}$, $\mathbf{b}^{(0)} \leftarrow \mathbf{0}$, $\tau^{(0)} \leftarrow 1/\max_{j \leq p}|\widetilde{z}_j|$, $k \leftarrow 1$.
*Iteration*:

(3.16)          Find $\boldsymbol{\eta}^{(k)}$ by (3.9) and (3.10),

(3.17)          Find $\mathbf{s}^{(k)}$ by (3.12),

(3.18)          Find $\tau^{(k)}$ by (3.13), (3.14) and (3.15),

(3.19)          $\mathbf{b}^{(k)} \leftarrow \mathbf{b}^{(k-1)} + (\tau^{(k)}-\tau^{(k-1)})\mathbf{s}^{(k)}$,

                $k \leftarrow k+1$.

*Termination*: (3.16) has no solution for $k = k^*+1$ or $\tau^{(k^*)} = \infty$.
*Output*: $\tau^{(0)}$, $\mathbf{b}^{(0)}$, $\boldsymbol{\eta}^{(k)}$, $\mathbf{s}^{(k)}$, $\tau^{(k)}$, $\mathbf{b}^{(k)}$, $k = 1, 2, \ldots, k^*$.

3.3. *The existence and uniqueness of the PLUS path.* We prove in this subsection that for the MCP the PLUS algorithm computes the main branch (2.7) of the solution graph of (2.6) and that the main branch is unique almost everywhere in $(\mathbf{X}, \mathbf{y})$.

*Nondegenerate designs.* The design matrix $\mathbf{X}$ in (1.2) is nondegenerate if for all $A \subset \{1, \ldots, p\}$ of size $|A| = n \wedge p - 1$ and $\eta_j \in \{-1, 0, 1\}$, $j \leq p$, the

$n \wedge p$ *vectors*

(3.20) $$\left\{ \mathbf{x}_j, j \in A, \sum_{k \notin A} \eta_k \mathbf{x}_k \right\} \text{ are linearly independent.}$$

*For $p \leq n$, $\mathbf{X}$ is nondegenerate iff $\operatorname{rank}(\mathbf{X}) = p$.*

THEOREM 3. *Suppose the MCP is used in the PLUS algorithm. Let $\mathbf{Q}(\boldsymbol{\eta}^{(k)})$ be as in (3.12).*

(i) *Suppose the design matrix $\mathbf{X}$ is nondegenerate in the sense of (3.20). Given $\mathbf{X}$, there exists a finite set $\Gamma_0(\mathbf{X})$ such that for all $\gamma \notin \Gamma_0(\mathbf{X})$, a path of the form (3.8) exists with $\det(\mathbf{Q}(\boldsymbol{\eta}^{(k)})) \neq 0$ for $k \leq k^*$ and perfect fit $\mathbf{X}'(\mathbf{y} - \mathbf{X}\widehat{\boldsymbol{\beta}}^{(k^*)}) = \mathbf{0}$ at a finite final step $k^*$.*

(ii) *For fixed $\gamma > 0$, the design matrix $\mathbf{X}$ is nondegenerate and $\gamma \notin \Gamma_0(\mathbf{X})$ almost everywhere in $\mathbb{R}^{n \times p}$ under the Lebesgue measure.*

(iii) *For fixed positive $\gamma \neq 1$, the design matrix $\mathbf{X}$ is nondegenerate and $\gamma \notin \Gamma_0(\mathbf{X})$ almost everywhere under the product of $p$ Haar measures in the $(n-1)$-sphere $\{\mathbf{x} : \|\mathbf{x}\|^2 = n\}$.*

(iv) *Suppose $\gamma \notin \Gamma_0(\mathbf{X})$. Then, almost everywhere in $\widetilde{\mathbf{z}} = \mathbf{X}'\mathbf{y}/n \in \mathbb{R}^p$, the graph of (2.7) is unique and the PLUS algorithm computes (2.7) within a finite step $k^*$ and ends with an optimal fit satisfying $\mathbf{X}'(\mathbf{y} - \mathbf{X}\widehat{\boldsymbol{\beta}}^{(k^*)}) = \mathbf{0}$. Consequently, for all $0 \leq k \leq k^*$ the path (3.8) is one-at-a-time in the sense of* (a) *the uniqueness and validity of (3.9), (3.10), (3.12), (3.13) and (3.15) and* (b) *the positiveness of $\Delta^{(k)}$ and $\tau^{(k)}$ in (3.15).*

(v) *If $\mathbf{Q}(\boldsymbol{\eta}^{(k)})$ is positive-definite and $\ell(\boldsymbol{\eta}^{(k)}|\widetilde{\mathbf{z}})$ in (3.6) does not live in the boundary of $S(\boldsymbol{\eta}^{(k)})$ in (3.5), then $\widehat{\boldsymbol{\beta}}^{(x)}$ is a local minimizer of $L(\mathbf{b}; \lambda)$ in (1.1) at $\lambda = \lambda^{(x)}$, $k - 1 < x < k$.*

Theorem 3(ii) and (iii) ensure that $\gamma \notin \Gamma_0(\mathbf{X})$ almost everywhere in $\mathbf{X}$ for all fixed $\{n, p, \gamma\}$. The condition of $\gamma \notin \Gamma_0(\mathbf{X})$ is not necessary for the MC+ path to end with an optimal fit. For example, if $\mathbf{x}_{j_0} = \pm\mathbf{x}_{k_0}$, the PLUS path uses at most one design vector $\mathbf{x}_{j_0}$ or $\mathbf{x}_{k_0}$ in any step, so that it behaves as if one of them never exists. Theorem 3(iv) guarantees that the PLUS algorithm yields an entire path of solutions (2.7) covering all $0 \leq \lambda < \infty$. Theorem 3(v) implies that the estimator $\widehat{\boldsymbol{\beta}}(\lambda)$ is a local minimizer under (2.5) whenever $\#\{j : \widehat{\beta}_j(\lambda) \neq 0\} \leq d^*$, as guaranteed by the conditions of Theorems 1, 2, 5 and 6. For simplicity, we omit an extension of Theorem 3 to the PLUS with general quadratic penalty (3.1).

We note that the map $\lambda^{(x)} \to \widehat{\boldsymbol{\beta}}^{(x)}$ is potentially many-to-one in the PLUS path due to the possible concavity of the penalized loss, since $\tau^{(k)} < \tau^{(k-1)}$ is allowed as (3.8) traverses through the solution graph. Theorem 3 does not guarantee that the PLUS path contains all solutions of (2.6) due to loops outside its path, as Figure 3 demonstrates. However, such multiplicity of branches is less severe for sparse data. In the example in Figure 4, the convex penalized loss with $\gamma = 2$ yields identical MC+ path as the nonconvex one with $\gamma = 1/2$ for sparse data outside regions where the the projections of the parallelograms $S(\boldsymbol{\eta})$ fold severely in the **z**-space for $\gamma = 1/2$. This should be compared with the dramatic difference between $\gamma = 2$ and $\gamma = 1/2$ in Figures 2 and 3 for dense data.

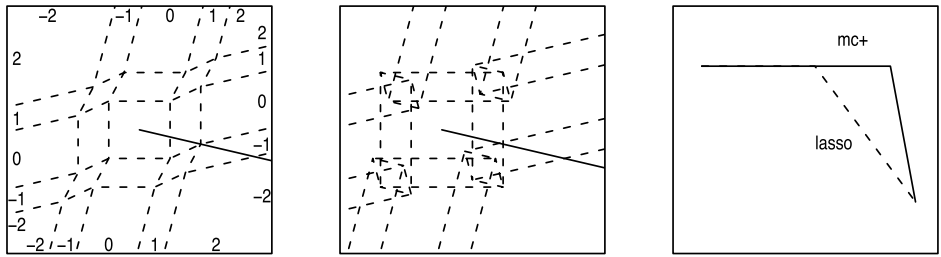

FIG. 4.   *The same type of plots as in Figures 2 and 3 for the same* **X** *and more sparse* $(\tilde{z}_1, \tilde{z}_2) = (1, -1/2)$. *From the left: the* **z**-*space plot for MC+ with* $\gamma = 2$; *MC+ with* $\gamma = 1/2$; *the MC+ (same for both* $\gamma = 2$ *and* $\gamma = 1/2$) *and LASSO paths in the* $\boldsymbol{\beta}$-*space. The loop disappears since the solid line* $\tau \tilde{\mathbf{z}}$ *does not pass through the places where the projection of S folds in two different directions.*

**4. Selection consistency for general penalty.** We provide in Section 4.1 two sets of lower bounds for the probability of correct selection for general penalized LSE: one for the global minimizer of (1.1) in the regular case of rank(**X**) = $p$ ($p \le n$ necessarily) and one for the local solution (2.8) in the case of rank(**X**) < $p$ (including $p \gg n$). These lower bounds imply the sign consistency $P\{\text{sgn}(\widehat{\boldsymbol{\beta}}) = \text{sgn}(\boldsymbol{\beta})\} \to 1$ and thus the selection consistency (1.4) as $\max(n, p) \to \infty$. As a crucial element of our proof and a matter of independent interest, we also provide in Section 4.2 upper bounds of the false positive for any given oracular set $B$ of interest and a general class of penalties.

4.1. *Probability bounds for selection consistency.* Our selection consistency results are proved by showing that the global minimizer of (1.1) or the local solution (2.8) are identical to the oracle LSE (2.10) with high probability. Let

(4.1)            $(w_j^o, j \in A^o)' =$ the diagonal elements of $\boldsymbol{\Sigma}_{A^o}^{-1}$,

so that $\text{Var}(\widehat{\beta}_j^o) = w_j^o \sigma^2 / n \, \forall j \in A^o$ for the oracle LSE $\widehat{\boldsymbol{\beta}}^o$ with $B = A^o$. We first present nonasymptotic bounds for selection consistency under the fol-lowing global convexity condition:

$$(4.2) \quad c_{\min}(\mathbf{\Sigma}) + \{\dot{\rho}(t_2; \lambda) - \dot{\rho}(t_1; \lambda)\}/(t_2 - t_1) > 0 \qquad \forall 0 < t_1 < t_2,$$

where $\mathbf{\Sigma} \equiv \mathbf{X}'\mathbf{X}/n$. Under (4.2), (2.6) is the KKT condition and its solution is unique, so that the estimator (2.8) is the global minimizer of (1.1). Let $\Phi(\cdot)$ be the $N(0,1)$ distribution.

THEOREM 4. *Suppose (2.3) and (4.2) hold for* $\lambda_1 \leq \lambda \leq \lambda_2$. *Let* $\widehat{\boldsymbol{\beta}}(\lambda)$ *be as in (2.8) for each* $\lambda > 0$ *and* $\widehat{\boldsymbol{\beta}} = \widehat{\boldsymbol{\beta}}(\widehat{\lambda})$ *for a deterministic or random penalty level* $\widehat{\lambda}$. *Let* $A^o$, $d^o$, $\widehat{A}$ *and* $\beta_* \equiv \min_{\beta_j \neq 0} |\beta_j|$ *be as in (1.3), (1.4) and (2.9) and* $\widehat{\boldsymbol{\beta}}^o$ *be as in (2.10) with* $B = A^o$. *Suppose* $\beta_* \geq \gamma \lambda_2$ *and* $P\{\lambda_1 \leq \widehat{\lambda} \leq \lambda_2\} = 1$. *Then*

$$(4.3) \quad P\{\widehat{A} \neq A^o\} \leq P\{\widehat{\boldsymbol{\beta}} \neq \widehat{\boldsymbol{\beta}}^o \ or \ \mathrm{sgn}(\widehat{\boldsymbol{\beta}}) \neq \mathrm{sgn}(\boldsymbol{\beta})\} \leq \pi_{n,1}(\lambda_1) + \pi_{n,2}(\lambda_2),$$

*where* $\pi_{n,1}(\lambda) \equiv 2\sum_{j \notin A^o} \Phi(-n\lambda/(\sigma \|\mathbf{x}_j\|))$ *and* $\pi_{n,2}(\lambda) \equiv \sum_{j \in A^o} \Phi((\gamma\lambda - |\beta_j|)/(\sigma(w_j^o/n)^{1/2}))$.

COROLLARY 2. *Suppose (2.3), (4.2),* $\|\mathbf{x}_j\|^2 = n$ *and* $|\beta_j| \geq \gamma\lambda + \sigma \times \sqrt{w_j^o(2/n)\log a_n}$ *for all* $j \in A^o$ *with* $a_n \geq d^o$ *and* $\lambda \geq \lambda_{1,1} \equiv \sigma \times \sqrt{(2/n)\log(p - d^o)}$. *Then, for large* $\sqrt{n}\lambda/\sigma$ *and* $a_n$,

$$(4.4) \quad P\{\widehat{\boldsymbol{\beta}}(\lambda) \neq \widehat{\boldsymbol{\beta}}^o \ or \ \mathrm{sgn}(\widehat{\boldsymbol{\beta}}(\lambda)) \neq \mathrm{sgn}(\boldsymbol{\beta})\} \to 0.$$

For the MC+, (4.2) is equivalent to $c_{\min}(\mathbf{\Sigma}) > 1/\gamma$, and $\beta_* \geq (\gamma + \sqrt{w^o})\lambda_{\mathrm{univ}}$ with $p \to \infty$ suffices for (4.4), where $w^o \equiv \max_{j \in A^o} w_j^o$ is as in Theorem 1. For the SCAD, we need the larger $\gamma > 1 + 1/c_{\min}(\mathbf{\Sigma})$ for (4.2). For $d^o \ll p$ and $\|\mathbf{x}_j\|^2 = n$, (4.4) provides theoretical support to the heuristic condition (2.9) for the selection consistency at $\lambda = \lambda_{\mathrm{univ}}$.

We now consider selection consistency for general $p$, including $p \gg n$. For $c^* \geq c_* \geq \kappa \geq 0$ and $0 < \alpha < 1$, define $w \equiv w_{c_*, c^*, \kappa, \alpha} \equiv (2-\alpha)/(c_* c^*/\kappa^2 - 1)$ and

$$(4.5) \quad \begin{aligned} K_* &\equiv K_{c_*, c^*, \kappa, \alpha} \\ &\equiv \inf_{0 < t < (2/w + 1 + \alpha)/\alpha} \frac{(1 + w\{1 + (\alpha/t)/(1-\alpha)\})c^*/c_* - 1}{\{2 + w(1 + \alpha - t\alpha)\}(1-\alpha)}. \end{aligned}$$

THEOREM 5. *Let* $\rho(t; \lambda)$ *be a penalty satisfying* $\dot{\rho}(0+; \lambda) = \lambda$, $\dot{\rho}(t; \lambda) \leq \lambda I\{t \leq \gamma\lambda\}$ *and* $\ddot{\rho}(t; \lambda) \geq -\kappa$ *for all* $t > 0$ *and* $\lambda \geq \lambda_1$. *Let* $A^o$, $d^o$, $\widehat{\boldsymbol{\beta}} = \widehat{\boldsymbol{\beta}}(\widehat{\lambda})$, $\widehat{A}$, $\beta_*$, $\widehat{\boldsymbol{\beta}}^o$, $w_j^o$, $\pi_{n,1}(\lambda)$ *and* $\pi_{n,2}(\lambda)$ *be as in Theorem 4. Suppose (2.11) holds* with certain rank $d^*$ and $c^* \geq c_* > \kappa$. For these $\{c_*, c^*, \kappa\}$ and $0 < \alpha < 1$, let $K_*$ be as in (4.5). Suppose (1.2) holds with $d^o \leq d_* = d^*/(1 + K_*)$. Let $\pi_{n,3}(\lambda) \equiv \binom{p-d^o}{m} P\{\sigma^2 \chi_m^2 > m\lambda\}$ with $m = d^* - d^o$.

TABLE 2
*Example configurations of $\{c_*, c^*, \kappa, \alpha\}$ for fixed $K_*$ $c^* = 1 + \delta$, $c_* = 1 - \delta$, optimal $t = \sqrt{(c^*/c_*)/K_*}/(1-\alpha)$ in (4.5)*

| $K_* = 1/2$ | | | $K_* = 1$ | | | $K_* = 2$ | | | $K_* = 3$ | | |
|---|---|---|---|---|---|---|---|---|---|---|---|
| $\delta$ | $\alpha$ | $1/\kappa \geq$ | $\delta$ | $\alpha$ | $1/\kappa \geq$ | $\delta$ | $\alpha$ | $1/\kappa \geq$ | $\delta$ | $\alpha$ | $1/\kappa \geq$ |
| 1/4 | 1/5 | 4.84 | 2/5 | 1/5 | 4.14 | 1/2 | 1/3 | 3.30 | 1/2 | 1/2 | 2.98 |
| 1/5 | 2/7 | 3.73 | 1/3 | 1/3 | 3.57 | 1/3 | 1/2 | 2.32 | 1/3 | 1/2 | 1.73 |
| 1/6 | 1/3 | 3.28 | 1/4 | 2/5 | 2.65 | 1/4 | 1/2 | 1.86 | 1/4 | 1/2 | 1.49 |

(i) Let $\lambda_2 \geq \max\{\lambda_1, (\sqrt{c^*}/\alpha)\lambda_3\}$. Suppose $\beta_* \geq \gamma\lambda_2$ and $P\{\lambda_1 \leq \widehat{\lambda} \leq \lambda_2\} = 1$. Then

$$(4.6) \quad P\{\widehat{A} \neq A^o\} \leq P\{\widehat{\boldsymbol{\beta}} \neq \widehat{\boldsymbol{\beta}}^o \text{ or } \text{sgn}(\widehat{\boldsymbol{\beta}}) \neq \text{sgn}(\boldsymbol{\beta})\} \leq \sum_{k=1}^{3} \pi_{n,k}(\lambda_k).$$

(ii) Let $\lambda_{1,\epsilon} \equiv \sigma\sqrt{(2/n)\log((p-d^o)/\epsilon)}$, $\lambda_{3,\epsilon} = \sigma\sqrt{(2/n)\log\widetilde{p}_\epsilon}$ with $\widetilde{p}_\epsilon$ in (2.12), $\lambda_{2,\epsilon} \geq \max\{\lambda_{1,\epsilon}, (\sqrt{c^*}/\alpha)\lambda_{3,\epsilon}\}$ and $a_n \geq d^o$. Suppose $|\beta_j| \geq \gamma\lambda_{2,\epsilon} + \sigma\sqrt{w_j^o(2/n)\log(a_n/\epsilon)}$ for $j \in A^o$ and $\|\mathbf{x}_j\|^2 = n$. If $P\{\lambda_{1,\epsilon} \leq \widehat{\lambda} \leq \lambda_{2,\epsilon}\} = 1$, then

$$
\begin{aligned}
P\{\widehat{A} \neq A^o\} &\leq P\{\widehat{\boldsymbol{\beta}} \neq \widehat{\boldsymbol{\beta}}^o \text{ or } \text{sgn}(\widehat{\boldsymbol{\beta}}) \neq \text{sgn}(\boldsymbol{\beta})\} \\
(4.7) \qquad &\leq \epsilon\left\{\frac{1}{1 \vee J_1} + \frac{d^o/(2a_n)}{1 \vee J_2} + \frac{(4\log\widetilde{p}_\epsilon)^{-1/2}}{1 \vee J_3}\right\} \\
&\leq \left(\frac{3}{2} + \frac{1}{\sqrt{2}}\right)\epsilon
\end{aligned}
$$

with $J_1 = \sqrt{\pi\log((p-d^o)/\epsilon)}$, $J_3 = \{2\log\widetilde{p}_\epsilon - 1 + 1/m\}\sqrt{m\pi}/\sqrt{4\log\widetilde{p}_\epsilon}$ and $J_2 = \sqrt{\pi\log(a_n/\epsilon)}$. Consequently, (1.4) holds as $\epsilon^{-1} \vee \min(J_1, J_2, J_3) \to \infty$ and $P\{\lambda_1 \leq \widehat{\lambda} \leq \lambda_2\} \to 1$.

REMARK 5. A convenient choice is $\alpha = 1/2$ and $t = 3$ in (4.5) which leads to $K_* \leq \{1 + 2/(c_*c^*/\kappa^2 - 1)\}c^*/c_* - 1$. In Theorems 1 and 2, $1/\kappa = \gamma \geq c_*^{-1}\sqrt{4 + c_*/c^*}$, so that $K_* \leq c^*/c_* - 1/2$. For the LASSO, $\kappa = 0 = w$ and $K_* = (c^*/c_* - 1)/(2 - 2\alpha)$. Some other configurations of $\{c_*, c^*, \kappa, \alpha\}$ are given in Table 2.

REMARK 6. Theorem 5(i) is applicable to the problem of finding a sparse solution $\boldsymbol{\beta}$ of $\mathbf{y} = \mathbf{X}\boldsymbol{\beta}$ with $p > n$, i.e., $\boldsymbol{\varepsilon} = 0$ in (1.2). With $\lambda_2 = \lambda^{(k^*)}$ (nearly zero) and $\sigma = \lambda_1 = \lambda_3 = \alpha = 0$, it asserts $\widehat{\boldsymbol{\beta}}^{(k^*)} = \boldsymbol{\beta}$ at the last step of the PLUS algorithm whenever $\beta_* > \gamma\lambda^{(k^*)}$ and $d^o < d^*/(K_* + 1)$, where $K_* + 1 = (c^*/c_* + 1)/\{2 - \kappa^2/(c^*c_*)\}$. See Section 6.5.

REMARK 7. Consider the MC+ and LASSO. For $\beta_* > \gamma\lambda_{\text{univ}}$, the oracle $\tau(\widetilde{\mathbf{z}} \oplus \widehat{\boldsymbol{\beta}}^o)$ has a high probability of solving (3.5) for the parallelepiped $S(\boldsymbol{\eta})$ with $\boldsymbol{\eta} = 2\,\text{sgn}(\boldsymbol{\beta})$. Such parallelepipeds are unbiased, since they involve regions with $u(\pm 2) = v(\pm 2) = 0 = \dot{\rho}_2(|b_j|)$ in (3.5). An extension of Theorem 5 to biased $S(\boldsymbol{\eta})$ requires $\text{sgn}(\beta_j)(1 + I\{|\beta_j| > \gamma\lambda\}) = \eta_j$ with a larger $\lambda$. Such an extension with $\max_j |\beta_j| < \gamma\lambda$ and $\boldsymbol{\eta} = \text{sgn}(\boldsymbol{\beta})$ would match the theory of selection consistency for the LASSO with uniformity in a neighborhood of $\gamma = \infty$.

Compared with Theorem 4, an obvious advantage of Theorem 5 is its applicability to $p > n$. In the case of $p > n$, Theorem 5 still allows $c_* > c_{\min}(\boldsymbol{\Sigma})$ and thus smaller $\gamma = 1/\kappa$ and $\beta_*$ for the MC+ than Theorem 4 does. With $\kappa = 1/\gamma$, the MCP allows the smallest $\gamma$ and thus the smallest possible $\beta_*$ in Theorem 5.

4.2. *An upper bound for the false positive.* Given a target set $B \subset \{1, \ldots, p\}$, we provide upper bounds for the false positive $\#\{j \notin B : |\widehat{\beta}_j(\lambda)| > 0\}$ for the selector (2.8) with a general class of penalties. See Remark 4 for examples of $B$.

THEOREM 6. *Suppose* (2.11) *holds with certain* $d^*$ *and* $c^* \geq c_* \geq \kappa \geq 0$. *For these* $\{c_*, c^*, \kappa\}$ *and an* $\alpha \in (0, 1)$, *let* $K_*$ *be as in* (4.5). *Let* $B$ *be a deterministic subset of* $\{1, \ldots, p\}$ *with* $|B| = d^o \leq d_* = d^*/(K_* + 1)$. *Let* $\lambda_1 > 0$. *Suppose* $\rho(t; \lambda)$ *satisfy* $\lambda(1 - \kappa t/\lambda)_+ \leq \dot{\rho}(t; \lambda) \leq \lambda$ *for* $t > 0$ *and* $\lambda \geq \lambda_1$. *Let* $\widehat{\lambda} \geq \lambda_1 \vee \{(\sqrt{c^*}/\alpha)(\sigma\sqrt{(2/n)\log\widetilde{p}_\epsilon} + \theta_B/\sqrt{m})\}$ *with the* $\theta_B$ *in Theorem* 2, $m = d^* - d^o$ *and* $\widetilde{p}_\epsilon$ *in* (2.12). *Let* $\widehat{\boldsymbol{\beta}} = \widehat{\boldsymbol{\beta}}(\widehat{\lambda})$ *with the* $\widehat{\boldsymbol{\beta}}(\lambda)$ *in* (2.8). *Then*

(4.8)
$$P\{\#(j \notin B : \widehat{\beta}_j \neq 0) \geq 1 \vee (K_*|B|)\}$$
$$\leq \epsilon(\log\widetilde{p}_\epsilon)^{-1/2}e^{\mu^2/2}\Phi(-\mu) \leq \epsilon/\sqrt{2},$$

*where* $\mu = \{2\log\widetilde{p}_\epsilon - 1 + 1/m\}\sqrt{m}/\sqrt{2\log\widetilde{p}_\epsilon}$ *and* $\Phi(x)$ *is the* $N(0, 1)$ *distribution function.*

This theorem is an extension of the upper bound on $|\widehat{A}|$ in Zhang and Huang (2008) from the LASSO to a general continues path of penalized LSE. Since it is relatively easy to find sharp conditions for the oracle LSE (2.10) to be a solution of (2.6), the upper bounds in Theorem 6 is a crucial element in our proof of selection consistency. Remark 5 applies to Theorem 6.

**5. The MSE, degrees of freedom and noise level.**    In this section, we consider the estimation of the estimation and prediction risks for general penalized LSE and the noise level in (1.2). Formulas for the *degrees of freedom* and unbiased risk estimators are derived and justified via Stein's (1981) unbiased risk estimation (SURE). Necessary and sufficient conditions are provided for the continuity of the penalized LSE.

5.1. *The estimation of MSE and degrees of freedom.*    The formulas derived here are based on Stein's (1981) theorem for the unbiased estimation of the MSE of almost differentiable estimators of a mean vector. A map $\mathbf{h} : \mathbb{R}^p \to \mathbb{R}^p$ is almost differentiable if

$$(5.1) \qquad \mathbf{h}(\mathbf{z} + \mathbf{v}) = \mathbf{h}(\mathbf{z}) + \left\{ \int_0^1 \mathbf{H}(\mathbf{z} + x\mathbf{v}) \, dx \right\} \mathbf{v} \qquad \forall \mathbf{v} \in \mathbb{R}^p,$$

for a certain map $\mathbf{H} : \mathbb{R}^p \to \mathbb{R}^{p \times p}$. Suppose in this subsection that $\rho(t; \lambda)$ is almost twice differentiable in $t > 0$, or equivalently

$$(5.2) \qquad \dot{\rho}(t; \lambda) \equiv \frac{\partial}{\partial t} \rho(t; \lambda) = \dot{\rho}(1; \lambda) + \int_1^t \ddot{\rho}(x; \lambda) \, dx \qquad \forall t > 0,$$

for a certain function $\ddot{\rho}(x; \lambda)$. Under this condition, $\ddot{\rho}(t; \lambda) = (\partial/\partial t)\dot{\rho}(t; \lambda)$ almost everywhere in $(0, \infty)$ and the maximum concavity (2.2) can be written as $\kappa(\rho; \lambda) = \|(\ddot{\rho}(t; \lambda))_-\|_\infty$.

For multivariate normal vectors $\mathbf{z} \sim N(\boldsymbol{\mu}, \mathbf{V})$, Stein's theorem can be stated as

$$(5.3) \qquad E\mathbf{h}(\mathbf{z})(\mathbf{z} - \boldsymbol{\mu})' = E\mathbf{H}(\mathbf{z})\mathbf{V},$$

provided (5.1) and the integrability of all the elements of $\mathbf{H}(\mathbf{z})$. This applies to the penalized LSE. Let $\boldsymbol{\Sigma}_A$ be as in (2.4). We extend (3.11) to general penalties $\rho(t; \lambda)$ as follows:

$$(5.4) \qquad \mathbf{Q}(\boldsymbol{\beta}; \lambda) \equiv \boldsymbol{\Sigma}_{\{j \, : \, \beta_j \neq 0\}} + \mathrm{diag}(\ddot{\rho}(|\beta_j|; \lambda), \beta_j \neq 0),$$

$$d(\boldsymbol{\beta}) \equiv \#\{j : \beta_j \neq 0\}.$$

THEOREM 7.    *Let $\lambda > 0$ be fixed and $\widehat{\boldsymbol{\beta}} \equiv \widehat{\boldsymbol{\beta}}(\lambda) \equiv \arg\min_{\mathbf{b}} L(\mathbf{b}; \lambda)$ with the data $(\mathbf{X}, \mathbf{y})$ in (1.2) and $L(\mathbf{b}; \lambda)$ in (1.1). Suppose (2.5) holds with $d^* = p$.*

Let $\mathbf{\Sigma} \equiv \mathbf{X}'\mathbf{X}/n$ and $\widehat{\mathbf{P}}$ be the $d(\widehat{\boldsymbol{\beta}}) \times p$ matrix giving the projection $\widehat{\mathbf{P}}\mathbf{b} = (b_j : \widehat{\beta}_j \neq 0)'$ as in (3.12). Then

$$
\begin{aligned}
(5.5) \quad & E(\widehat{\boldsymbol{\beta}} - \boldsymbol{\beta})(\widehat{\boldsymbol{\beta}} - \boldsymbol{\beta})' \\
& = E\left\{ (\widehat{\boldsymbol{\beta}} - \widetilde{\boldsymbol{\beta}})(\widehat{\boldsymbol{\beta}} - \widetilde{\boldsymbol{\beta}})' + \frac{2\sigma^2}{n}\widehat{\mathbf{P}}'\mathbf{Q}^{-1}(\widehat{\boldsymbol{\beta}}; \lambda)\widehat{\mathbf{P}} \right\} - \frac{\sigma^2}{n}\mathbf{\Sigma}^{-1},
\end{aligned}
$$

where $\widetilde{\boldsymbol{\beta}} \equiv \mathbf{\Sigma}^{-1}\mathbf{X}'\mathbf{y}/n$ is the ordinary LSE of $\boldsymbol{\beta}$. In particular, for all $\mathbf{a} \in \mathbb{R}^p$,

$$
(5.6) \quad |\mathbf{a}'(\widehat{\boldsymbol{\beta}} - \widetilde{\boldsymbol{\beta}})|^2 + \frac{2\widehat{\sigma}^2}{n}(\widehat{\mathbf{P}}\mathbf{a})'\mathbf{Q}^{-1}(\widehat{\boldsymbol{\beta}}; \lambda)(\widehat{\mathbf{P}}\mathbf{a}) - \frac{\widehat{\sigma}^2}{n}\mathbf{a}'\mathbf{\Sigma}^{-1}\mathbf{a}
$$

is an unbiased estimator of the MSE $E|\mathbf{a}'(\widehat{\boldsymbol{\beta}} - \boldsymbol{\beta})|^2$, provided $\widehat{\sigma}^2 = \sigma^2$ in the case of known $\sigma^2$ or $\widehat{\sigma}^2 = \|\mathbf{y} - \mathbf{X}\widetilde{\boldsymbol{\beta}}\|^2/(n-p)$ in the case of $p < n$. Consequently,

$$
(5.7) \quad E\left\{ \|\widehat{\boldsymbol{\beta}} - \widetilde{\boldsymbol{\beta}}\|^2 + \frac{2\widehat{\sigma}^2}{n}\operatorname{trace}(\mathbf{Q}^{-1}(\widehat{\boldsymbol{\beta}}; \lambda)) - \frac{\widehat{\sigma}^2}{n}\operatorname{trace}(\mathbf{\Sigma}^{-1}) \right\} = E\|\widehat{\boldsymbol{\beta}} - \boldsymbol{\beta}\|^2.
$$

REMARK 8. Condition (2.5) with $d^* = p$ asserts $c_{\min}(\mathbf{\Sigma}) > \kappa(\rho; \lambda)$, which is slightly stronger than the global convexity condition (4.2). We prove in the next subsection that (4.2) is a necessary and sufficient condition for the continuity of $\widehat{\boldsymbol{\beta}}$, which is weaker than the almost differentiability of $\widehat{\boldsymbol{\beta}}$. Thus, the conditions of Theorem 7 are nearly sharp for the application of the SURE. In the $k$th segment of the PLUS path, $\mathbf{Q}(\widehat{\boldsymbol{\beta}}(\lambda); \lambda) = \mathbf{Q}(\boldsymbol{\eta}^{(k)})$ as in (3.11).

Let $\boldsymbol{\mu} \equiv E\mathbf{y} = \mathbf{X}\boldsymbol{\beta}$ and $\widehat{\boldsymbol{\mu}} = \mathbf{X}\widehat{\boldsymbol{\beta}}$ with the penalized LSE in Theorem 7. Let $\widetilde{\boldsymbol{\mu}}$ and $\widehat{\boldsymbol{\mu}}^o$ be the orthogonal projections of $\mathbf{y}$ to the linear spans of $\{\mathbf{x}_j, j \leq p\}$ and $\{\mathbf{x}_j, \beta_j \neq 0\}$, respectively. For uncorrelated errors with common variance $\sigma^2$, the degrees of freedom for $\widehat{\boldsymbol{\mu}}^o$ is $\sum_{j=1}^p \operatorname{Cov}(\widetilde{\mu}_j, \widehat{\mu}_j^o)/\sigma^2 = \operatorname{rank}(\mathbf{x}_j : \beta_j \neq 0)$. Thus, since $E\|\widetilde{\boldsymbol{\mu}} - \boldsymbol{\mu}\|^2 = \sigma^2 \operatorname{rank}(\mathbf{X})$ and $\|\widehat{\boldsymbol{\mu}} - \boldsymbol{\mu}\|^2 + \|\widetilde{\boldsymbol{\mu}} - \boldsymbol{\mu}\|^2 - \|\widetilde{\boldsymbol{\mu}} - \widehat{\boldsymbol{\mu}}\|^2 = 2(\widetilde{\boldsymbol{\mu}} - \boldsymbol{\mu})'(\widehat{\boldsymbol{\mu}} - \boldsymbol{\mu})$,

$$
(5.8) \quad \operatorname{df}(\widehat{\boldsymbol{\mu}}) \equiv \sum_{j=1}^p \frac{\operatorname{Cov}(\widetilde{\mu}_j, \widehat{\mu}_j)}{\sigma^2} = \frac{1}{2}E\left( \operatorname{rank}(\mathbf{X}) - \frac{\|\widetilde{\boldsymbol{\mu}} - \widehat{\boldsymbol{\mu}}\|^2}{\sigma^2} + \frac{\|\widehat{\boldsymbol{\mu}} - \boldsymbol{\mu}\|^2}{\sigma^2} \right)
$$

extends the notion of degrees of freedom. This also provides the $C_p$-type risk estimate

$$
(5.9) \quad \widehat{C}_p \equiv \widehat{C}_p(\lambda) \equiv \|\widetilde{\boldsymbol{\mu}} - \widehat{\boldsymbol{\mu}}\|^2 + \widehat{\sigma}^2\{2\widehat{\operatorname{df}} - \operatorname{rank}(\mathbf{X})\} \approx \|\widehat{\boldsymbol{\mu}} - \boldsymbol{\mu}\|^2.
$$

Theorem 7 suggests the unbiased estimator for the degrees of freedom (5.8) as

$$\widehat{\mathrm{df}} \equiv \widehat{\mathrm{df}}(\lambda) \equiv \mathrm{trace}(\mathbf{Q}^{-1}(\widehat{\boldsymbol{\beta}}; \lambda)\widehat{\mathbf{P}}\boldsymbol{\Sigma}\widehat{\mathbf{P}}') \tag{5.10}$$

and the related $C_p$-type estimator of the MSE $E\|\widehat{\boldsymbol{\mu}} - \boldsymbol{\mu}\|^2$ via (5.9). We refer to Efron (1986) and Meyer and Woodroofe (2000) for more discussions about (5.8) and (5.9). We will present in Section 6 simulation results to demonstrate that (5.9) provides a reasonable risk estimator. The following theorem asserts the unbiasedness of (5.8) and (5.9).

THEOREM 8. *Suppose (2.5) holds with $d^* = p$. Then, the SURE method provides unbiased estimators for the degrees of freedom and the $\ell_2$ risk for the estimation of the mean vector,*

$$E(\widehat{\mathrm{df}}) = \mathrm{df}(\widehat{\boldsymbol{\mu}}), \qquad E\widehat{C}_p = E\|\widehat{\boldsymbol{\mu}} - \boldsymbol{\mu}\|^2, \tag{5.11}$$

*in the linear model (1.2), where $\mathrm{df}(\widehat{\boldsymbol{\mu}})$, $\widehat{\mathrm{df}}$ and $\widehat{C}_p$ are, respectively, given by (5.8), (5.10) and (5.9), the $\widehat{\sigma}^2$ in (5.9) is as in (5.6), and $\widehat{\boldsymbol{\mu}} = \mathbf{X}\widehat{\boldsymbol{\beta}}$ is as in Theorem 7 with a fixed $\lambda$. Furthermore, if $\rho(t; \lambda) = \lambda t$ for the LASSO or $|\widehat{\beta}_j| > \gamma\lambda$ for all $\widehat{\beta}_j \neq 0$ under (2.3), then*

$$\widehat{\mathrm{df}} = \#\{j : \widehat{\beta}_j \neq 0\}. \tag{5.12}$$

Under a positive cone condition on $\mathbf{X}$, Efron et al. (2004) proved the unbiasedness of $\#\{j : \widehat{\beta}_j \neq 0\}$ as an estimator for the degrees of freedom for the LARS estimator (not the LASSO) at a fixed step $k$. Our definition of the degrees of freedom and $C_p$ is slightly different, since we use $\|\widetilde{\boldsymbol{\mu}} - \widehat{\boldsymbol{\mu}}\|^2$ and $\mathrm{rank}(\mathbf{X})$ in (5.8) and (5.9) for variance reduction, instead of $\|\mathbf{y} - \widehat{\boldsymbol{\mu}}\|^2$ and $n$. We prove $E\#\{j : \widehat{\beta}_j \neq 0\} = \mathrm{df}(\widehat{\boldsymbol{\mu}})$ for the LASSO for fixed $\lambda$ without requiring the positive cone condition, but not for fixed $k$ with a stochastic $\lambda$. The performances of $\widehat{C}_p$ for the LASSO and MC+ are similar in our simulation experiments.

5.2. *Estimation of noise level.* Consider throughout this subsection standardized designs with $\|\mathbf{x}_j\|^2 = n$ for all $j \leq p$ in (1.2). We have shown in Theorem 1 and Table 1 that the MC+ at $\lambda_{\mathrm{univ}} \equiv \sigma\sqrt{(2/n)\log p}$ works well for variable selection. In practice, this requires a reasonable estimate of the noise level $\sigma$. For $p < n$, the mean residual squares $\|\mathbf{y} - \widetilde{\boldsymbol{\mu}}\|^2/\{n - \mathrm{rank}(\mathbf{X})\}$ for the full model provides an unbiased estimator of $\sigma^2$ as in Table 1, where $\widetilde{\boldsymbol{\mu}}$ is the orthogonal projection of $\mathbf{y}$ to the linear span of $\{\mathbf{x}_j, j \leq p\}$. However, the estimation of $\sigma^2$ is a more delicate problem for $p > n$ or small $n - p > 0$. Here, we present a simple estimator of $\sigma^2$ in such cases based on Theorem 8.

Since (2.8) provides estimates $\widehat{\boldsymbol{\mu}}(\lambda) \equiv \mathbf{X}\widehat{\boldsymbol{\beta}}(\lambda)$ of the mean $\boldsymbol{\mu} \equiv \mathbf{X}\boldsymbol{\beta}$, we may use

$$\widehat{\sigma}^2(\lambda) \equiv \|\mathbf{y} - \widehat{\boldsymbol{\mu}}(\lambda)\|^2 / \{n - \widehat{\mathrm{df}}(\lambda)\} \tag{5.13}$$

to estimate $\sigma^2$, with the $\widehat{\mathrm{df}}(\lambda)$ in (5.10) as an adjustment for the degrees of freedom. Still, good $\widehat{\sigma}^2(\lambda)$ requires a consistent $\widehat{\boldsymbol{\mu}}(\lambda)$, which depends on the choice of a suitable $\lambda$ of the order $\sigma\sqrt{(\log p)/n}$. This circular estimation problem can be solved with

$$\widehat{\sigma} \equiv \widehat{\sigma}(\widehat{\lambda}), \qquad \widehat{\lambda} \equiv \min\{\lambda \geq \lambda_* : \widehat{\sigma}^2(\lambda) \leq n\lambda^2/(r_0 \log p)\}, \tag{5.14}$$

for suitable $r_0 \leq 2$ and $\lambda_* > 0$. Here, $\lambda_*$ could be preassigned or determined by upper bounds on $\widehat{\mathrm{df}}(\lambda)$ or the dimension $\#\{j : \widehat{\beta}_j(\lambda) \neq 0\}$. In principle, we may also use in (5.14) estimates $\widehat{\sigma}^2(\lambda)$ based on cross-validation or bootstrap, but the computationally much simpler (5.13) turns out to have the best overall performance in our simulation experiments.

5.3. *Convexity, continuity and almost differentiability.* Here, we consider the continuity and almost differentiability of a penalized LSE $\widehat{\boldsymbol{\beta}}$, which the proof of Theorems 7 and 8 require.

The continuity of $\widehat{\boldsymbol{\beta}}$, demanded by Stein (1981), is a property of independent interest on its own right for robust estimation [Fan and Li (2001)]. For full rank designs, we provide here the equivalence of the continuity of the penalized LSE and the global convexity condition (4.2). We have considered (2.3) for unbiased selection. For the continuity of $\widehat{\boldsymbol{\beta}}$, we only need

$$\lim_{t \to \infty} \rho(t; \lambda)/t^2 = 0, \qquad 0 \leq \dot{\rho}(0+; \lambda) < \infty. \tag{5.15}$$

THEOREM 9. *Let $\lambda$ be fixed. Suppose $\rho(t; \lambda)$ is continuously differentiable in $t > 0$, (5.15) holds, and $\mathrm{rank}(\mathbf{X}) = p$. Then the following three statements are equivalent to each other:*

  (i) *The global minimizer $\widehat{\boldsymbol{\beta}}$ of (1.1) is unique and continuous in $\mathbf{y} \in \mathbb{R}^n$.*
  (ii) *The global convexity condition (4.2) holds.*
  (iii) *The penalized loss $L(\mathbf{b}; \lambda)$ in (1.1) is strictly convex in $\mathbf{b} \in \mathbb{R}^p$.*

For $p > n$, an implication of Theorem 9 is the continuity of solution $\widehat{\boldsymbol{\beta}}$ of the estimating equation (2.6) subject to $\{j : \widehat{\beta}_j \neq 0\} \subset A$ for all fixed $\lambda$ and $A$ with $|A| \leq d^*$, provided the sparse convexity (2.5). Thus, minimizing the maximum concavity allows the broadest extent for such sparse continuity of solutions of (2.6). The most difficult part of the proof of Theorem 9 is (i) $\Rightarrow$ (ii), which is done by showing $(x, x, \ldots, x)' = x\mathbf{1}$ is in the range of $\widehat{\boldsymbol{\beta}}$ for all $x > 0$. Since the penalized loss attains minimum at $\widehat{\boldsymbol{\beta}}$, $\mathbf{Q}(\widehat{\boldsymbol{\beta}}; \lambda)$ in (5.4) is positive definite for smooth penalties, and the positive-definiteness of $\mathbf{Q}(t\mathbf{1}; \lambda)$ gives $c_{\min}(\boldsymbol{\Sigma}) > \ddot{\rho}(t; \lambda)$.

The application of SURE in Theorems 7 and 8 also requires the almost differentiability of $\widehat{\boldsymbol{\beta}}$. In the following proposition, we establish the stronger Liptchitz condition for $\widehat{\boldsymbol{\beta}}$ under the conditions of Theorem 7.

PROPOSITION 2. *Let $\lambda$ and $\mathbf{X}$ be fixed and treat $\widehat{\boldsymbol{\beta}}$ in (2.8) as a function of $\mathbf{y}$. Suppose (2.5) holds with $d^* = p$. Then $\widehat{\boldsymbol{\beta}} = \mathbf{h}(\widetilde{\mathbf{z}})$ for $\widetilde{\mathbf{z}} = \mathbf{X}'\mathbf{y}/n \in \mathbb{R}^p$ and a certain almost differentiable function $\mathbf{h} : \mathbb{R}^p \to \mathbb{R}^p$, such that for all $\mathbf{z}$ and $\mathbf{v}$ in $\mathbb{R}^p$*

$$(5.16) \qquad \mathbf{h}(\mathbf{z} + \mathbf{v}) = \mathbf{h}(\mathbf{z}) + \left\{ \int_0^1 (\mathbf{P}'\mathbf{Q}^{-1}\mathbf{P})(\mathbf{h}(\mathbf{z} + x\mathbf{v}); \lambda)\, dx \right\} \mathbf{v},$$

*where $\mathbf{Q}$ is as in (5.4) and $\mathbf{P}(\boldsymbol{\beta}; \lambda) : \mathbf{b} \to (b_j : \beta_j \neq 0)'$ is as in (3.12). Consequently, $\mathbf{h}(\mathbf{z})$ satisfies the Lipschitz condition $\|\mathbf{h}(\mathbf{z} + \mathbf{v}) - \mathbf{h}(\mathbf{z})\| \leq \|\mathbf{v}\|/ \{c_{\min}(\boldsymbol{\Sigma}) - \kappa(\rho; \lambda)\}$.*

**6. More simulation results.** In this section, we present simulation results along with some discussion on the performance of the LASSO, MC+ and SCAD+ in selection consistency and estimation of $\boldsymbol{\beta}$ and $\boldsymbol{\mu} \equiv \mathbf{X}\boldsymbol{\beta}$, sparse recovery, the computational complexity and the scalability of the PLUS algorithm, the choice of the tuning parameter $\gamma$, the estimation of the noise level $\sigma$ and the risk, and the sparse Riesz condition.

6.1. *Selection consistency.* For the MC+, the tuning parameter $\gamma$ regulates its computational complexity and bias level. We study its effects through three experiments, say experiments 1, 2 and 3, including cases where $\gamma$ is smaller than the "smallest" theoretical value $1/(1 - \max_{j \neq k} |\mathbf{x}_j'\mathbf{x}_k|/n)$ with $d^* = 2$ in (2.5) and $\lambda < \lambda_{\text{univ}}$.

Experiment 1, summarized in Table 1 in Section 2, illustrates the superior selection accuracy of the MC+ for sparse $\boldsymbol{\beta}$, compared with the LASSO and SCAD+. Experiment 2, summarized in Table 3, shows the effects of the regularization parameter $\gamma$ on selection accuracy and computational complexity of the MC+. Experiment 3, summarized in Table 4, demonstrates the scalability of the PLUS algorithm for large $p$. The design matrix $\mathbf{X}$ has the same distribution in experiments 1 and 2. For each replication, we generate a $300 \times 600$ random matrix as the difference of two independent random matrices, the first with i.i.d. unit exponential entries and the second i.i.d. $\chi_1^2$ entries. We normalize the 600 columns of this difference matrix to summation zero and Euclidean length $\sqrt{n}$. We then sequentially sample groups of 10 vectors from this pool of normalized columns. For the $m$th

TABLE 3

*Performance of MC+ with different $\gamma$ in experiment 2 100 replications, $n = 300$, $p = 200$, $d^o = 30$, $\beta_* = 3/8$, LASSO for $\gamma = \infty$ $CS \equiv I\{\widehat{A} = A^o\}$, $SE_{\boldsymbol{\beta}} \equiv \|\widehat{\boldsymbol{\beta}} - \boldsymbol{\beta}\|^2$, $SE_{\boldsymbol{\mu}} \equiv \|\mathbf{X}(\widehat{\boldsymbol{\beta}} - \boldsymbol{\beta})\|$, $K \equiv \#(steps)$*

|  | $\gamma$ | **1.01** | **1.4** | **1.7** | **2.0** | **2.4** | **2.7** | **3.0** | **5.0** | **∞** |
|---|---|---|---|---|---|---|---|---|---|---|
| $\lambda_{\text{univ}}$ $= 0.188$ | $\overline{\text{CS}}$ | 0.81 | **0.82** | 0.66 | 0.53 | 0.34 | 0.35 | 0.27 | 0.11 | 0.00 |
|  | $\overline{\text{SE}}_{\boldsymbol{\beta}}$ | 0.136 | **0.128** | 0.265 | 0.495 | 0.729 | 0.801 | 0.817 | 1.007 | 1.420 |
|  | $\overline{\text{SE}}_{\boldsymbol{\mu}}$ | 0.117 | **0.112** | 0.205 | 0.358 | 0.510 | 0.564 | 0.583 | 0.761 | 1.123 |
|  | $\overline{k}$ | 561 | 98 | 62 | 47 | 36 | 33 | 32 | 32 | 34 |
| $\lambda$ for max $\overline{\text{CS}}$ | $\lambda$ | 0.195 | 0.182 | 0.175 | 0.164 | 0.164 | 0.158 | 0.169 | 0.188 | 0.201 |
|  | $\overline{\text{CS}}$ | **0.83** | 0.82 | 0.75 | 0.63 | 0.46 | 0.36 | 0.27 | 0.11 | 0.02 |
|  | $\overline{k}$ | 561 | 98 | 65 | 57 | 45 | 40 | 34 | 32 | 33 |
| $\lambda$ for min $\overline{\text{SE}}_{\boldsymbol{\beta}}$ | $\lambda$ | 0.182 | 0.175 | 0.153 | 0.138 | 0.120 | 0.108 | 0.101 | 0.094 | 0.050 |
|  | $\overline{\text{SE}}_{\boldsymbol{\beta}}$ | 0.132 | **0.117** | 0.119 | 0.124 | 0.133 | 0.140 | 0.149 | 0.255 | 0.394 |
|  | $\overline{k}$ | 562 | 98 | 68 | 64 | 65 | 67 | 68 | 47 | 84 |
| $\lambda$ for min $\overline{\text{SE}}_{\boldsymbol{\mu}}$ | $\lambda$ | 0.182 | 0.175 | 0.153 | 0.138 | 0.120 | 0.108 | 0.101 | 0.094 | 0.050 |
|  | $\overline{\text{SE}}_{\boldsymbol{\mu}}$ | 0.115 | **0.104** | 0.106 | 0.110 | 0.117 | 0.124 | 0.130 | 0.201 | 0.278 |
|  | $\overline{k}$ | 562 | 98 | 68 | 64 | 65 | 67 | 68 | 47 | 84 |

group, we sample from the remaining $610 - 10m$ columns one member as $\mathbf{x}_{10m-9}$ and 9 more to maximize the absolute correlation $|\mathbf{x}_j' \mathbf{x}_{10m-9}|/n$, $j = 10m - 8, \ldots, 10m$. In experiment 3, $\mathbf{X}$ are generated in the same way for each replication with groups of size 50 from a pool of 6000 i.i.d. columns. In all the three experiments, $\beta_j = \pm \beta_*$ for $j \in A^o$ and $\boldsymbol{\varepsilon} \sim N(0, \mathbf{I}_n)$.

Strong effects of bias on selection accuracy is observed in all three tables. In Table 1 where $\beta_* \approx \sqrt{10}\lambda_{\text{univ}}$, the selection accuracy of the LASSO clearly deteriorates as $d^o$ increases. In Tables 3 and 4, the unbiasedness criterion $\beta_* > \gamma \lambda_{\text{univ}}$ in (2.9) matches the best selection results well, with $1.7\lambda_{\text{univ}} < \beta_* < 2\lambda_{\text{univ}}$ in Table 3 and $2\lambda_{\text{univ}} < \beta_* < 2.4\lambda_{\text{univ}}$ in Table 4. In Table 1, $\gamma \lambda_{\text{univ}}/\sigma \approx 1/2 = \beta_*/\sigma$, but slightly larger $\lambda$ yields the largest $\overline{\text{CS}}$. Comparison between the results for $\lambda_{\text{univ}}$ and $\arg\max_\lambda \overline{\text{CS}}$ in all three tables demonstrates that $\lambda_{\text{univ}}$ is a reasonable choice for variable selection with $\|\mathbf{x}_j\|^2 = n$, especially when $\beta_*$ is near the minimum for accurate selection as in Tables 3 and 4.

An interesting phenomenon exhibited in experiments 2 and 3 is that the observed selection accuracy $\overline{\text{CS}}$ is always decreasing in $\gamma$. Despite the computational complexity for small $\gamma$, the MC+ still recovers the true $A^o$ among so many parallelepipeds it traverses through. This suggests that the interference of the bias, not the complexity of the path or the lack of the convexity of the penalized loss, is a dominant factor in variable selection. Of course, bias reduction does not always provide accurate variable selection. When the signal is reduced to $\beta_* = 1/4$ from $\beta_* = 3/8$ in experiment 2, the selection accuracy suddenly drops to $\overline{\mathrm{CS}} \leq 0.11$ for all values of $(\lambda, \gamma)$.

6.2. *Estimation of regression coefficients and the mean responses.* Tables 1, 3 and 4 also report results for the estimation of regression coefficients $\boldsymbol{\beta}$ with the square error $\mathrm{SE}_{\boldsymbol{\beta}} \equiv \|\widehat{\boldsymbol{\beta}} - \boldsymbol{\beta}\|^2$. The MC+ and SCAD+ clearly outperform the LASSO in these settings. In Table 4, the minimum $\overline{\mathrm{SE}}_{\boldsymbol{\beta}}$ for the SCAD+ are 2.5% and 6.2% smaller than the MC+ with matching $\gamma = 2.4$ and 2.7, respectively, while those of the MC+ are 14% and 16% smaller than the SCAD+ with matching maximum concavity $\kappa(\rho)$ ($\gamma = 1.4$ and 1.7 for the MC+ versus $\gamma = 2.4$ and 2.7 for the SCAD+, respectively). The SCAD penalty requires $\gamma > 2$. The results for the SCAD+ in experiment 2 are not reported since they show a similar pattern as experiment 3. Results for the estimation of the mean $\boldsymbol{\mu} \equiv \mathbf{X}\boldsymbol{\beta}$ with the average squared error $\mathrm{SE}_{\boldsymbol{\mu}} \equiv \|\mathbf{X}\widehat{\boldsymbol{\beta}} - \mathbf{X}\boldsymbol{\beta}\|^2/n$ are similar to those for the estimation of $\boldsymbol{\beta}$ in Tables 3 and 4.

6.3. *Computational complexity and choice of $\gamma$.* As expected, we observe in Tables 3 and 4 that the MC+ with smaller $\gamma$ is computationally more costly. Dramatic rise in the number of needed PLUS steps is observed when $\gamma$ decreases to 1/2 in experiment 2. We avoid $\gamma = 1$, since it produces the singularity $\det(\mathbf{Q}(\boldsymbol{\eta})) = 0$ for (3.12) whenever $\sum_{j=1}^{p} |\eta_j| = 1$ for the MC+ with the standardization $\|\mathbf{x}_j\|^2 = n$.

Table 4

*Performance of MC+ and SCAD with $p > n$ in experiment 3 100 replications, $n = 300$, $p = 2000$, $d^o = 30$, $\beta_* = 1/2$, SCAD+ for $\gamma*$, LASSO with $\gamma = \infty$*

|  | $\gamma$ | **1.4** | **1.7** | **2.0** | **2.4** | **2.7** | **2.4\*** | **2.7\*** | **$\infty$** |
|---|---|---|---|---|---|---|---|---|---|
| $\lambda_{\mathrm{univ}}$ | $\overline{\mathrm{CS}}$ | **0.99** | **0.99** | 0.96 | 0.80 | 0.56 | 0.00 | 0.00 | 0.00 |
| $= 0.225$ | $\overline{\mathrm{SE}}_{\boldsymbol{\beta}}$ | **0.109** | 0.116 | 0.205 | 0.534 | 0.712 | 2.703 | 2.764 | 2.640 |
|  | $\overline{\mathrm{SE}}_{\boldsymbol{\mu}}$ | **0.098** | 0.103 | 0.170 | 0.395 | 0.515 | 1.602 | 1.661 | 1.785 |
|  | $\overline{k}$ | 119 | 76 | 62 | 46 | 41 | 130 | 84 | 56 |
| $\lambda$ | $\lambda$ | 0.241 | 0.225 | 0.225 | 0.225 | 0.210 | 0.177 | 0.171 |  |
| for max $\overline{\mathrm{CS}}$ | $\overline{\mathrm{CS}}$ | **1.00** | 0.99 | 0.96 | 0.80 | 0.60 | 0.08 | 0.02 | 0.00 |
|  | $\overline{k}$ | 118 | 76 | 62 | 46 | 44 | 255 | 169 |  |
| $\lambda$ | $\lambda$ | 0.225 | 0.203 | 0.183 | 0.165 | 0.149 | 0.134 | 0.129 | 0.069 |
| for min $\overline{\mathrm{SE}}_{\boldsymbol{\beta}}$ | $\overline{\mathrm{SE}}_{\boldsymbol{\beta}}$ | **0.109** | 0.112 | 0.117 | 0.127 | 0.138 | 0.124 | 0.130 | 1.292 |
|  | $\overline{k}$ | 119 | 77 | 69 | 71 | 76 | 279 | 200 | 181 |
| $\lambda$ | $\lambda$ | 0.225 | 0.203 | 0.183 | 0.165 | 0.149 | 0.143 | 0.134 | 0.069 |
| for min $\overline{\mathrm{SE}}_{\boldsymbol{\mu}}$ | $\overline{\mathrm{SE}}_{\boldsymbol{\mu}}$ | **0.098** | 0.100 | 0.104 | 0.112 | 0.122 | 0.112 | 0.118 | 0.563 |
|  | $\overline{k}$ | 119 | 77 | 69 | 71 | 76 | 273 | 197 | 181 |

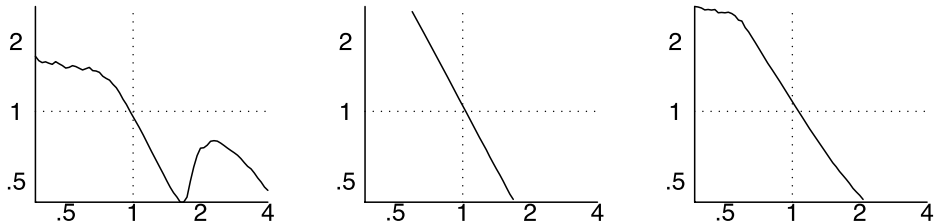

Fig. 5. *The median of $\widehat{\sigma}^2(\lambda)/(n\lambda^2/\log p)$ as a function of $\lambda/\sqrt{(\log p)/n} \in [2^{-3/2}, 4]$ based on 100 replications. Left: experiment 4 with $n = 300$, $p = 2000$ and $d^o = 30$. Middle and right: experiment 5 with high and low correlations, respectively, $n = 600$, $p = 3000$ and $d^o = 35$. For $1/1.5 \le r_0 \le 2$ in (5.14), $\widehat{\sigma}^2(\lambda)/(n\lambda^2/\log p) \approx 1/r_0$ matches $\lambda/\sqrt{(\log p)/n} = \sqrt{r_0}$ reasonably well to provide $\widehat{\sigma}^2(\lambda) \approx \sigma^2 = 1$. This is especially the case for $r_0 = 1$ as indicated by the dotted lines.*

Table 4 shows that the PLUS algorithm scales well for $p > n$. Comparisons between Tables 3 and 4 demonstrate that for similar $d^o$ and SNR $\beta_*/\lambda_{\mathrm{univ}}$, the computational complexity of the MC+ is insensitive to $p$ as measured by the average number of steps $\overline{k}$.

In practice, full implementation of the MC+ requires a specification of $\gamma$ and possibly a stopping rule for large $(n, p)$, say $k = k_{\max} \wedge k^*$, to allow the algorithm to end before it reaches the perfect fit at $k = k^*$. As we have discussed in the Introduction, large $\gamma$ provides computational simplicity but may harm selection consistency with larger bias. Our simulation results in Tables 3 and 4 demonstrate robust selection accuracy for smaller-than-necessary $\gamma > 0$ at the universal penalty level. Thus, the choice of $\gamma$ should largely be determined by the available computational resources as long as the MC+ path reaches a sufficiently small $\lambda$. In our simulations, $k_{\max} = 5000$, and all replications failing to reach $\lambda < \lambda_*/1.2$ occur only for unreasonably small $\gamma = 1/2$, where $\lambda_*$ is (much) smaller than the smallest reported penalty level in each experiment. Since $\hat{\sigma}$ in (5.13) is based on the beginning segments of the PLUS path, we "know" whether the desired penalty level is reached.

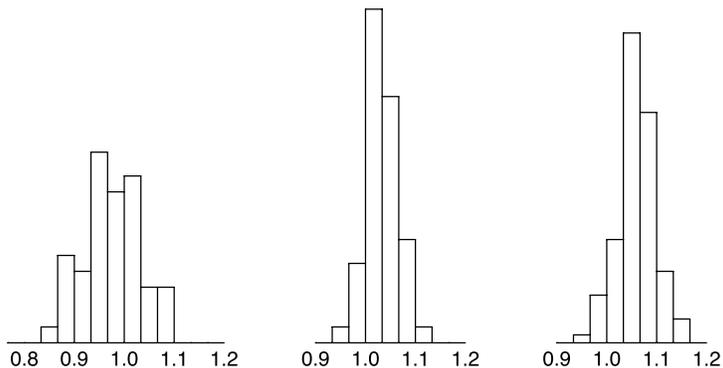

Fig. 6. *Histograms of $\widehat{\sigma}$ at $r_0 = 1$ for the same simulations as in Figure 5 with respective means and standard deviations $0.971 \pm 0.057$, $1.033 \pm 0.032$ and $1.060 \pm 0.039$ from the left to the right. It turns out that the MSE for $\widehat{\sigma}$ is of the same order as $n^{-1/2}$ in these simulations.*

6.4. *Estimation of noise.* In Figures 5 and 6, we present simulation results for the estimation of $\sigma$ in experiments 4 and 5 with the MC+ estimator $\hat{\boldsymbol{\mu}}(\lambda) = \mathbf{X}\hat{\boldsymbol{\beta}}(\lambda)$. In experiment 4, $(n, p) = (300, 2000)$, $\gamma = 1.7$, $\beta_* = 1/2$, $\boldsymbol{\beta}$ is generated every 10 replications and $\mathbf{X}$ is fixed. Its configurations are otherwise identical to that of experiment 3 reported in Table 4. In experiment 5, $(n, p) = (600, 3000)$, $\mathbf{x}_j$ are normalized columns from a Gaussian random matrix with i.i.d. rows and the correlation $\sigma_{j,k} = \sigma_{1,2}^{|k-j|}$ among entries within each row, $\gamma = 2/(1 - \max_{j>k} |\mathbf{x}_k' \mathbf{x}_j|/n)$ as in experiment 1, the nonzero $\beta_j$ are composed of 5 blocks of $\beta_*(1, 2, 3, 4, 3, 2, 1)'$ centered at random multiples $j_1, \ldots, j_5$ of 25, $\beta_*$ sets $\|\mathbf{X}\boldsymbol{\beta}\|^2/n = 3$, $\boldsymbol{\varepsilon} \sim N(0, \mathbf{I}_n)$, and $\{\mathbf{X}, \boldsymbol{\beta}\}$ are generated every 10 replications. It has two settings: $\sigma_{1,2} = 0.9$ for high correlation and $\sigma_{1,2} = 0.1$ for low correlation. We set $\lambda_* = \{2^{-3}(\log p)/n\}^{1/2}$ in both experiments 4 and 5.

Figure 5 plots the median of $\hat{\sigma}^2(\lambda)/(n\lambda^2/\log p)$ versus $\lambda/\sqrt{(\log p)/n}$ in the simulations described above. Since all three curves cross the level $\hat{\sigma}^2(\lambda)/(n\lambda^2/\log p) = 1$ at approximately $\lambda/\sqrt{(\log p)/n} = 1$, the estimation equation (5.14) provides approximately the right answer $\hat{\sigma}^2 \approx 1$ for $r_0 = 1$. We solve (5.14) for individual replications and plot the histograms of $\hat{\sigma}$ in Figure 6. These simulation results suggest that the MSE for $\hat{\sigma}$ is of the same order as $n^{-1/2}$ for sparse $\boldsymbol{\beta}$.

6.5. *Sparse recovery.* Our variable selection theorems are applicable to sparse recovery in the noiseless case of $\sigma = 0$ as we mentioned in Remark 6. Table 5 reports simulation results to show that the LASSO ($y = \infty$) may miss up to about 45% of nonzero $\beta_j$, while the MC+ ($\gamma = 3$) still manages to recover the true $\boldsymbol{\beta}$. For $(n, p, d^o) = (100, 2000, 28)$ and $(200, 10,000, 40)$, the LASSO does not capture most of the nonzero $\beta_j$ before falsely selected variables manage to perfectly fit $\mathbf{y} = \mathbf{X}\boldsymbol{\beta}$ at the last step of the LARS, at the expense of substantially many additional computation steps.

6.6. *Estimation of risk.* We summarize in Figure 7 the performance of $\hat{C}_p$ in (5.9) for the MC+ in experiments 4 and 5, with the df in (5.10) and the $\hat{\sigma}$ in (5.14). For each of the three settings, $E\|\hat{\boldsymbol{\mu}}(\lambda) - \boldsymbol{\mu}\|^2$ and $E\hat{C}_p(\lambda)$ are approximated by the averages in 100 replications and the expected conditional variance $E \operatorname{Var}(\hat{C}_p(\lambda)|\mathbf{X}, \boldsymbol{\beta})$ is approximated by the within-group variance, since $(\mathbf{X}, \boldsymbol{\beta})$ is unchanged in every 10 replications in each of the three settings. From Figure 7, we observe that the MSE $E\|\widehat{\boldsymbol{\mu}}(\lambda) - \boldsymbol{\mu}\|^2$ is reasonably approximated by $\widehat{C}_p(\lambda)$ for $p > n$, at least before the MC+ starts to over fit with small $\lambda$.

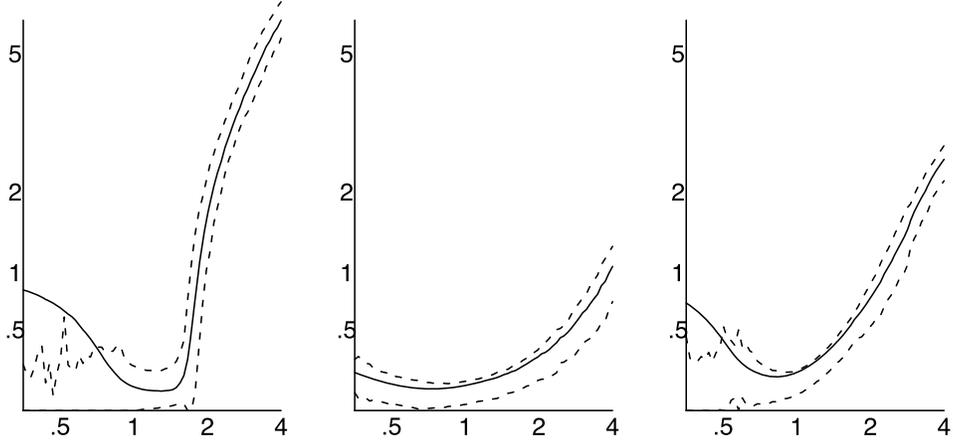

Fig. 7. *Approximations of $E\|\widehat{\boldsymbol{\mu}}(\lambda) - \boldsymbol{\mu}\|^2/n$ (solid) and $E\widehat{C}_p(\lambda)/n \pm 2\{E\,\mathrm{Var}(\widehat{C}_p(\lambda)/n|\mathbf{X}, \boldsymbol{\beta})\}^{1/2}$ (dashed) as functions of $\lambda/\sqrt{(\log p)/n}$ for the MC+ based on the same simulations as in Figure 5. The MSE $E\|\widehat{\boldsymbol{\mu}}(\lambda) - \boldsymbol{\mu}\|^2$ is reasonably approximated by $\widehat{C}_p(\lambda)$ in these experiments with $p > n$, at least before the MC+ starts to over fit with small $\lambda$.*

6.7. *The sparse Riesz condition.* The SRC (2.11) and constant factors in Theorems 1, 2, 4 and 5 are quite conservative compared with our simulation results. Technically, this is probably due to the following two reasons: (i) the sparse minimum and maximum eigenvalues, or $c_*$ and $c^*$, respectively, in (2.11), are used to bound the effects of matrix operations in the worst case scenario given the dimension/rank of the matrix; (ii) we use the con-servative bound $c_{\min}(\boldsymbol{\Sigma}_{j:\eta_j \neq 0} - \mathrm{diag}(1/\gamma, |\eta_j| = 1)) \geq c_{\min}(\boldsymbol{\Sigma}_{j:\eta_j \neq 0}) - 1/\gamma$ to ensure sparse convexity in the $k$th segment $\boldsymbol{\eta}^{(k)}$ of the MC+ path, but $\#\{j : |\eta_j^{(k)}| = 1\}$ could be much smaller than $\#\{j : \eta_j^{(k)} \neq 0\}$. These considerations suggest that the penalized loss (1.1) with the MCP (2.1) possesses sufficient convexity if

$$(6.1) \qquad P^* \{ c_{\min}(\boldsymbol{\Sigma}_A) \geq \kappa(\rho_2) = 1/\gamma \,||\, |A| = d, \mathbf{X} \} \approx 1$$

at a reasonable dimension $d$, where $P^*$ is the probability under which $A$ is a random subset of $\{1, \ldots, p\}$. In practice, we may substitute the SRC (2.11) with (6.1) and a similar probabilistic upper bound on $c_{\max}(\boldsymbol{\Sigma}_A)$ under $P^*$, which are weaker and much easier to check. Figure 8 plots the mean and a lower confidence bound of $c_{\min}(\boldsymbol{\Sigma}_A)$ under $P^*$ as functions of given $d = |A|$. We observe that (6.1) holds for quite a few possible combinations of $(d, \gamma)$ in our experiments, in view of Tables 1, 3 and 4.

Table 5
*Sparse recovery with MC+ at the last PLUS step $k^*$. Entries of $\mathbf{X}$ and nonzero $\beta_j$ are i.i.d. $N(0,1)$, $\boldsymbol{\varepsilon} = 0$; FN $\equiv \#\{j : \widehat{\beta}_j^{(k^*)} = 0 \neq \beta_j\}$*

| $(n, p, d^o)$ | $(100, 2000, 15)$ | | $(100, 2000, 28)$ | | $(200, 10{,}000, 40)$ | |
|---|---|---|---|---|---|---|
| $\gamma$ | **3** | $\boldsymbol{\infty}$ | **3** | $\boldsymbol{\infty}$ | **3** | $\boldsymbol{\infty}$ |
| $\%\{\widehat{\boldsymbol{\beta}}^{(k^*)} = \boldsymbol{\beta}\}$ | 100 | 51 | 73 | 0 | 100 | 0 |
| mean(FN$|\widehat{\boldsymbol{\beta}}^{(k^*)} \neq \boldsymbol{\beta}$) | | 2 | 19 | 13 | | 18 |
| mean($k^*|\widehat{\boldsymbol{\beta}}^{(k^*)} = \boldsymbol{\beta}$) | 32 | 65 | 87 | | 102 | |
| mean($k^*|\widehat{\boldsymbol{\beta}}^{(k^*)} \neq \boldsymbol{\beta}$) | | 144 | 513 | 153 | | 311 |

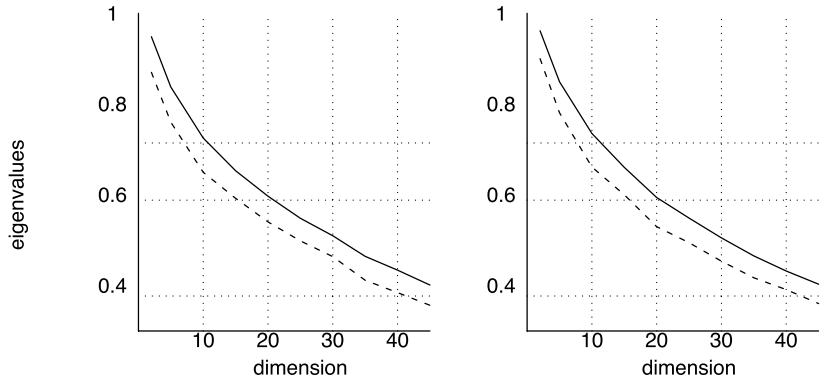

FIG. 8. *The mean (solid) of the minimum eigenvalue $c_{\min}(\mathbf{X}'_A \mathbf{X}_A / n)$ for a random set $A$ of design vectors and the mean minus two standard deviations (dashed) as functions of the dimension $|A|$, each point based on 100 replications, with horizontal dotted lines at $\kappa(\rho_2) = 1/\gamma$ for $\gamma \in \{1.4, 1.7, 2.652\}$. Left: the design $\mathbf{X}$ in experiments 1 and 2. Right: the design $\mathbf{X}$ in experiments 3 and 4.*

## 7. Discussion.
We have introduced and studied the MC+ methodology for unbiased penalized selection. Our theoretical and simulation results have shown the superior selection accuracy of this method and the computational efficiency of the PLUS algorithm. We have provided an oracle inequality to demonstrate the advantage of the MC+ for the estimation of regression coefficients and proved its convergence at certain minimax rates in $\ell_r$ balls. We have also discussed unbiased estimation of the risk, estimation of the noise level in the linear model in the case of $p > n$, and the necessary and sufficient conditions for the continuity of the penalized LSE. In this section, we briefly discuss the choice among multiple solutions in the PLUS path, the one-at-a-time condition with the PLUS algorithm, the penalized LSE for orthogonal designs, adaptive penalty, general loss and sub-Gaussian errors.

7.1. *Choice among multiple solutions in the path.* In (2.8), $\widehat{\boldsymbol{\beta}}(\lambda)$ is taken as the $\widehat{\boldsymbol{\beta}}^{(x)}$ when $\lambda^{(x)}$ first reaches a level no greater than $\lambda$. An alternative choice [Zhang (2007b)] is to pick $\widehat{\boldsymbol{\beta}}(\lambda)$ as the sparsest $\widehat{\boldsymbol{\beta}}^{(x)}$ in (2.7) with $\lambda^{(x)} = \lambda$. Theorems 1, 4 and 5 holds verbatim for the sparsest solution, while Theorem 2 holds with a smaller $d_* = d^*/(c^*/c_* + 3/2)$. Our simulation experiments yield nearly identical results among the two choices. A significant reason for using (2.8) is its simplicity in implementation since it does not require the entire path to compute $\widehat{\boldsymbol{\beta}}(\lambda)$ for given penalty levels $\lambda$.

7.2. *The one-at-a-time condition with the PLUS algorithm.* The formulas (3.16)–(3.19) provide a simplified version of the PLUS algorithm dealing with the one-at-a-time scenario in which every intermediate turning point in the PLUS path is the intersection of exactly two line segments of positive length. Although the one-at-a-time condition holds almost everywhere, numerical ties do occur in applications. When the one-at-a-time condition fails, the main branch (2.7) is a limit path of one-at-a-time paths, so that it is a graph with no dead end. The difference here is that when the PLUS path reaches a more-than-two-way intersection, say at step $k$, it must checked the indicators $\boldsymbol{\eta}^{(\ell)}, 0 \le \ell < k$, to avoid infinite looping with the covered segments . The computational cost for checking the indicators is $O(k)$ if $\boldsymbol{\eta}$ are efficiently coded, which is small compared with the cost $O(np)$ for finding the exit time (3.15). See Zhang (2007b) for details.

7.3. *Orthonormal designs and more discussion on penalties.* For orthonormal designs $\mathbf{x}_j'\mathbf{x}_k/n = I\{j = k\}$, the penalized estimation problem is reduced to the case of $p = 1$. For $\rho(t; \lambda) = \lambda^2 \rho_m(t/\lambda)$ with the quadratic spline penalties (3.1),

$$(7.1) \quad \widehat{\beta}_j = \lambda b(\mathbf{x}_j'\mathbf{y}/(n\lambda)) \qquad \text{where } b(z) \equiv \arg\min_b \{(z - b)^2/2 + \rho_m(|b|)\}.$$

For $p = 1$ and the MCP with $\kappa(\rho_2) = 1/\gamma < 1$, the solution of (7.1) is

$$b_f(z) = \text{sgn}(z) \min\{|z|, \gamma(|z| - \lambda)_+/(\gamma - 1)\},$$

which turns out to be the firm threshold estimator of Gao and Bruce (1997). The firm threshold estimator is always between the soft threshold estimator $b_s(z) \equiv \text{sgn}(z)(|z| - \lambda)_+$ and the hard threshold estimator $b_h(z) \equiv zI\{|z| > \lambda\}$. Actually, $b_s(z) \le b(z) \le b_f(z) \le b_h(z)$ for $z > 0$ and the opposite inequalities hold for $z < 0$ for all solutions of (7.1), given a fixed $\gamma\lambda$ in (2.3) or a fixed maximum concavity $\kappa(\rho_m) = 1/\gamma$ with $\gamma > 1$. We plot these univariate estimators in Figure 9 along with the univariate SCAD estimator. For $p = 1$ and $\kappa(\rho_2) = 1/\gamma \geq 1$, the MC+ path (2.7) has three segments and (2.8), identical to the hard threshold estimator, globally minimizes the penalized loss. See Figure 9 on the left. Antoniadis and Fan (2001) observed that in the orthonormal case, the global minimizer (7.1) for the penalty (2.1) with $\gamma = 1/2$ yields the hard threshold estimator. In fact, in the univariate case, any penalty function with concave derivative $\dot{\rho}(t; \lambda)$ and $\gamma \leq 1$ in (2.3) yields the hard threshold estimator as the global minimizer in (7.1).

The analytical and computational properties of penalized estimation and selection for general correlated $\mathbf{X}$ and concave penalty is much more complicated than the case of $p = 1$, since they are determined in many ways by the interplay between the penalty and the design. To a large extent, the effects of the penalty can be summarized by the threshold factor $\gamma$ for the unbiasedness in (2.3), the maximum concavity $\kappa(\rho; \lambda)$ in (2.2) and their relationships to the correlations of the design vectors. This naturally leads to our choice of the MCP as the minimizer of $\kappa(\rho; \lambda)$ given the threshold factor $\gamma$ and the role of $\gamma = 1/\kappa(\rho_1)$ as the regularization parameter for the bias and computational complexity of the MC+.

### 7.4. Adaptive penalty.

The PLUS algorithm applies to the penalized loss

$$(7.2) \qquad (2n)^{-1}\|\mathbf{y} - \mathbf{X}\boldsymbol{\beta}\|^2 + \sum_{j=1}^{p} \lambda^2 \rho_m(|\beta_j|r_j/\lambda), \qquad r_j > 0 \;\forall j,$$

through the scale change $\{\mathbf{x}_j, \beta_j\} \to \{\mathbf{x}_j r_j, \beta_j/r_j\}$. It can be easily modified to accommodate different quadratic $\rho_m$ of the form (3.1) for different $j$. For example, different $\gamma = \gamma_j$ can be used with the MC+, so that the $j$th path $\widehat{\beta}_j(\lambda)$ reaches the unbiased region when $|\widehat{\beta}_j(\lambda)|r_j/\lambda \geq \gamma_j$. This allows $r_j$ and $\gamma_j$ to be data dependent. For $r_j = 1$, the unbiasedness condition $\gamma_j \lambda \leq |\beta_j|$ allows a higher level of convexity than (2.9) does.

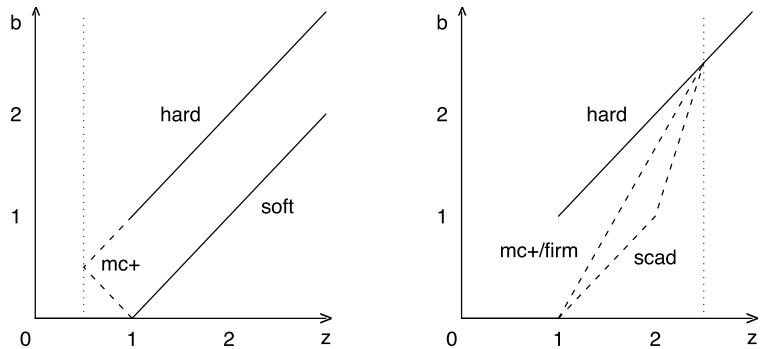

Fig. 9. *Left: the univariate hard, soft and MC+ paths in $z \oplus b \in H \oplus H^* = \mathbb{R}^2$ with a vertical dotted line at $z = \gamma = 1/2$. Right: the hard, MC+/firm and SCAD paths for $p = 1$ with $\gamma = 5/2$. Hard and soft path in solid, and additional segments of MC+ and SCAD in dashed lines.*

Zou (2006) proposed an adaptive LASSO with $\lambda^2 \rho_1(|\beta_j|r_j/\lambda) = \lambda r_j|\beta_j|$, where $r_j$ is a decreasing function of an initial estimate of $\beta_j$. The idea is to reduce the penalty level or the bias for large/nonzero $|\beta_j|$, but its effectiveness for selection consistency essentially requires the initial estimator to be larger than a (possibly unspecified and random) threshold for most large/nonzero $|\beta_j|$ and smaller than the same threshold for most small/zero $|\beta_j|$. This approach was proven for bounded $p = \mathrm{rank}(\mathbf{X})$ to provide selection consistency (1.4) in Zou (2006) and Zou and Li (2008). Marginal regression $\mathbf{x}_j \mathbf{y}/\|\mathbf{x}_j\|^2$ can be used as an initial estimate of $\beta_j$ and is proved to result in the selection consistency of the adaptive LASSO under a certain partial orthogonality condition on the pairwise correlations among vectors $\{\mathbf{y}, \mathbf{x}_1, \ldots, \mathbf{x}_p\}$ [Huang, Ma and Zhang (2008)].

### 7.5. General loss functions.

Consider the general penalized loss $L(\boldsymbol{\beta}; \lambda) \equiv \psi(\boldsymbol{\beta}) + \sum_{j=1}^{p} \rho(|\beta_j|; \lambda)$, where $\psi(\boldsymbol{\beta}) \equiv \psi_n(\boldsymbol{\beta}; \mathbf{X}, \mathbf{y})$ is a convex function of $\boldsymbol{\beta} \in \mathbb{R}^p$ given data $(\mathbf{X}, \mathbf{y})$. In generalized linear models, $n\psi_n(\boldsymbol{\beta}; \mathbf{X}, \mathbf{y})$ is the negative log-likelihood. With $(\dot{\psi}_j)_{p \times 1}$ and $(\ddot{\psi}_{j\ell})_{p \times p}$ being the gradient and Hessian of $\psi$, (2.7) must satisfy

$$(7.3) \quad \begin{cases} \dot{\psi}_j(\widehat{\boldsymbol{\beta}}^{(x)}) + \mathrm{sgn}(\widehat{\beta}_j^{(x)})\dot{\rho}(|\widehat{\beta}_j^{(x)}|; \lambda^{(x)}) = 0, & \widehat{\beta}_j^{(x)} \neq 0, \\ |\dot{\psi}_j(\widehat{\boldsymbol{\beta}}^{(x)})| \leq \lambda^{(x)}, & \widehat{\beta}_j^{(x)} = 0. \end{cases}$$

Let $\mathbf{s}^{(x)} \equiv d\widehat{\boldsymbol{\beta}}^{(x)}/d\lambda^{(x)}$ and $\widehat{A}^{(x)} \equiv \{j : \widehat{\beta}_j^{(x)} \neq 0\}$. Differentiation of (7.3) yields

$$(7.4) \quad \left\{ \sum_{\ell \in \widehat{A}^{(x)}} \ddot{\psi}_{j\ell}(\widehat{\boldsymbol{\beta}}^{(x)}) s_\ell^{(x)} \right\} + \ddot{\rho}(|\widehat{\beta}_j^{(x)}|; \lambda^{(x)}) s_j^{(x)} = a(\widehat{\beta}_j^{(x)}; \lambda^{(x)})$$

for $j \in \widehat{A}^{(x)}$ and $s_j^{(x)} = 0$ for $j \notin \widehat{A}^{(x)}$, where $a(t; \lambda) = -\mathrm{sgn}(t)(\partial/\partial\lambda)\dot{\rho}(|t|; \lambda)$. This provides the local direction of the next move and thus allows an extension of the PLUS algorithm. The main difference of such an extension from (3.8) is that the step size has to be small when $\psi(\cdot)$ is not a quadratic spline. The main difference of such an extension from the computation of the LASSO for the generalized linear models [Genkin, Lewis and Madigan (2004), Zhao and Yu (2007) and Park and Hastie (2007)] is the possibility of the sign change $d\lambda^{(x)}/dx$ to allow the path to traverse from one local minimum to another. Extensions of the LARS with large step size $\Delta^{(k)}$ have been considered by Rosset and Zhu (2007) for support vector machine and by Zhang (2007a) for continuous generalized gradient descent.

7.6. *Sub-Gaussian errors.* Remark 3 in Section 2 mentions the validity of our theorems when the normality condition $\boldsymbol{\varepsilon} \sim N(0, \sigma^2 \mathbf{I}_n)$ in (1.2) is replaced by a sub-Gaussian condition on the error vector. Here, we provide some details.

PROPOSITION 3. *Let $\boldsymbol{\varepsilon} \in \mathbb{R}^n$ be a random vector satisfying the sub-Gaussian condition $E \exp(\mathbf{x}' \boldsymbol{\varepsilon}) \leq e^{\|\mathbf{x}\|^2 \sigma_1^2/2}$ for all $\mathbf{x} \in \mathbb{R}^n$. Then for projections $\mathbf{P}$ of rank $m$*

$$P\left\{ \frac{\|\mathbf{P}\boldsymbol{\varepsilon}\|^2}{m\sigma_1^2} \geq \frac{1+x}{\{1 - 2/(e^{x/2}\sqrt{1+x} - 1)\}_+^2} \right\} \leq e^{-mx/2}(1+x)^{m/2} \qquad \forall x > 0.$$

The normality condition is used in our proofs only to provide upper bounds for the tail probabilities of $\mathbf{u}' \boldsymbol{\varepsilon}$ and $\|\mathbf{P}\boldsymbol{\varepsilon}\|^2/m$. The sub-Gaussian condition implies $P\{\mathbf{u}' \boldsymbol{\varepsilon}/\sigma_1 > t\} \leq e^{-t^2/2} \leq (t+1/t)\Phi(-t)$ for $t > 0$ and $\|\mathbf{u}\| = 1$, comparable to the normal tail probability. Proposition 3 is comparable to the $\chi_m^2/m$ tail probability bound needed in our proofs.

## APPENDIX

In this appendix, we provide all the proofs. Theorem 1 is a special case of Theorem 5 and Theorem 2 concerns estimation in the same special case. The proof of Theorem 5 requires Theorem 6 and the proof of Theorem 7 requires Theorem 9 and Proposition 2. Thus, the proofs are given in the following order: Theorems 3, 4, 6, 5, 1, 2 and 9, Proposition 2, Theorems 7 and 8 and then Proposition 3. Two lemmas, needed in the proof of Theorems 6, 5 and 2, are stated before the proof of Theorem 6 and proved at the end of the Appendix.

PROOF OF THEOREM 3. Let $\mathbf{X}$ be fixed. Define $d_k(\boldsymbol{\eta}) \equiv \#\{j : |\eta_j| = k\}$, $k = 1, 2$. We consider three types of indicators $\boldsymbol{\eta} \in \{-2, -1, 0, 1, 2\}^p$ with $\boldsymbol{\eta} = \mathbf{0}$ as type-1.

Type-2: $d_2(\boldsymbol{\eta}) \geq n \wedge p$. Let $(\tau \widetilde{\mathbf{z}}) \oplus \mathbf{b} \in S(\boldsymbol{\eta})$ as in (3.5), so that (3.3) holds with $z_j = \tau \widetilde{z}_j = \tau \mathbf{x}_j' \mathbf{y}/n$. Since $\dot{\rho}_2(|b_j|) = 0$ for $|\eta_j| = 2$, (3.3) implies $\mathbf{x}_j'(\tau \mathbf{y} - \mathbf{X}\mathbf{b}) = 0$ for all $|\eta_j| = 2$. Since $\tau \mathbf{y} - \mathbf{X}\mathbf{b} \in \mathbb{R}^n$ and $\{\mathbf{x}_j, |\eta_j| = 2\}$ contains at least $n \wedge p$ linearly independent vectors, by (3.6)

$$\text{(A.1)} \qquad \begin{cases} d_2(\boldsymbol{\eta}) \geq n \wedge p \\ (\tau \widetilde{\mathbf{z}}) \oplus \mathbf{b} \in \ell(\boldsymbol{\eta}|\widetilde{\mathbf{z}}) \end{cases} \Rightarrow \quad \mathbf{X}'(\tau \mathbf{y} - \mathbf{X}\mathbf{b}) = \mathbf{0}.$$

Type-3: $d_2(\boldsymbol{\eta}) < n \wedge p$ and $\boldsymbol{\eta} \neq \mathbf{0}$. If (3.3) holds for $\mathbf{z} = \mathbf{0}$, then $b_j \mathbf{x}_j' \mathbf{X}\mathbf{b}/n = b_j \boldsymbol{\chi}_j' \mathbf{b} = -|b_j| \dot{\rho}_2(|b_j|)$ for all $b_j \neq 0$, so that $\|\mathbf{X}\mathbf{b}\|^2/n = -\sum_j |b_j| \dot{\rho}_2(|b_j|) = 0$ due to $\dot{\rho}_2(|b_j|) \geq 0$. Since $\dot{\rho}_2(|b_j|) = (1 - |b_j|/\gamma)_+ > 0$ for $|b_j| < \gamma$ and $|b_j| \leq \gamma$ for $|\eta_j| = 1$, $b_j$ equals either $0$ or $\gamma \eta_j$ for $|\eta_j| = 1$ in such cases. Therefore, $\mathbf{X}\mathbf{b} = \sum_{|\eta_j|=2} b_j \mathbf{x}_j + \gamma \sum_{b_k \neq 0, |\eta_k| < 2} \eta_k \mathbf{x}_k = 0$. This is impossible for nondegenerate $\mathbf{X}$ since $\gamma > 0$ and $d_2(\boldsymbol{\eta}) < n \wedge p$. Thus, $\mathbf{0} \oplus \mathbf{b} \notin S(\boldsymbol{\eta})$ for indicators $\boldsymbol{\eta}$ of type-3.

We now consider the choice of $\gamma$ for the MCP. It follows from (3.2) and (3.11) that the determinant $\det(\mathbf{Q}(\boldsymbol{\eta}))$ is a polynomial of $v_1 = 1/\gamma$ with $\det(\boldsymbol{\Sigma}_{j\,:\,|\eta_j|=2})(-v_1)^{d_1(\boldsymbol{\eta})}$ as the leading term and $\det(\boldsymbol{\Sigma}_{j\,:\,|\eta_j|=2}) \neq 0$ for type-3 $\boldsymbol{\eta}$ by (3.20). Let $\Gamma_0(\mathbf{X})$ be the finite set of all reciprocals of the real roots of such polynomials with type-3 $\boldsymbol{\eta}$. We choose $\gamma \notin \Gamma_0(\mathbf{X})$ hereafter, so that $\det(\mathbf{Q}(\boldsymbol{\eta})) \neq 0$ for all $\boldsymbol{\eta}$ of type-3. Since $\det(\mathbf{Q}(\boldsymbol{\eta})) \neq 0$, in $S(\boldsymbol{\eta})$ the vector $(b_j, \eta_j \neq 0)'$ is a linear function of $\mathbf{z}$ by (3.12), so that by (3.6) and the discussion in the previous paragraph

$$(A.2) \quad \begin{cases} d_2(\boldsymbol{\eta}) < n \wedge p, \\ \boldsymbol{\eta} \neq \mathbf{0}, \end{cases} \quad \Rightarrow \quad \begin{cases} \det(\mathbf{Q}(\boldsymbol{\eta})) \neq 0, \\ \ell(\boldsymbol{\eta}|\mathbf{z}) \text{ is a generalized line segment}, \\ \mathbf{0} \oplus \mathbf{b} \notin \ell(\boldsymbol{\eta}|\mathbf{z}) \ \forall \mathbf{b}. \end{cases}$$

Here, a generalized line segment includes the empty set, single points in $H \oplus H^* = \mathbb{R}^{2p}$, and line segments of finite or infinite length.

For each nonzero $\mathbf{z} \in H \equiv \mathbb{R}^p$, we define a graph $G(\mathbf{z})$ with $\ell(\boldsymbol{\eta}|\mathbf{z})$ of positive length and type-3 $\boldsymbol{\eta}$ as edges and the end points of edges as vertices. The graph $G(\mathbf{z})$ is not necessarily connected. A vertex in $G(\mathbf{z})$ is terminal if it is also a boundary point of $S(\boldsymbol{\eta})$ for some $\boldsymbol{\eta}$ of type-2. If the PLUS path reaches a terminal vertex $(\tau \widetilde{\mathbf{z}}) \oplus \mathbf{b}$, then $\mathbf{b}/\tau$ provides an optimal fit by (A.1). The degree of a vertex in $G(\mathbf{z})$ is the number of edges connected to it.

Suppose $\widetilde{\mathbf{z}} \neq \mathbf{0}$. At step $k = 0$, the MC+ path reaches $(\tau^{(0)}\widetilde{\mathbf{z}}) \oplus \mathbf{b}^{(0)}$ as a boundary point of $S(\mathbf{0})$. Since the $p$-parallelepipeds (3.5) are contiguous, $(\tau^{(0)}\widetilde{\mathbf{z}}) \oplus \mathbf{b}^{(0)}$ is also a boundary point of $S(\boldsymbol{\eta}^{(1)})$ for some $\boldsymbol{\eta}^{(1)}$ satisfying either (A.1) or (A.2) with $\mathbf{z} = \widetilde{\mathbf{z}}$. If $\boldsymbol{\eta}^{(1)}$ is of type-2, then $\mathbf{b}^{(0)}/\tau^{(0)}$ gives an optimal fit and the MC+ path ends with $k^* = 0$. Otherwise, the MC+ path enters the graph $G(\widetilde{\mathbf{z}})$ at the initial vertex $(\tau^{(0)}\widetilde{\mathbf{z}}) \oplus \mathbf{b}^{(0)}$. If the degree of the initial vertex is odd and the degrees of all other nonterminal vertices are even, then the MC+ path traverses through $G(\widetilde{\mathbf{z}})$ and eventually reaches a terminal vertex in one pass. This is simply an Euler's Konigsberg problem.

Let $S_0$ be the union of all intersections of three or more distinct $p$-parallelepipeds $S(\boldsymbol{\eta})$, $\boldsymbol{\eta} \in \{-2, -1, 0, 1, 2\}^p$, and $H_0 \equiv \{\mathbf{z} : (\tau \mathbf{z}) \oplus \mathbf{b} \in S_0$ for some $\tau$ and $\mathbf{b}\}$. Since the interiors of the $p$-parallelepipeds $S(\boldsymbol{\eta})$ do not intersect, the intersections of three distinct $S(\boldsymbol{\eta})$ are $(p-2)$-parallelepipeds, so that the projection of $S_0$ to the $(p-1)$-sphere $\{\mathbf{z} : \|\mathbf{z}\| = 1\}$ along the rays $\{\tau \mathbf{z}, \tau > 0\}$ has Haar measure zero. Consequently, $H_0$ has Lebesgue measure zero in $H \equiv \mathbb{R}^p$.

For $\mathbf{z} \notin H_0$, each vertex in $G(\mathbf{z})$ is a boundary point of exactly two $p$-parallelepipeds $S(\boldsymbol{\eta})$, so that the initial vertex has degree 1 and other nonterminal vertices have degree 2 in $G(\mathbf{z})$. Thus, the initial vertex is connected to a terminal vertex in $G(\widetilde{\mathbf{z}})$ in a unique way for $\widetilde{\mathbf{z}} \notin H_0$, and the conclusions of part (i) holds by (A.2).

For $\widetilde{\mathbf{z}} \in H_0$, the initial vertex is still connected to at lease one terminal vertex in $G(\widetilde{\mathbf{z}})$ since $H_0^c$ is dense in $H \equiv \mathbb{R}^p$, and the limits of $G(\mathbf{z})$ as $\mathbf{z} \to \widetilde{\mathbf{z}}$ are subgraphs of $G(\widetilde{\mathbf{z}})$. Hence, the conclusion of part (i) hold in either cases.

Parts (ii) and (iii) hold since $\det(\mathbf{Q}(\boldsymbol{\eta})) \neq 0$ almost everywhere for fixed $\gamma$ and type-3 $\boldsymbol{\eta}$. For part (iv), we consider $\widetilde{\mathbf{z}} \notin H_0$. We have already proved the uniqueness of the graph and that the path ends with perfect fit at a type-2 $\boldsymbol{\eta}$. Since the vertex $(\tau^{(k-1)}\widetilde{\mathbf{z}}) \oplus \mathbf{b}^{(k-1)}$ is a boundary point of exactly two $p$-parallelepipeds $S(\boldsymbol{\eta}^{(k-1)})$ and $S(\boldsymbol{\eta}^{(k)})$, $j^{(k-1)}$ in (3.9) uniquely indicates the side of the intersection. Since the edges must pass through the interior of the $p$-parallelepipeds, $\boldsymbol{\eta}^{(k)}$ and $\xi^{(k)}$ are given by (3.10) and (3.13). The slope $\mathbf{s}^{(k)}$ is uniquely determined by (3.12) due to $\det(\mathbf{Q}(\boldsymbol{\eta}^{(k)})) \neq 0$. The hitting time $\Delta_j$ in (3.14) is computed from the current position $(\tau^{(k-1)}\widetilde{\mathbf{z}}) \oplus \mathbf{b}^{(k-1)}$, the slope $\mathbf{s}^{(k)}$ and the inequalities for the boundary of the $p$-parallelepiped $S(\boldsymbol{\eta}^{(k)})$ in (3.5). Since the length of the edge is positive, $\Delta^{(k)} > 0$ in (3.15). Since the path does not return to $\mathbf{0}$, $\tau^{(k)} > 0$ in (3.15).

For part (v), $(\widetilde{\mathbf{z}} \oplus \widehat{\boldsymbol{\beta}}^{(x)})/\lambda$ is in the interior of $S(\boldsymbol{\eta}^{(k)})$ at $\lambda = \lambda^{(x)}$ for $k - 1 < x < k$, so that (3.3) and thus (2.6) hold with strict inequality. By (3.11),

$$
\begin{aligned}
(\partial/\partial t)L&(\widehat{\boldsymbol{\beta}}^{(x)} + t\mathbf{b}; \lambda) \\
&= t\mathbf{b}_1'\mathbf{Q}(\boldsymbol{\eta}^{(k)})\mathbf{b}_1 + \sum_{\eta_j^{(k)} = 0} |b_j|\{\lambda - \operatorname{sgn}(b_j)\mathbf{x}_j'(\mathbf{y} - \mathbf{X}\widehat{\boldsymbol{\beta}}^{(x)})/n + O(t)\}
\end{aligned}
$$

is positive for small $t > 0$, where $\mathbf{b}_1 = (b_j, \eta_j^{(k)} \neq 0)'$. Thus, $\widehat{\boldsymbol{\beta}}^{(x)}$ is a local minimizer. $\square$

PROOF OF THEOREM 4. Since $\widehat{\boldsymbol{\beta}}^o$ is the oracle LSE, $\mathbf{x}_j'(\mathbf{y} - \mathbf{X}\widehat{\boldsymbol{\beta}}^o) = 0$ for $j \in A^o$. If $|\widehat{\beta}_j^o| \geq \lambda\gamma$, then $\dot{\rho}(|\widehat{\beta}_j^o|; \lambda) = 0$ by (2.3). Thus, $\widehat{\boldsymbol{\beta}}^o$ is a solution of (2.6) and $\operatorname{sgn}(\widehat{\boldsymbol{\beta}}^o) = \operatorname{sgn}(\boldsymbol{\beta})$ for all $\lambda_1 \leq \lambda \leq \lambda_2$ in the intersection of

$$
\begin{aligned}
(A.3) \qquad \Omega_1^o(\lambda_1) &\equiv \Big\{\max_{j \notin A^o}|\mathbf{x}_j(\mathbf{y} - \mathbf{X}\widehat{\boldsymbol{\beta}}^o)/n| < \lambda_1\Big\}, \\
\Omega_2^o(\lambda_2) &\equiv \Big\{\min_{j \in A^o} \operatorname{sgn}(\beta_j)\widehat{\beta}_j^o > \gamma\lambda_2\Big\}.
\end{aligned}
$$

Moreover, since the solution of (2.6) is unique, $\widehat{\boldsymbol{\beta}}^o = \widehat{\boldsymbol{\beta}}(\lambda)$ for all $\lambda_1 \leq \lambda \leq \lambda_2$ in this case.

Let $\mathbf{P}_1^o$ be the orthogonal projection from $\mathbb{R}^n$ to the linear span of $\{\mathbf{x}_j, j \in A^o\}$. Since $\mathbf{y} - \mathbf{X}\widehat{\boldsymbol{\beta}}^o = (\mathbf{I}_n - \mathbf{P}_1^o)\boldsymbol{\varepsilon}$, $\mathbf{x}_j'(\mathbf{y} - \mathbf{X}\widehat{\boldsymbol{\beta}}^o)/n$ are normal variables with zero mean and variance bounded by $\sigma^2 \|\mathbf{x}_j\|^2/n^2$, so that $1 - P\{\Omega_1^o(\lambda_1)\} \leq \pi_{n,1}(\lambda_1)$. By (2.10) and (4.1), $\widehat{\beta}_j^o \sim N(\beta_j, \sigma^2 w_j^o/n)$ for all $j \in A^o$. Since $|\beta_j| \geq \beta_* \geq \gamma\lambda_2$, we have $1 - P\{\Omega_2^o(\lambda)\} \leq \pi_{n,2}(\lambda_2)$. Inequality (4.3) follows by combining the above two probability bounds. □

Let us state the two lemmas. For $m \geq 1$ and $B \subset \{1, \ldots, p\}$, define semi-norms

$$\text{(A.4)} \quad \zeta(\mathbf{v}; m, B) \equiv \max\left\{\frac{\|(\mathbf{P}_A - \mathbf{P}_B)\mathbf{v}\|}{(mn)^{1/2}} : B \subseteq A \subseteq \{1, \ldots, p\}, |A| = m + |B|\right\}$$

for $\mathbf{v} \in \mathbb{R}^n$, where $\mathbf{P}_A$ is the orthogonal projection from $\mathbb{R}^n$ to the span of $\{\mathbf{x}_j : j \in A\}$.

LEMMA 1. *Suppose (2.11) holds for $\mathbf{X}$ with certain $d^*$ and $c^* \geq c_* \geq \kappa \geq 0$. Let $K_*$ be as in (4.5) with an $\alpha \in (0,1)$, and $B \subset \{1, \ldots, p\}$ with $|B| \leq d^*/(K_* + 1)$. Let $\lambda > 0$ be fixed and $\rho(t; \lambda)$ be a penalty satisfying $\lambda(1 - \kappa t/\lambda)_+ \leq \dot{\rho}(t; \lambda) \leq \lambda$ for all $t > 0$. Let $1 \leq m \leq m^* \equiv d^* - |B|$ and $\mathbf{y} \in \mathbb{R}^n$ with $(\sqrt{c^*}/\alpha)\zeta(\mathbf{y}; m, B) \leq \lambda$, where $\zeta(\cdot; m, B)$ is as in (A.4). Let $\lambda \oplus \widehat{\boldsymbol{\beta}}$ be a solution of (2.6), $B \cup \{j : \widehat{\beta}_j \neq 0\} \subseteq A_1 \subseteq B \cup \{j : |\mathbf{x}_j'(\mathbf{y} - \mathbf{X}\widehat{\boldsymbol{\beta}})/n| = \dot{\rho}(|\widehat{\beta}_j|; \lambda)\}$ and $\widehat{\boldsymbol{\beta}}^o$ be as in (2.10). If $|A_1| = |B| + m$, then*

$$\text{(A.5)} \quad |A_1| - |B| < K_* \sum_{j \in B} \dot{\rho}^2(|\widehat{\beta}_j|; \lambda)/\lambda^2 \leq K_*|B|.$$

*If $\lambda \geq (\sqrt{c^*}/\alpha)\zeta(\mathbf{y}; m^*, B)$, then $\#\{j \notin B : |\mathbf{x}_j'(\mathbf{y} - \mathbf{X}\widehat{\boldsymbol{\beta}})/n| = \dot{\rho}(|\widehat{\beta}_j|; \lambda)\} < 1 \vee (K_*|B|)$ and*

$$\text{(A.6)} \quad \begin{aligned} c_* \|\widehat{\boldsymbol{\beta}} - \widehat{\boldsymbol{\beta}}^o\| &\leq \sqrt{c_*/n}\|\mathbf{X}(\widehat{\boldsymbol{\beta}} - \widehat{\boldsymbol{\beta}}^o)\| \\ &\leq \left\{\sum_{j \in B} \dot{\rho}^2(|\widehat{\beta}_j|; \lambda)\right\}^{1/2} + \alpha\lambda\sqrt{K_*|B|c_*/c^*}. \end{aligned}$$

LEMMA 2. *Let $\zeta(\mathbf{v}; m, B)$ be as in (A.4) with deterministic $m$ and $B$. Let $\widetilde{p}_\epsilon \geq \sqrt{e}$ be the solution of (2.12) with $d^o = |B|$. Suppose $\boldsymbol{\varepsilon} \sim N(0, \sigma^2 I_n)$. Then*

$$\text{(A.7)} \quad P\{\zeta(\boldsymbol{\varepsilon}; m, B) \geq \sigma\sqrt{(2/n)\log\widetilde{p}_\epsilon}\} \leq \frac{\epsilon e^{\mu^2/2}\Phi(-\mu)}{\sqrt{\log\widetilde{p}_\epsilon}} \leq \frac{\epsilon/2}{\sqrt{\log\widetilde{p}_\epsilon}} \leq \frac{\epsilon}{\sqrt{2}},$$

*where $\mu = \{2\log\widetilde{p}_\epsilon - 1 + 1/m\}\sqrt{m}/\sqrt{2\log\widetilde{p}_\epsilon}$.*

PROOF OF THEOREM 6. Let $d_1^{(x)} \equiv \#\{j : j \in B \text{ or } |\mathbf{x}_j'(\mathbf{y} - \mathbf{X}\widehat{\boldsymbol{\beta}}^{(x)})/n| = \dot{\rho}(|\widehat{\beta}_j^{(x)}|; \lambda^{(x)})\}$ and $x_1 = \inf\{x \geq 0 : \lambda^{(x)} < \lambda_1 \text{ or } \lambda^{(x)} < (\sqrt{c^*}/\alpha)\zeta(\mathbf{y}; m, B)\}$

with the $\lambda^{(x)} \oplus \widehat{\boldsymbol{\beta}}^{(x)}$ in (2.7) and $m = d^* - |B|$. We first prove $d_1^{(x)} < d^*$ for $0 \leq x \leq x_1$. Let $A_1^{(x)}$ be any set satisfying

(A.8)
$$B \cup \{j : \widehat{\beta}_j^{(x)} \neq 0\}$$
$$\subseteq A_1^{(x)} \subseteq B \cup \{j : |\mathbf{x}_j'(\mathbf{y} - \mathbf{X}\widehat{\boldsymbol{\beta}}^{(x)})/n| = \dot{\rho}(|\widehat{\beta}_j^{(x)}|; \lambda^{(x)})\}.$$

By (2.6), the left-hand side is always a subset of the right-hand side in (A.8). Moreover, since $\widehat{\boldsymbol{\beta}}^{(x)}$ is continuous in $x$, $\operatorname{sgn}(\widehat{\beta}_j^{(x-)}) = \operatorname{sgn}(\widehat{\beta}_j^{(x)}) = \operatorname{sgn}(\widehat{\beta}_j^{(x+)})$ fails to hold only if $\widehat{\beta}_j^{(x)} = 0$ and $|\mathbf{x}_j'(\mathbf{y} - \mathbf{X}\widehat{\boldsymbol{\beta}}^{(x)})/n| = \dot{\rho}(0; \lambda^{(x)}) = \lambda^{(x)}$, so that we are allowed to add variables to $|A_1^{(x)}|$ one-at-a-time. Thus, since $\widehat{\boldsymbol{\beta}}^{(0)} = \mathbf{0}$, if $d_1^{(x_2)} \geq d^*$ for some $0 \leq x_2 \leq x_1$, there must be a choice of $A_1^{(x)}$ with $|A_1^{(x)}| = d^*$ and $0 \leq x \leq x_2$. On the other hand, it follows from Lemma 1 that $\lambda^{(x)} \geq \lambda_1 \vee \{(\sqrt{c^*}/\alpha)\zeta(\mathbf{y}; m^*, B)\}$ and $|B| < |A_1^{(x)}| = d^*$ imply $|A_1^{(x)}| < (K_* + 1)|B| \leq d^*$, where $m^* = m$. Thus, $|B| < |A_1^{(x)}| = d^*$ can never be attained for $0 \leq x \leq x_1$. It follows that $\#\{j \notin B : \widehat{\beta}_j^{(x)} \neq 0\} \leq |A_1^{(x)}| - |B| < 1 \vee (K_*|B|)$ for all $0 \leq x \leq x_1$ by Lemma 1.

Let $\lambda_4 = \sigma\sqrt{(2/n)\log \widetilde{p}_\epsilon} + \theta_B/\sqrt{m}$. By (2.8), $\widehat{\lambda} \oplus \widehat{\boldsymbol{\beta}} = \lambda^{(x)} \oplus \widehat{\boldsymbol{\beta}}^{(x)}$ with a certain $\lambda^{(x)} \geq \lambda_1 \vee (\lambda_4\sqrt{c^*}/\alpha)$, so that the left-hand side of (4.8) is no greater than $P\{\Omega_4^c\}$ with $\Omega_4 = \{\zeta(\mathbf{y}; m, B) \leq \lambda_4\}$. Since $\zeta(\mathbf{X}\boldsymbol{\beta}; m, B) \leq \|(\mathbf{I}_n - \mathbf{P}_B)\mathbf{X}\boldsymbol{\beta}\|/\sqrt{nm}$ by (A.4) and $\theta_B \equiv \|\mathbf{X}\boldsymbol{\beta} - \mathbf{X}E\widehat{\boldsymbol{\beta}}^o\| = \|(\mathbf{I}_n - \mathbf{P}_B)\mathbf{X}\boldsymbol{\beta}\|$ by (2.10), $\zeta(\mathbf{y}; m, B) \leq \zeta(\boldsymbol{\varepsilon}; m, B) + \theta_B/\sqrt{m}$. Thus, $P\{\Omega_4^c\} \leq P\{\zeta(\boldsymbol{\varepsilon}; m, B) > \sigma\sqrt{(2/n)\log \widetilde{p}_\epsilon}\}$. The conclusion follows from Lemma 2. $\square$

PROOF OF THEOREM 5. Consider the event $\Omega = \bigcap_{j=1}^3 \Omega_j(\lambda_j)$, where $\Omega_j(\lambda_j), j = 1, 2,$ are as in (A.3) and $\Omega_3(\lambda_3) \equiv \{\zeta(\boldsymbol{\varepsilon}; m, A^o) \leq \lambda_3\}$. It follows from the proof of Theorem 4 that $\lambda \oplus \widehat{\boldsymbol{\beta}}^o$ is a solution of (2.6) for all $\lambda_1 \leq \lambda \leq \lambda_2$. Since $\kappa(\rho; \lambda_2) \leq \kappa < c_*$, the sparse convex condition (2.5) holds with rank $d^*$, so that $\lambda \oplus \widehat{\boldsymbol{\beta}}^o$ is the unique solution of (2.6) subject to $\lambda_1 \leq \lambda \leq \lambda_2$ and $\#\{j : |\beta_j| + |\widehat{\beta}_j| > 0\} \leq d^*$. Since $\zeta(\mathbf{y}; m, B) = \zeta(\boldsymbol{\varepsilon}; m, A^o)$ with $B = A^o$ in (A.4), we also have $\lambda_2 \geq (\sqrt{c^*}/\alpha)\lambda_3 \geq (\sqrt{c^*}/\alpha)\zeta(\mathbf{y}; m, B)$ in $\Omega$.

In the event $\Omega$, consider the path $\lambda^{(x)} \oplus \widehat{\boldsymbol{\beta}}^{(x)}$ with $0 \leq x \leq x_1 \equiv \inf\{x : \lambda^{(x)} < \lambda_1\}$. Let $A_1^{(x)}$ be as in (A.8). If $|A_1^{(x)}| = d^*$ and $\lambda^{(x)} \geq \lambda_2$, then Lemma 1 provide $|A_1^{(x)}| < d^*$. If $|A_1^{(x)}| = d^*$ and $\lambda_1 \leq \lambda^{(x)} \leq \lambda_2$, then the uniqueness of $\lambda \oplus \widehat{\boldsymbol{\beta}}^o$ implies $\widehat{\boldsymbol{\beta}}^{(x)} = \widehat{\boldsymbol{\beta}}^o$. Since $|\mathbf{x}_j'(\mathbf{y} - \mathbf{X}\widehat{\boldsymbol{\beta}}^o)/n| < \lambda_1 \leq \dot{\rho}(0; \lambda^{(x)})$ for $j \notin A^o$, we have $A_1^{(x)} = A^o$. Thus, $|A_1^{(x)}| = d^*$ can never be attained with $0 \leq x \leq x_1$, and $\widehat{\boldsymbol{\beta}}(\lambda) = \widehat{\boldsymbol{\beta}}^o$ for all $\lambda_1 \leq \lambda \leq \lambda_2$ in the event $\Omega$.

We still need to bound $1 - P\{\Omega\}$. The proof of Theorem 4 provides $1 - P\{\Omega_j(\lambda_j)\} \leq \pi_{n,j}(\lambda_j)$ for $j = 1, 2$. By (A.4), $\zeta(\varepsilon; m^*, A^o)$ is the maximum of $\binom{p - d^o}{m} \chi_m^2$ variables, so that $1 - P\{\Omega_3(\lambda_3)\} \leq \pi_{n,3}(\lambda_3)$. Thus, (4.6) holds. Finally, (4.7) follows from (4.6) with applications of the inequality $e^{t^2/2}\Phi(-t) \leq \min\{1/2, 1/(t\sqrt{2\pi})\}$ and Lemma 2. $\square$

PROOF OF THEOREM 1. Theorem 1 follows from Theorem 5 with $\alpha = 1/2$, since $\gamma = 1/\kappa \geq c_*^{-1}\sqrt{4 + c_*/c^*}$ implies $K_* + 1 \leq c^*/c_* + 1/2$ in (4.5) as in Remark 5. $\square$

PROOF OF THEOREM 2. As in the proof of Theorem 1, we have $K_* \leq c^*/c_* - 1/2$ in (4.5) with $\alpha = 1/2$. Let $m = m^* \equiv d^* - d^o \geq (c^*/c_* - 1/2)d^o \geq d^o/2$. As in the proof of Theorem 6, for the $\lambda$ in part (i), Lemma 2 gives $P\{2\sqrt{c_*}\zeta(\mathbf{y}; m, B) > \lambda\} \leq \epsilon/\sqrt{4\log \widetilde{p_\epsilon}}$. Thus, (2.16) follows from (A.6). It remains to prove (2.17) with $\widetilde{p_1}$ in (2.16).

We first bound $\widetilde{p_1}$. Since $m! \geq (m/e)^m$ and $m \geq d^o = R^r/\lambda_{\mathrm{mm}}^r$, by (2.15)

$$\frac{2}{m}\log\binom{p}{m} \leq 2\log\left(\frac{ep}{m}\right) \leq 2\log\left(\frac{2ep\lambda_{\mathrm{mm}}^r}{R^r}\right) = n\lambda_{\mathrm{mm}}^2/\sigma^2 + r\log\left(\frac{n\lambda_{\mathrm{mm}}^2}{\sigma^2(2e)^{-2/r}}\right).$$

Thus, by (2.12), $\sigma\sqrt{(2/n)\log\widetilde{p_1}} \leq \lambda_{\mathrm{mm}} + \epsilon_1\sigma/\sqrt{n}$ for large $n\lambda_{\mathrm{mm}}^2/\sigma^2$.

Let $\boldsymbol{\beta} \in \widetilde{\Theta}_{r,R}$ and $B_k$ be the set of $j$ for the $d^o$ largest $|\beta_j|$ with $j \notin B_0 \cup \cdots \cup B_{k-1}$, $k \geq 1$, with $B_0 = \varnothing$. Let $B = B_1$ and $v_j = |\beta_j| \wedge \lambda_{\mathrm{mm}}$. Since $|\beta_j| \leq \|\mathbf{v}_{B_{k-1}}\|_1/d^o$ for $j \in B_k$ and $k \geq 2$, $\sum_{k\geq 2}\|\boldsymbol{\beta}_{B_k}\|/\sqrt{d^o} \leq \sum_{k\geq 2}\|\mathbf{v}_{B_{k-1}}\|_1/d^o = \|\mathbf{v}\|_1/d^o \leq R^r\lambda_{\mathrm{mm}}^{1-r}/d^o = \lambda_{\mathrm{mm}}$. Thus, $\theta_B = \|(\mathbf{I}_n - \mathbf{P}_B)\mathbf{X}\boldsymbol{\beta}\|/\sqrt{n} \leq \sum_{k\geq 2}\|\mathbf{X}_{B_k}\boldsymbol{\beta}_{B_k}\|/\sqrt{n} \leq \sqrt{c^*d^o}\lambda_{\mathrm{mm}}$ by (2.11). Since $c^*d^o \leq 2c_*m$, $\theta_B/\sqrt{m} + \sigma\sqrt{(2/n)\log\widetilde{p_1}} \leq (\sqrt{2c_*} + 1)\lambda_{\mathrm{mm}} + \epsilon_1\sigma/\sqrt{n} = \lambda/(2\sqrt{c^*})$, so that

$$(A.9) \qquad \sup_{\boldsymbol{\beta} \in \widetilde{\Theta}_{r,R}} P\{c_*\|\widehat{\boldsymbol{\beta}}(\lambda) - \widehat{\boldsymbol{\beta}}^o\| \geq (3/2)\lambda\sqrt{d^o}\} \to 0$$

by (2.16). Since $\lambda = 2\sqrt{c^*}\lambda_{\mathrm{mm}}(1 + \sqrt{2c_*} + o(1))$, by the Hölder inequality, (A.9) and (A.5) imply that $\|\widehat{\boldsymbol{\beta}}(\lambda) - \widehat{\boldsymbol{\beta}}^o\|_q^q \leq |A_1|^{1-q/2}\|\widehat{\boldsymbol{\beta}}(\lambda) - \widehat{\boldsymbol{\beta}}^o\|^q \leq M_{1,q}^q\lambda_{\mathrm{mm}}^q d^o$ with large probability. Moreover, we have $\|\widehat{\boldsymbol{\beta}}^o - E\widehat{\boldsymbol{\beta}}^o\|^2 \leq O_P(1)d^o\sigma^2/(nc_*) = o_P(\lambda_{\mathrm{mm}}^2 d^o/c_*)$ and $\|E\widehat{\boldsymbol{\beta}}_B^o - \boldsymbol{\beta}_B\|^2 = \|\boldsymbol{\Sigma}_B^{-1}\boldsymbol{\Sigma}_{B,B^c}\boldsymbol{\beta}_{B^c}\|^2 \leq (c^*/c_*)(\sum_{k\geq 2}\|\boldsymbol{\beta}_{B_k}\|)^2 \leq (c^*/c_*)\lambda_{\mathrm{mm}}^2 d^o$. Since $\|\boldsymbol{\beta}_{B^c}\|_q^q \leq R^r\|\boldsymbol{\beta}_{B^c}\|_\infty^{q-r} \leq \lambda_{\mathrm{mm}}^q d^o$, these inequalities imply that $\|\widehat{\boldsymbol{\beta}}^o - \boldsymbol{\beta}\|_q^q = \|\widehat{\boldsymbol{\beta}}_B^o - \boldsymbol{\beta}_B\|_q^q + \|\boldsymbol{\beta}_{B^c}\|_q^q \leq |B|^{1-q/2}\{(o_P(1/c_*) + c^*/c_*)\lambda_{\mathrm{mm}}^2 d^o\}^{q/2} + \lambda_{\mathrm{mm}}^q d^o \leq M_{2,q}^q\lambda_{\mathrm{mm}}^q d^o$ with large probability. We obtain (2.17) by combining the upper bounds for $\|\widehat{\boldsymbol{\beta}}(\lambda) - \widehat{\boldsymbol{\beta}}^o\|_q^q$ and $\|\widehat{\boldsymbol{\beta}}^o - \boldsymbol{\beta}\|_q^q$. $\square$

PROOF OF THEOREM 9. (ii) $\Rightarrow$ (iii): let $\lambda$ be fixed and $\lambda_0 \equiv \dot{\rho}(0+; \lambda)$. Define $h(t) \equiv \kappa(\rho; \lambda)t^2/2 + \rho(|t|; \lambda) - \lambda_0|t|$. Since $\kappa(\rho; \lambda)$ is the maximum concavity in (2.2), $h(|t|)$ is a continuously differentiable convex function in $\mathbb{R}$. It follows that the penalized loss

$$L(\mathbf{b}; \lambda) = \left\{ \frac{1}{2n} \|\mathbf{y} - \mathbf{X}\mathbf{b}\|^2 - \frac{\kappa(\rho; \lambda)}{2} \|\mathbf{b}\|^2 \right\} + \sum_{j=1}^{p} \{\lambda_0 |b_j| + h(|b_j|)\}$$

is a sum of two convex functions, with the first one being strictly convex for $c_{\min}(\boldsymbol{\Sigma}) > \kappa(\rho; \lambda)$ and the second one being strictly convex otherwise.

(iii) $\Rightarrow$ (i): since the penalized loss $L(\mathbf{b}; \lambda)$ is $\|\mathbf{y}\|^2/(2n)$ for $\mathbf{b} = \mathbf{0}$, $\mathbf{y} \to \widehat{\boldsymbol{\beta}}$ maps bounded sets of $\mathbf{y}$ in $\mathbb{R}^n$ to bounded sets of $\widehat{\boldsymbol{\beta}}$ in $\mathbb{R}^p$. Since $L(\mathbf{b}; \lambda)$ is continuous in both $\mathbf{y}$ and $\mathbf{b}$ and strictly convex in $\mathbf{b}$ for each $\mathbf{y}$, its global minimum is unique and continuous in $\mathbf{y}$.

(i) $\Rightarrow$ (ii): since $\widehat{\boldsymbol{\beta}}$ depends on $\mathbf{y}$ only through $\widetilde{\mathbf{z}} = \mathbf{X}'\mathbf{y}/n$ and $\mathbf{X}$ is of rank $p$, the map $\widetilde{\mathbf{z}} \to \widehat{\boldsymbol{\beta}}$ is continuous from $\mathbb{R}^p$ to its range $\mathscr{I}$. Since $\widehat{\boldsymbol{\beta}}$ is the global minimum, (2.6) must hold and the inverse $\widehat{\boldsymbol{\beta}} \to \widetilde{\mathbf{z}} = \boldsymbol{\Sigma}\widehat{\boldsymbol{\beta}} + \operatorname{sgn}(\widehat{\boldsymbol{\beta}})\dot{\rho}(|\widehat{\boldsymbol{\beta}}|; \lambda)$ is continuous for $\widehat{\boldsymbol{\beta}} \in (0, \infty)^p \cap \mathscr{I}$, with per component application of functions and the product operation. It follows that $(0, \infty)^p \cap \mathscr{I}$ is open and does not have a boundary point in $(0, \infty)^p$. Let $\mathbf{1} \equiv (1, \ldots, 1)' \in \mathbb{R}^p$. For $\widetilde{\mathbf{z}} = x\boldsymbol{\Sigma}\mathbf{1}$ with $x > 0$, $L(x\mathbf{1}; \lambda) = o(x^2)$ for the ordinary LSE $x\mathbf{1}$ by the first condition of (5.15), and $L(\mathbf{b}; \lambda)$ is at least $c_{\min}(\boldsymbol{\Sigma})x^2$ for any $\mathbf{b}$ outside $(0, \infty)^p$. Thus, $(0, \infty)^p \cap \mathscr{I}$ is not empty. As the only nonempty set without any boundary point in $(0, \infty)^p$, $(0, \infty)^p \cap \mathscr{I} = (0, \infty)^p$. Moreover, the map $\widetilde{\mathbf{z}} \to \widehat{\boldsymbol{\beta}}$ is one-to-one for $\widehat{\boldsymbol{\beta}} \in (0, \infty)^p$.

We have proved that all points $\boldsymbol{\beta}$ in $(0, \infty)^p$ are unique global minimum of (1.1) for some $\mathbf{z} \in \mathbb{R}^p$. Let $\widehat{\boldsymbol{\beta}} = x\mathbf{1} \in (0, \infty)^p$ and $\mathbf{b}$ be the eigenvector with $\boldsymbol{\Sigma}\mathbf{b} = c_{\min}(\boldsymbol{\Sigma})\mathbf{b}$ and $\|\mathbf{b}\| = 1$. The quantity

$$
\begin{aligned}
(A.10) \quad & t^{-1} \frac{\partial}{\partial t} L(\widehat{\boldsymbol{\beta}} + t\mathbf{b}; \lambda) \\
& = \|\mathbf{X}\mathbf{b}\|^2 + \sum_{j=1}^{p} t^{-1} \operatorname{sgn}(\widehat{\beta}_j) b_j \{\dot{\rho}(|\widehat{\beta}_j + tb_j|; \lambda) - \dot{\rho}(|\widehat{\beta}_j|; \lambda)\} \\
& = c_{\min}(\boldsymbol{\Sigma}) + \sum_{j=1}^{p} t^{-1} b_j \{\dot{\rho}(x + tb_j; \lambda) - \dot{\rho}(x; \lambda)\}
\end{aligned}
$$

must have nonnegative lower limit as $t \to 0+$. Integrating over $x \in [t_1, t_2]$ and then taking the limit, we find

$$
\begin{aligned}
(A.11) \quad & c_{\min}(\boldsymbol{\Sigma})(t_2 - t_1) + \dot{\rho}(t_2; \lambda) - \dot{\rho}(t_1; \lambda) \\
& = \lim_{t \to 0+} \int_{t_1}^{t_2} t^{-1} \frac{\partial}{\partial t} L(x\mathbf{1} + t\mathbf{b}; \lambda) \, dx \geq 0.
\end{aligned}
$$

It remains to prove that (A.11) holds with strict inequality. If (A.11) holds with equality for certain $0 < t_1 < t_2$, then for $t_1 < x < t_2$ and small $t$ (A.10) becomes

$$t^{-1} \frac{\partial}{\partial t} L(\widehat{\boldsymbol{\beta}} + t\mathbf{b}; \lambda) = c_{\min}(\boldsymbol{\Sigma}) + \sum_{j=1}^{p} t^{-1} b_j \{-c_{\min}(\boldsymbol{\Sigma}) t b_j\} = 0.$$

This is contradictory to the uniqueness of $\widehat{\boldsymbol{\beta}}$. $\square$

PROOF OF PROPOSITION 2. Let $\widehat{\mathbf{P}}$ be as in Theorem 7. We write (2.6) as

(A.12)
$$\begin{cases} \widehat{\mathbf{P}}\boldsymbol{\Sigma}\widehat{\boldsymbol{\beta}} + \widehat{\mathbf{P}}\operatorname{sgn}(\widehat{\boldsymbol{\beta}})\dot{\rho}(|\widehat{\boldsymbol{\beta}}|; \lambda) = \widehat{\mathbf{P}}\widetilde{\mathbf{z}}, \\ |\widetilde{z}_j - \mathbf{x}_j'\mathbf{X}\widehat{\boldsymbol{\beta}}/n| \leq \lambda, \forall j. \end{cases}$$

Let $\boldsymbol{\eta} \in \{-1, 0, 1\}^p$ be fixed (not confused with the $\boldsymbol{\eta}$ in Section 3). It follows from Theorem 9 that the map $\widehat{\mathbf{P}}\widetilde{\mathbf{z}} \to \widehat{\mathbf{P}}\widehat{\boldsymbol{\beta}}$ is continuous in $\widetilde{\mathbf{z}} \in \mathbb{R}^p$ and continuously invertible given a fixed $\operatorname{sgn}(\widehat{\boldsymbol{\beta}}) = \boldsymbol{\eta}$. Let $H(\boldsymbol{\eta}) \equiv \{\widetilde{\mathbf{z}} : \operatorname{sgn}(\widehat{\boldsymbol{\beta}}) = \boldsymbol{\eta}\}$. The boundary of $H(\boldsymbol{\eta})$ has zero Lebesgue measure, since it is contained in the set of $\widetilde{\mathbf{z}}$ satisfying $\eta_j \widehat{\beta}_j = 0+$ for $\eta_j \neq 0$ or $\widetilde{z}_j - \mathbf{x}_j'\mathbf{X}\widehat{\boldsymbol{\beta}}/n = \pm\lambda$ for $\eta_j = 0$, $j = 1, \ldots, p$, according to (A.12). In the interior of $H(\boldsymbol{\eta})$, (A.12) gives $(\partial/\partial\widetilde{z}_j)\widehat{\boldsymbol{\beta}} = \mathbf{0}$ and $(\partial/\partial\widetilde{\mathbf{z}})\widehat{\beta}_j = 0$ for $\eta_j = 0$ and

$$\widehat{\mathbf{P}}\frac{\partial}{\partial\widehat{\boldsymbol{\beta}}}(\widehat{\mathbf{P}}\widetilde{\mathbf{z}})' = \widehat{\mathbf{P}}\boldsymbol{\Sigma}\widehat{\mathbf{P}} + \widehat{\mathbf{P}}\operatorname{diag}(\ddot{\rho}(|\widehat{\beta}_j|; \lambda))\widehat{\mathbf{P}}' = \mathbf{Q}(\widehat{\boldsymbol{\beta}}; \lambda).$$

Since (2.5) holds with $d^* = p$, $c_{\min}(\mathbf{Q}(\boldsymbol{\beta}; \lambda)) \geq c_{\min}(\boldsymbol{\Sigma}) - \kappa(\rho; \lambda) > 0$ for all $\boldsymbol{\beta} \neq \mathbf{0}$. Thus, the differentiation of the inverse map yields $(\partial/\partial\widetilde{\mathbf{z}})\widehat{\boldsymbol{\beta}}' = \widehat{\mathbf{P}}'\mathbf{Q}^{-1}(\widehat{\boldsymbol{\beta}}; \lambda)\widehat{\mathbf{P}}$. $\square$

PROOF OF THEOREM 7. It follows from Proposition 2 that $\widehat{\boldsymbol{\beta}} - \boldsymbol{\Sigma}^{-1}\widetilde{\mathbf{z}}$ is almost differentiable in $\widetilde{\mathbf{z}}$ with derivative

$$\frac{\partial}{\partial\widetilde{\mathbf{z}}}(\widehat{\boldsymbol{\beta}} - \boldsymbol{\Sigma}^{-1}\widetilde{\mathbf{z}})' = \widehat{\mathbf{P}}'\mathbf{Q}^{-1}(\widehat{\boldsymbol{\beta}}; \lambda)\widehat{\mathbf{P}} - \boldsymbol{\Sigma}^{-1}.$$

Since $\widetilde{\mathbf{z}} \equiv \mathbf{X}'\mathbf{y}/n \sim N(\boldsymbol{\Sigma}\boldsymbol{\beta}, \boldsymbol{\Sigma}\sigma^2/n)$, this and (5.3) imply

$$E(\widehat{\boldsymbol{\beta}} - \boldsymbol{\Sigma}^{-1}\widetilde{\mathbf{z}})(\boldsymbol{\Sigma}^{-1}\widetilde{\mathbf{z}} - \boldsymbol{\beta})' = E(\widehat{\boldsymbol{\beta}} - \boldsymbol{\Sigma}^{-1}\widetilde{\mathbf{z}})(\widetilde{\mathbf{z}} - \boldsymbol{\Sigma}\boldsymbol{\beta})'\boldsymbol{\Sigma}^{-1}$$

$$= \frac{\sigma^2}{n}\{E\widehat{\mathbf{P}}'\mathbf{Q}^{-1}(\widehat{\boldsymbol{\beta}}; \lambda)\widehat{\mathbf{P}} - \boldsymbol{\Sigma}^{-1}\}.$$

Since the ordinary LSE is $\widetilde{\boldsymbol{\beta}} = \boldsymbol{\Sigma}^{-1}\widetilde{\mathbf{z}} \sim N(\boldsymbol{\beta}, \boldsymbol{\Sigma}^{-1}\sigma^2/n)$, it follows that

$$E(\widehat{\boldsymbol{\beta}} - \boldsymbol{\beta})(\widehat{\boldsymbol{\beta}} - \boldsymbol{\beta})'$$

$$= E(\widehat{\boldsymbol{\beta}} - \widetilde{\boldsymbol{\beta}})(\widehat{\boldsymbol{\beta}} - \widetilde{\boldsymbol{\beta}})' - E(\boldsymbol{\beta} - \widetilde{\boldsymbol{\beta}})(\boldsymbol{\beta} - \widetilde{\boldsymbol{\beta}})' + 2E(\widehat{\boldsymbol{\beta}} - \widetilde{\boldsymbol{\beta}})(\widetilde{\boldsymbol{\beta}} - \boldsymbol{\beta})'$$

$$= E(\widehat{\boldsymbol{\beta}} - \widetilde{\boldsymbol{\beta}})(\widehat{\boldsymbol{\beta}} - \widetilde{\boldsymbol{\beta}})' + \frac{2\sigma^2}{n}\{E\widehat{\mathbf{P}}'\mathbf{Q}^{-1}(\widehat{\boldsymbol{\beta}}; \lambda)\widehat{\mathbf{P}} - \boldsymbol{\Sigma}^{-1}\} + \frac{\sigma^2}{n}\boldsymbol{\Sigma}^{-1}.$$

This proves (5.5). The rest of the theorem follows immediately. $\square$

PROOF OF THEOREM 8. Since trace($\mathbf{bb}'$) = $\|\mathbf{b}\|^2$, (5.5) gives

$$E\|\widehat{\boldsymbol{\mu}} - \boldsymbol{\mu}\|^2 = E\left\{\|\widehat{\boldsymbol{\mu}} - \widetilde{\boldsymbol{\mu}}\|^2 + \frac{2\sigma^2}{n}\operatorname{trace}(\mathbf{X}\widehat{\mathbf{P}}'\mathbf{Q}^{-1}(\widehat{\boldsymbol{\beta}}; \lambda)\widehat{\mathbf{P}}\mathbf{X}')\right\}$$

$$- \frac{\sigma^2}{n}\operatorname{trace}(\mathbf{X}\boldsymbol{\Sigma}^{-1}\mathbf{X}')$$

$$= E\{\|\widehat{\boldsymbol{\mu}} - \widetilde{\boldsymbol{\mu}}\|^2 + 2\sigma^2\widehat{\mathrm{df}} - \sigma^2\operatorname{rank}(\mathbf{X})\},$$

which implies (5.11) via (5.8). For (5.12), we observe that $\mathbf{Q}(\widehat{\boldsymbol{\beta}}; \lambda) = \widehat{\mathbf{P}}\boldsymbol{\Sigma}\widehat{\mathbf{P}}'$ by (5.4) when $\ddot{\rho}(|\widehat{\beta}_j|; \lambda) = 0$ for all $\widehat{\beta}_j \neq 0$. $\square$

PROOF OF PROPOSITION 3. Let $\mathbf{u}_1, \ldots, \mathbf{u}_N$ be vectors in the unit sphere $S^{m-1}$ of the range of $\mathbf{P}$ such that balls $\{\mathbf{v}: \|\mathbf{v} - \mathbf{u}_j\| \leq \epsilon_1\}$ are disjoint and $\bigcup_{j=1}^N \{\mathbf{v}: \|\mathbf{v} - \mathbf{u}_j\| \leq 2\epsilon_1\} \supset S^{m-1}$. Volume comparison yields $N\epsilon_1^m \leq (1 + \epsilon_1)^m - (1 - \epsilon_1)^m$. Since $\mathbf{v}'\boldsymbol{\varepsilon} = \mathbf{u}'\boldsymbol{\varepsilon} + (\mathbf{v} - \mathbf{u})'\boldsymbol{\varepsilon}$, $\|\mathbf{P}\boldsymbol{\varepsilon}\| = \max_{\mathbf{v}\in S^{m-1}} \mathbf{v}'\boldsymbol{\varepsilon} \leq \max_{j\leq N} \mathbf{u}_j'\boldsymbol{\varepsilon} + 2\epsilon_1\|\mathbf{P}\boldsymbol{\varepsilon}\| \leq \max_{j\leq N} \mathbf{u}_j'\boldsymbol{\varepsilon}/(1 - 2\epsilon_1)_+$. It follows that $P\{\|\mathbf{P}\boldsymbol{\varepsilon}\| > \sigma_1 t\} \leq (1 + 1/\epsilon_1)^m e^{-(1-2\epsilon_1)^2 t^2/2}$. Taking $t^2 = m(1+x)/(1 - 2\epsilon_1)^2$, we find

$$P\left\{\|\mathbf{P}\boldsymbol{\varepsilon}\|^2/\sigma_1^2 \geq \frac{m(1+x)}{(1 - 2\epsilon_1)_+^2}\right\} \leq (1 + 1/\epsilon_1)^m e^{-m(1+x)/2} \leq e^{-mx/2}(1+x)^{m/2}$$

for $(1 + 1/\epsilon_1)^2 = (1+x)e^x$. This proves the proposition since $\epsilon_1 = 1/(e^{x/2} \times \sqrt{1+x} - 1)$. $\square$

PROOF OF LEMMA 1. Let $\mathbf{X}_1 \equiv \mathbf{X}_{A_1}$ as in (2.4) and $\boldsymbol{\Sigma}_{11} \equiv \mathbf{X}_1'\mathbf{X}_1/n$. Since $|A_1| \leq d^*$,

$$(\mathrm{A}.13) \quad c_* \leq \frac{\|\boldsymbol{\Sigma}_{11}\mathbf{v}\|^2}{\|\mathbf{v}\|} \leq c^*, \qquad \frac{1}{c^*} \leq \frac{\|\boldsymbol{\Sigma}_{11}^{-1}\mathbf{v}\|^2}{\|\mathbf{v}\|} \leq \frac{1}{c_*} \qquad \forall 0 \neq \mathbf{v} \in \mathbb{R}^{|A_1|},$$

by (2.11). Set $A_2 \equiv \{1, \ldots, p\} \setminus A_1$, $A_3 \equiv B$ and $A_4 \equiv A_1 \setminus B$. Define $\mathbf{b}_k \equiv (b_j, j \in A_k)$ for $\mathbf{b} \in \mathbb{R}^p$ and $k = 1, 2, 3, 4$. For $k = 3, 4$, let $\mathbf{Q}_k$ be the matrix representing the selection of variables in $A_k$ from $A_1$, defined as $\mathbf{Q}_k\mathbf{b}_1 = \mathbf{b}_k$.

Let $\widehat{\boldsymbol{\beta}}^o$ be the oracle LSE in (2.10) and $\widetilde{\boldsymbol{\varepsilon}} \equiv \mathbf{y} - \mathbf{X}\widehat{\boldsymbol{\beta}}^o = (\mathbf{I}_n - \mathbf{P}_B)\mathbf{y}$. Since $\widehat{\boldsymbol{\beta}}_2 = \widehat{\boldsymbol{\beta}}_2^o = \mathbf{0}$, the $A_1$ components of the negative gradient

$$(\mathrm{A}.14) \qquad\qquad \mathbf{g} \equiv \mathbf{X}'(\mathbf{y} - \mathbf{X}\widehat{\boldsymbol{\beta}})/n$$

must satisfy $\mathbf{g}_1 = \mathbf{X}_1'(\mathbf{y} - \mathbf{X}_1\widehat{\boldsymbol{\beta}}_1)/n = \mathbf{X}_1'\widetilde{\boldsymbol{\varepsilon}}/n + \boldsymbol{\Sigma}_{11}(\widehat{\boldsymbol{\beta}}_1^o - \widehat{\boldsymbol{\beta}}_1)$, so that

$$(A.15) \qquad \boldsymbol{\Sigma}_{11}^{-1}\mathbf{g}_1 + (\widehat{\boldsymbol{\beta}}_1 - \widehat{\boldsymbol{\beta}}_1^o) = \boldsymbol{\Sigma}_{11}^{-1}\mathbf{X}_1'\widetilde{\boldsymbol{\varepsilon}}/n.$$

Let $\mathbf{v}_1 \equiv \boldsymbol{\Sigma}_{11}^{-1/2}\mathbf{g}_1$ and $\mathbf{v}_k \equiv \boldsymbol{\Sigma}_{11}^{-1/2}\mathbf{Q}_k'\mathbf{g}_k$, $k = 3, 4$. Let $\mathbf{P}_1 \equiv \mathbf{X}_1\boldsymbol{\Sigma}_{11}^{-1}\mathbf{X}_1'/n = \mathbf{P}_{A_1}$ be the projection to the range of $\mathbf{X}_1$ as in (A.4). Since $A_1 \supset B$, $\mathbf{P}_1\widetilde{\boldsymbol{\varepsilon}} = (\mathbf{P}_1 - \mathbf{P}_B)\mathbf{y}$, so that $\|\mathbf{P}_1\widetilde{\boldsymbol{\varepsilon}}\|^2/n \le |A_4|\zeta^2(\mathbf{y}; |A_4|, B)$ by (A.4). Thus, for $\lambda \ge (\sqrt{c^*}/\alpha)\zeta(\mathbf{y}; |A_4|, B)$ as provided,

$$(A.16) \qquad \mathbf{g}_k'\mathbf{Q}_k\boldsymbol{\Sigma}_{11}^{-1}\mathbf{X}_1'\widetilde{\boldsymbol{\varepsilon}}/n \le \|\mathbf{v}_k\|\|\mathbf{P}_1\widetilde{\boldsymbol{\varepsilon}}\|/\sqrt{n} \le \|\mathbf{v}_k\|\alpha\lambda\sqrt{|A_4|/c^*}.$$

Since $\mathbf{Q}_4(\widehat{\boldsymbol{\beta}}_1^o - \widehat{\boldsymbol{\beta}}_1) = \widehat{\boldsymbol{\beta}}_4^o - \widehat{\boldsymbol{\beta}}_4 = -\widehat{\boldsymbol{\beta}}_4$ and $\mathbf{v}_3 = \mathbf{v}_1 - \mathbf{v}_4$, by (A.15) we have

$$\|\mathbf{v}_4\|^2 - \|\mathbf{v}_3\|^2 + \|\mathbf{v}_1\|^2 = 2\mathbf{v}_4'\mathbf{v}_1 = 2\mathbf{g}_4'\mathbf{Q}_4\boldsymbol{\Sigma}_{11}^{-1}\mathbf{g}_1 = 2\mathbf{g}_4'\mathbf{Q}_4\boldsymbol{\Sigma}_{11}^{-1}\mathbf{X}_1'\widetilde{\boldsymbol{\varepsilon}}/n - 2\mathbf{g}_4'\widehat{\boldsymbol{\beta}}_4.$$

Since $2\|\mathbf{v}_4\|\lambda\sqrt{|A_4|/c^*} \le \|\mathbf{v}_4\|^2 + \lambda^2|A_4|/c^*$, the above identity and (A.16) yield

$$(1-\alpha)\|\mathbf{v}_4\|^2 + \|\mathbf{v}_1\|^2 + 2\mathbf{g}_4'\widehat{\boldsymbol{\beta}}_4 \le \|\mathbf{v}_3\|^2 + \alpha\lambda^2|A_4|/c^*.$$

Similarly, it follows from (A.15) and (A.16) that

$$\begin{aligned}
\|\mathbf{v}_4\|^2 &+ 2\mathbf{g}_4'\widehat{\boldsymbol{\beta}}_4 + \|\boldsymbol{\Sigma}_{11}^{1/2}(\widehat{\boldsymbol{\beta}}_1 - \widehat{\boldsymbol{\beta}}_1^o)\|^2 \\
&= \|\mathbf{v}_4\|^2 + 2\mathbf{g}_4'\widehat{\boldsymbol{\beta}}_4 + \|\mathbf{v}_1\|^2 - 2\mathbf{g}_1'\boldsymbol{\Sigma}_{11}^{-1}\mathbf{X}_1'\widetilde{\boldsymbol{\varepsilon}}/n + \|\mathbf{P}_1\widetilde{\boldsymbol{\varepsilon}}\|^2/n \\
&= \|\mathbf{v}_3\|^2 + 2\mathbf{g}_4'\mathbf{Q}_4\boldsymbol{\Sigma}_{11}^{-1}\mathbf{X}_1'\widetilde{\boldsymbol{\varepsilon}}/n - 2\mathbf{g}_1'\boldsymbol{\Sigma}_{11}^{-1}\mathbf{X}_1'\widetilde{\boldsymbol{\varepsilon}}/n + \|\mathbf{P}_1\widetilde{\boldsymbol{\varepsilon}}\|^2/n \\
&= \|\mathbf{v}_3\|^2 - 2\mathbf{g}_3'\mathbf{Q}_3\boldsymbol{\Sigma}_{11}^{-1}\mathbf{X}_1'\widetilde{\boldsymbol{\varepsilon}}/n + \|\mathbf{P}_1\widetilde{\boldsymbol{\varepsilon}}\|^2/n \\
&\le \|\mathbf{v}_3\|^2 + 2\|\mathbf{v}_3\|\alpha\lambda\sqrt{|A_4|/c^*} + \alpha^2\lambda^2|A_4|/c^*
\end{aligned}$$

due to $\mathbf{g}_1' = \mathbf{g}_3'\mathbf{Q}_3 + \mathbf{g}_4'\mathbf{Q}_4$. For the $w \equiv (2-\alpha)/(c_*c^*/\kappa^2 - 1)$ in (4.5), the $\{1, w\}$ weighted sum of the above two inequalities yields

$$(A.17) \qquad \begin{aligned}
\text{LHS} &\equiv (1-\alpha+w)\|\mathbf{v}_4\|_2^2 + \|\mathbf{v}_1\|^2 + (1+w)2\mathbf{g}_4'\widehat{\boldsymbol{\beta}}_4 \\
&\quad + w\|\boldsymbol{\Sigma}_{11}^{1/2}(\widehat{\boldsymbol{\beta}}_1 - \widehat{\boldsymbol{\beta}}_1^o)\|^2 \\
&\le (1+w)\|\mathbf{v}_3\|^2 + (\alpha+w\alpha^2)\lambda^2|A_4|/c^* \\
&\quad + 2w\|\mathbf{v}_3\|\alpha\lambda\sqrt{|A_4|/c^*}.
\end{aligned}$$

Note that (A.17) holds with equality only in the following scenario: $\|\mathbf{v}_4\|^2 = \lambda^2|A_4|/c^*$ and (A.16) holds with equalities for both $\mathbf{v}_3$ and $\mathbf{v}_4$. Since $|A_4| = |A_1| - |B| > 0$ and $\boldsymbol{\Sigma}_{11}^{1/2}\mathbf{v}_k = \mathbf{Q}_k\mathbf{g}_k$ have different support for $k \in \{3, 4\}$, this scenario could happen only if $\|\mathbf{v}_3\| = 0$. Thus, (A.17) holds strictly unless $\|\mathbf{v}_3\| = \|\mathbf{g}_3\| = 0$.

We first bound the LHS. Since $\lambda(1 - \kappa|t|/\lambda)_+ \leq \dot{\rho}(|t|;\lambda) \leq \lambda$, (2.6) and (A.14) provide

$$(A.18) \quad \frac{\|\mathbf{g}_4\|^2}{\lambda^2} = \sum_{j \in A_4} \frac{\dot{\rho}^2(|\widehat{\beta}_j|;\lambda)}{\lambda^2} \geq \sum_{j \in A_4} \left(1 - \kappa\frac{|\widehat{\beta}_j|}{\lambda}\right)_+^2,$$

$$\frac{\|\mathbf{g}_3\|^2}{\lambda^2} \leq \sum_{j \in B} \frac{\dot{\rho}^2(|\widehat{\beta}_j|;\lambda)}{\lambda^2} \leq |B|,$$

in view of the second condition on $A_1 = A_4 \cup B$. We also have $\widehat{\beta}_j^o = 0$ and $\widehat{\beta}_j g_j = |\widehat{\beta}_j g_j|$ for $j \in A_4$. Thus, by (A.13) and the definition $\mathbf{v}_k \equiv \mathbf{\Sigma}_{11}^{-1/2}\mathbf{Q}_k'\mathbf{g}_k$,

$$
\begin{aligned}
\text{LHS} &\equiv (1 - \alpha + w)\|\mathbf{v}_4\|_2^2 + \|\mathbf{v}_1\|^2 + (1+w)2\mathbf{g}_4'\widehat{\boldsymbol{\beta}}_4 \\
&\quad + w\|\mathbf{\Sigma}_{11}^{1/2}(\widehat{\boldsymbol{\beta}}_1 - \widehat{\boldsymbol{\beta}}_1^o)\|^2 \\
&\geq (1 - \alpha + w)\|\mathbf{g}_4\|_2^2/c^* + \|\mathbf{g}_1\|^2/c^* \\
&\quad + (1+w)2\mathbf{g}_4'\widehat{\boldsymbol{\beta}}_4 + wc_*\|\widehat{\boldsymbol{\beta}}_1 - \widehat{\boldsymbol{\beta}}_1^o\|^2 \\
(A.19) \quad &\geq \lambda^2 \sum_{j \in A_4}\{(2 - \alpha + w)(1 - \kappa t_j)_+^2/c^* \\
&\qquad + (1+w)2(1 - \kappa t_j)_+ t_j + wc_* t_j^2\} + \frac{\|\mathbf{g}_3\|^2}{c^*} \\
&\geq \lambda^2|A_4| \min_{0 \leq \kappa t \leq 1}\{(2 - \alpha + w)(1 - \kappa t)^2/c^* \\
&\qquad + (1+w)2t(1 - \kappa t) + wc_* t^2\} + \frac{\|\mathbf{g}_3\|^2}{c^*},
\end{aligned}
$$

where $t_j \equiv |\widehat{\beta}_j|/\lambda$. Since $c^* \geq c_* \geq \kappa$ and $w \equiv (2 - \alpha)/(c_* c^*/\kappa^2 - 1)$, we have

$$
\begin{aligned}
&(2 - \alpha + w)\kappa^2/c^* - (1+w)2\kappa + wc_* \\
&= 2\{wc_* - \kappa(1 + w)\} \\
&= 2\frac{(2 - \alpha)c_* - \kappa(c_* c^*/\kappa^2 + 1 - \alpha)}{c_* c^*/\kappa^2 - 1} \leq 0
\end{aligned}
$$

due to $\kappa\alpha - c_*\alpha \leq 0$ and $-c_* c^*/\kappa^2 - 1 + 2c_*/\kappa \leq -(c_*/\kappa - 1)^2$. Thus, the minimum in (A.19) is taken over a concave quadratic function with equal value at $\{0, 1/\kappa\}$, so that

$$(A.20) \quad \text{LHS} \geq \lambda^2|A_4|(2 - \alpha + w)/c^* + \|\mathbf{g}_3\|^2/c^*.$$

Inserting (A.20) into (A.17), we find

$$\lambda^2|A_4|\{2 - \alpha + w - (\alpha + w\alpha^2)\}/c^*$$

$$\leq (1+w)\|\mathbf{v}_3\|^2 - \|\mathbf{g}_3\|^2/c^* + 2w\|\mathbf{v}_3\|\alpha\lambda\sqrt{|A_4|/c^*}$$

$$\leq (1+w)\|\mathbf{v}_3\|^2 - \|\mathbf{g}_3\|^2/c^* + w\alpha\left(\frac{\|\mathbf{v}_3\|^2}{t(1-\alpha)} + t(1-\alpha)\lambda^2|A_4|/c^*\right)$$

and that the strict inequality holds unless $\|\mathbf{v}_3\| = \|\mathbf{g}_3\| = 0$. We move $w\alpha t(1-\alpha)\lambda^2|A_4|/c^*$ to the left-hand side and then multiply both sides by $c^*/\lambda^2$ to arrive at

$$\begin{aligned}
&\{2 + w(1+\alpha) - tw\alpha\}(1-\alpha)|A_4| \\
(\text{A.21}) \quad &< (1+w\{1+(\alpha/t)/(1-\alpha)\})c^*\|\mathbf{v}_3\|^2/\lambda^2 - \|\mathbf{g}_3\|^2/\lambda^2 \\
&\leq \{(1+w\{1+(\alpha/t)/(1-\alpha)\})c^*/c_* - 1\}\|\mathbf{g}_3\|^2/\lambda^2
\end{aligned}$$

due to $c_*\|\mathbf{v}_3\|^2 \leq \|\mathbf{g}_3\|^2$ by (A.13). The strict inequality holds above, since the equality would imply $\|\mathbf{g}_3\| = 0$ and then $|A_4| = 0$. This proves (A.5) via (A.18).

For (A.6), it follows from (A.15), $(\widehat{\boldsymbol{\beta}}_4' - \widehat{\boldsymbol{\beta}}_4^o)'\mathbf{g}_4 \geq 0$ and then (A.13) and (A.16) that

$$\begin{aligned}
&(\widehat{\boldsymbol{\beta}}_1 - \widehat{\boldsymbol{\beta}}_1^o)'\boldsymbol{\Sigma}_{11}(\widehat{\boldsymbol{\beta}}_1 - \widehat{\boldsymbol{\beta}}_1^o) \\
&= -(\widehat{\boldsymbol{\beta}}_1 - \widehat{\boldsymbol{\beta}}_1^o)'\mathbf{g}_1 + (\widehat{\boldsymbol{\beta}}_1 - \widehat{\boldsymbol{\beta}}_1^o)'\mathbf{X}_1'\widetilde{\boldsymbol{\varepsilon}}/n \\
&\leq \|\widehat{\boldsymbol{\beta}}_3 - \widehat{\boldsymbol{\beta}}_3^o\|\|\mathbf{g}_3\| + \|\mathbf{X}_1(\widehat{\boldsymbol{\beta}}_1 - \widehat{\boldsymbol{\beta}}_1^o)\|\|\mathbf{P}_1\widetilde{\boldsymbol{\varepsilon}}\|/n \\
&\leq \|\boldsymbol{\Sigma}_{11}^{1/2}(\widehat{\boldsymbol{\beta}}_1 - \widehat{\boldsymbol{\beta}}_1^o)\|\|\mathbf{g}_3\|/\sqrt{c_*} + \|\boldsymbol{\Sigma}_{11}^{1/2}(\widehat{\boldsymbol{\beta}}_1 - \widehat{\boldsymbol{\beta}}_1^o)\|\alpha\lambda\sqrt{|A_4|/c^*}.
\end{aligned}$$

Dividing both sides by $\|\boldsymbol{\Sigma}_{11}^{1/2}(\widehat{\boldsymbol{\beta}}_1 - \widehat{\boldsymbol{\beta}}_1^o)\|$, we find with another application of (A.13) that

$$c_*\|\widehat{\boldsymbol{\beta}} - \widehat{\boldsymbol{\beta}}^o\| \leq \sqrt{c_*}\|\boldsymbol{\Sigma}_{11}^{1/2}(\widehat{\boldsymbol{\beta}}_1 - \widehat{\boldsymbol{\beta}}_1^o)\| \leq \|\mathbf{g}_3\| + \alpha\lambda\sqrt{|A_4|c_*/c^*}.$$

Since $\|\mathbf{X}(\widehat{\boldsymbol{\beta}} - \widehat{\boldsymbol{\beta}}^o)\|/\sqrt{n} = \|\boldsymbol{\Sigma}_{11}^{1/2}(\widehat{\boldsymbol{\beta}} - \widehat{\boldsymbol{\beta}}^o)\|$, this proves (A.6) via (A.18) and (A.5). $\quad\square$

PROOF OF LEMMA 2. Since $m$ and $B$ are deterministic, $nm\zeta^2(\boldsymbol{\varepsilon}; m, B)/\sigma^2$ in (A.4) is the maximum of $\binom{p-d^o}{m}$ variables with the $\chi_m^2$ distribution, so that

$$(\text{A.22}) \quad P\{\zeta(\boldsymbol{\varepsilon}; m, B) \geq \sigma\sqrt{(2/n)\log\widetilde{p}_\epsilon}\} \leq \binom{p-d^o}{m} P\{\chi_m^2 \geq m(1+x)\}$$

with $x = 2\log\widetilde{p}_\epsilon - 1 > 0$. Since $\chi_m^2/(1+x)$ has the gamma$(m/2, (1+x)/2)$ distribution,

$$\begin{aligned}
(\text{A.23}) \quad &P\{\chi_m^2 > m(1+x)\} \\
&= \frac{e^{-m(1+x)/2}(1+x)^{m/2}}{\Gamma(m/2)2^{m/2}}\int_m^\infty t^{m/2-1}e^{-(1+x)(t-m)/2}\,dt.
\end{aligned}$$

Let $y = \sqrt{t}$ and $h(y) = (1+x)(y^2-m)/2 - (m-1)\log y$. Since $(d/dy)^2 h(y) \geq (1+x)$,

$$
(A.24) \quad
\begin{aligned}
&\int_m^\infty t^{m/2-1} e^{-(1+x)(t-m)/2}\, dt \\
&\qquad = \int_{\sqrt{m}}^\infty 2 e^{-h(y)}\, dy \leq \frac{2 e^{-h(\sqrt{m})}}{\sqrt{1+x}} \int_0^\infty e^{-\mu z - z^2/2}\, dz
\end{aligned}
$$

with $z = \sqrt{1+x}(y - \sqrt{m})$ and $\mu = (dh/dy)(\sqrt{m})/\sqrt{1+x} = (x+1/m)\sqrt{m}/\sqrt{1+x}$. Since

$$
\frac{e^{-m/2} 2 e^{-h(\sqrt{m})}}{\Gamma(m/2) 2^{m/2}} \leq \frac{e^{-m/2} 2 m^{(m-1)/2}}{(m/2)^{m/2-1/2} e^{-m/2}\sqrt{2\pi}\, 2^{m/2}} = \frac{1}{\sqrt{\pi}}
$$

by the Stirling formula and $x+1 = 2\log \widetilde{p}_\epsilon$, (A.23) and (A.24) imply

$$
P\{\chi_m^2 \geq m(1+x)\} \leq \frac{e^{-mx/2}(1+x)^{m/2}}{\sqrt{2\pi \log \widetilde{p}_\epsilon}} \int_0^\infty e^{-\mu z - z^2/2}\, dz.
$$

This and (A.22) imply (A.7), since $(2\pi)^{-1/2} \int_0^\infty e^{-\mu z - z^2/2}\, dz = e^{\mu^2/2}\Phi(-\mu) \leq 1/2$ and (2.12) implies $\binom{p-d^o}{m} e^{-mx/2}(1+x)^{m/2} = \epsilon$.    $\square$

DEPARTMENT OF STATISTICS
  AND BIOSTATISTICS
BUSCH CAMPUS
RUTGERS UNIVERSITY
PISCATAWAY, NEW JERSEY 08854
USA
E-MAIL: czhang@stat.rutgers.edu